\documentclass{compositio}

\font\tensans=cmss10 

\font\teneuf=eufm10
\font\seveneuf=eufm7
\font\fiveeuf=eufm5
\font\tenmsa=msam10
\font\sevenmsa=msam7
\font\fivemsa=msam5
\font\tenmsb=msbm10
\font\sevenmsb=msbm7
\font\fivemsb=msbm5
\newfam\euffam
\textfont\euffam=\teneuf
\scriptfont\euffam=\seveneuf
\scriptscriptfont\euffam=\fiveeuf
\def\hexnumber@#1{\ifcase#1 0\or1\or2\or3\or4\or5\or6\or7\or8\or9\or
    A\or B\or C\or D\or E\or F\fi }
\newfam\msafam
\newfam\msbfam
\textfont\msafam=\tenmsa  \scriptfont\msafam=\sevenmsa
  \scriptscriptfont\msafam=\fivemsa
\textfont\msbfam=\tenmsb  \scriptfont\msbfam=\sevenmsb
  \scriptscriptfont\msbfam=\fivemsb
\edef\msb@{\hexnumber@\msbfam}
\mathchardef\smallsetminus="2\msb@72
\mathchardef\varGamma="0100
\mathchardef\Delta="0101
\mathchardef\varTheta="0102
\mathchardef\varLambda="0103
\mathchardef\varXi="0104
\mathchardef\varPi="0105
\mathchardef\varSigma="0106
\mathchardef\varUpsilon="0107
\mathchardef\varPhi="0108
\mathchardef\varPsi="0109
\mathchardef\varOmega="010A

\def\lj{\mathsf{[}}
\def\rj{\mathsf{]}}

\def\bbQ{\mathbf{Q}}
\def\bbR{\mathbf{R}}
\def\rto{\mathbf{R}\hskip-.5pt^2}
\def\rtr{\mathbf{R}\hskip-.7pt^3}

\def\bbC{\mathbf{C}}

\def\bbCP{\mathbf{C}\mathrm{P}}
\def\bbZ{\mathbf{Z}}
\def\dimc{\dim_{\hskip.4pt\bbC\hskip-1.2pt}}

\def\sgn{\mathrm{sgn}\hs}
\def\emp{{\tensans\O}}
\def\el{d}
\def\ka{n}
\def\lk{\ell}
\def\lkm{\lk_m}
\def\k{Y}

\def\y{c}

\def\cj{c}

\def\dz{Z}
\def\zi{R}

\def\rh{\theta}
\def\dl{\Delta}
\def\reg{\mathcal{U}}
\def\hut{\mathcal{H}\,\cup\,\mathcal{T}}
\def\xah{{\,\mathsf{X}\,\cap\,\mathcal{H}}}
\def\xat{{\,\mathsf{X}\,\cap\,\mathcal{T}}}
\def\smi{\mathcal{S}_m^{\hs\mathsf{I}}}
\def\smh{\mathcal{S}_m^{\mathcal{H}}}
\def\smt{\mathcal{S}_m^{\mathcal{T}}}
\def\xo{S}
\def\ax{A}
\def\bx{B}
\def\cx{C}
\def\fy{f}
\def\fj{\phi}
\def\vh{\varphi}

\def\hi{\chi}
\def\ee{\tilde E} 
\def\ef{\tilde F} 

\def\hz{l}

\def\hs{\hskip.7pt}
\def\nh{\hskip-.7pt}
\def\ns{\hskip-1.2pt}
\def\vu{\vbox{\hbox{\hskip.7pt$\uparrow$\hskip-.8pt}\vskip-2pt}}
\def\vum{\mathrm{\vbox{\hbox{\hskip.7pt$\uparrow$\hskip-.8pt}\vskip-2pt}}}
\def\vd{\vbox{\hbox{\hskip.7pt$\downarrow$\hskip-.8pt}\vskip-2pt}}

\def\hyp{\hskip.5pt\vbox
{\hbox{\vrule width2.5ptheight0.8ptdepth0pt}\vskip2pt}\hskip.5pt}
\def\di{\hskip1pt\vbox
{\hbox{\hskip.08pt\vrule width.92ptheight5.4ptdepth0pt\hskip.08pt}\vskip-11.4pt
\hbox{\vrule width1.1ptheight1.2ptdepth-.1pt}}\hskip1pt}

\def\nea{\vbox{\hbox{\hskip.7pt$^\nearrow$\hskip-.8pt}\vskip-2pt}}
\def\neam{\mathrm{\vbox{\hbox{\hskip.7pt$^\nearrow$\hskip-.8pt}\vskip-2pt
}}}
\def\sea{\vbox{\hbox{\hskip.7pt$^\searrow$\hskip-.8pt}\vskip-2pt}}
\def\seam{\mathrm{\vbox{\hbox{\hskip.7pt$^\searrow$\hskip-.8pt}\vskip-2pt
}}}

\def\dsq{d^{\hskip.4pt2}\hskip-1pt}

\def\th{u}
\def\ha{v}

\def\sx{{s}}
\def\wu{\tilde\tw}
\def\tw{{w}}
\def\xj{x}
\def\yj{y}
\def\xh{\xi}
\def\yh{\eta}
\def\zh{\zeta}
\def\bv{\beta}
\def\av{\alpha}

\def\hh{H}
\def\hatt{\ps^*}
\def\hu{\th^*}
\def\hv{\ha^*}
\def\wt{\widetilde T}
\def\vx{\widetilde\varXi}
\def\vk{K}
\def\vy{\varXi}

\def\vrp{\varPhi}
\def\vi{\varPsi}

\def\ly{\psi}
\def\bda{\lambda}
\def\jv{\mathbf{V}}
\def\rih{\ri^{(h)}}
\def\roh{\rho^{(h)}}
\def\vp{{\tau\hskip-4.8pt\iota\hskip.6pt}}
\def\vpab{{\tau\hskip-4.8pt\iota\hskip.2pt}}
\def\mgmt{(M,g,m\hs,\hskip.2pt\vp)}
\def\hatmgmt{(\hat M,\hat g,\hat m\hs,\hskip.2pt\hat\vp)}
\def\vs{\varSigma}
\def\vg{\varGamma}
\def\si{\phi}
\def\ta{\psi}
\def\la{\lambda}
\def\my{\mu}
\def\ny{\nu}
\def\ph{\phi}
\def\ps{t}
\def\zy{\zeta}
\def\tj{\vartheta}
\def\zu{\tilde\tz}
\def\tz{{z}}
\def\zj{\tz^*}
\def\wj{\tw^*}

\def\p{\mathcal{P}}
\def\pj{W}
\def\pw{P}
\def\lf{L}
\def\lx{\lambda}
\def\mx{\mu}
\def\rx{r}
\def\rz{r}
\def\px{p}
\def\qx{q}

\def\gy{\gamma}
\def\gm{\gamma}

\def\ro{\rho}
\def\sj{\sigma}
\def\sa{\sigma}

\def\pv{\varPi}

\def\ah{a}

\def\dv{\delta}

\def\diml{-di\-men\-sion\-al}

\def\omh{\omega^{(h)}}
\def\ve{\varepsilon}
\def\kx{\kappa}

\def\ri{\mathrm{r}}

\def\navp{\nabla\hskip-.2pt\vp}

\usepackage{amssymb,latexsym,psfrag}   

\begin{document}

\newtheorem{thm}{Theorem}[section] 
\newtheorem{lem}[thm]{Lemma}
\newtheorem{prop}[thm]{Proposition}
\newtheorem{cor}[thm]{Corollary}

\theoremstyle{definition} 
\newtheorem{defn}[thm]{Definition} 

\theoremstyle{remark} 
\newtheorem{rem}[thm]{Remark}
\newtheorem{example}[thm]{Example}
\newtheorem{conj}{Conjecture}

\renewcommand{\theequation}{\thesection.\arabic{equation}}

\title
[Conformally-Einstein K\"ahler manifolds]  
{A moduli curve for compact\\ conformally-Einstein K\"ahler manifolds}

\author{A. Derdzinski}
\email{andrzej@math.ohio-state.edu}
\address{Dept. of Mathematics, The Ohio State University,
Columbus, OH 43210, USA}
\author{G. Maschler}
\email{maschler@math.toronto.edu}
\address{Department of Mathematics, University of Toronto, Canada M5S 3G3}
%
\thanks{The work of the second author was carried out in part at the Max Planck
Institute in Bonn, and he is also partially supported by an NSERC Canada
individual research grant.}
\begin{abstract}
We classify quadruples $(M,g,m\hs,\vpab)$ in which $\,(M,g)\,$ is a
compact K\"ahler manifold of complex dimension $\,m>2\,$ with a nonconstant
function $\,\vpab\,$ on $M\hs$ such that the conformally related metric
$\,g/\vpab^{\hs2}$, defined wherever $\,\vpab\ne0\hs$, is Einstein. It turns
out that $M\hs$ then is the total space of a holomorphic $\,\bbCP^1$ bundle
over a compact K\"ahler-Einstein manifold $\,(N,h)$. The quadruples in
question constitute four disjoint families: one, well-known, with K\"ahler
metrics $\,g\,$ that are locally reducible; a second, discovered by B\'erard
Bergery (1982), and having $\,\vpab\ne0\,$ everywhere; a third one, related to
the second by a form of analytic continuation, and analogous to some known
K\"ahler surface metrics; and a fourth family, present only in odd complex
dimensions $\,m\ge9$. Our classification uses a {\it moduli curve}, which is a
subset $\,\mathcal{C}$, depending on $\,m\hs$, of an algebraic curve in
$\,\rto\nh$. A point $\,(\th,\ha)\,$ in $\,\mathcal{C}\hs$ is naturally
associated with any $(M,g,m\hs,\vpab)$ having all of the above properties
except for compactness of $M$, replaced by a weaker requirement of
``vertical'' compactness. One may in turn reconstruct $\,M,g\,$ and
$\,\vpab\,$ from this $\,(\th,\ha)\,$ coupled with some other data, among them
a K\"ahler-Einstein base $\,(N,h)\,$ for the $\,\bbCP^1$ bundle $\,M$. The
points $\,(\th,\ha)\,$ arising in this way from $(M,g,m\hs,\vpab)$ with
compact $\,M\,$ form a countably infinite subset of $\,\mathcal{C}$.
\end{abstract}

\maketitle

\setcounter{section}{-1}
\section{Introduction}\label{intr}
\setcounter{equation}{0}
This paper may be treated as a sequel to \cite{local} -- \cite{potentials},
and provides a classification, up to biholomorphic isometries, of compact
K\"ahler manifolds in complex dimensions $\,m>2\,$ that are
al\-most-eve\-ry\-where conformally Einstein. Specifically, we describe all
quadruples $\,\mgmt\,$ in which
\begin{equation}
\begin{array}{l}
M\mathrm{\hskip5.7ptis\hskip5.5pta\hskip5.5ptcompact\hskip5.5ptcomplex
\hskip5.5ptmanifold\hskip5.5ptof\hskip5.5ptcomplex\hskip5.5ptdimension
\hskip6.7pt}m\ge\nh3\mathrm{\hskip6.2ptwith\hskip5.5pta}\\
\mathrm{K}\ddot\mathrm{a}\mathrm{hler\ metric}\hskip5ptg\hskip5pt\mathrm{
and\ a\ nonconstant}
\hskip4ptC^\infty\hskip3pt\mathrm{function}\hskip6pt\vp:\nh M\nh\to\bbR
\hskip5.5pt\mathrm{such\ that\ the}\\
\mathrm{conformally\hskip4.1ptrelated\hskip4.1ptmetric\hskip6pt}\tilde g
=g/\vp^2\nh\mathrm{,\hskip5ptdefined\hskip4ptwherever\hskip6pt}\vp
\ne0\hs\mathrm{,\hskip5ptis\hskip4ptEinstein.}
\end{array}
\label{zon}
\end{equation}
When $\,m=2$, we also classify quadruples $\,\mgmt\,$ with the following
property.
\begin{equation}
\begin{array}{l}
\mathrm{Condition\ (\ref{zon})\ holds\ except\ for\ the\ requirement\ that}\hskip4pt
m\ge3\mathrm{,\ now}\\
\mathrm{replaced\ by}\hskip4ptm=2\hs\mathrm{.\ In\ addition,}\hskip5pt
d\vp\wedge\hs d\Delta\vp=0\hskip5pt\mathrm{everywhere\ in}\hskip4ptM.
\end{array}
\label{zto}
\end{equation}
Our classification of (\ref{zon}) -- (\ref{zto}) is summarized in
Theorems~\ref{cnver},~\ref{asymp} and~\ref{ctbly}. By Theorem~\ref{cnver},
$\,M\,$ must be the total space of a holomorphic $\,\bbCP^1$ bundle over a
compact K\"ahler-Einstein manifold $\,(N,h)\,$ such that the metrics $\,g\,$
and $\,h\,$ make the bundle projection $\,M\to N\,$ a horizontally homothetic
submersion \cite{gudmundsson} with totally geodesic fibres.

Theorem~\ref{cnver} also implies that every quadruple $\,\mgmt\,$ satisfying
(\ref{zon}) or (\ref{zto}) belongs to one of the four disjoint families listed
below, cf.\ Remark~\ref{fourf}.

The first and simplest family of examples with (\ref{zon}) or (\ref{zto}) involves
{\it locally reducible\/} K\"ahler metrics $\,g$. They all have $\,\vp=0\,$
somewhere in $\,M$. See \S\ref{redu}.

By constructing the corresponding conformally-K\"ahler compact Einstein
manifolds, Page \cite{page} and B\'erard Bergery \cite{berard-bergery} obtained a second
family of quadruples $\,\mgmt\,$ with (\ref{zto}) or, respectively,
(\ref{zon}) for any $\,m\ge3$. This time, $\,\vp\ne0\,$ everywhere in $\,M$.
(See \S\ref{bbpg}
and \cite[Chapter 9, Section K]{besse}.) The K\"ahler metric conformal to
Page's metric was independently discovered by Calabi \cite{calabi}, \cite{calabi-ii},
\cite{chave-valent}. It follows from our classification that Page's and B\'erard
Bergery's examples just mentioned are the only quadruples $\,\mgmt\,$ with (\ref{zto}) or (\ref{zon}) which are {\it globally\/} conformally Einstein (that
is, $\,\vp\ne0\,$ everywhere in $\,M$).

More recently, Hwang and Simanca \cite{hwang-simanca} and T\o nnesen-Friedman
\cite{tonnesen-friedman} provided examples of (\ref{zto}) on minimal ruled surfaces $\,M\,$
with $\,\vp=0\,$ somewhere in $\,M$. We extend their construction by
describing, in every complex dimension $\,m\ge2$, a third family of quadruples
$\,\mgmt\,$ that satisfy (\ref{zon}) or (\ref{zto}) and have $\,\vp=0\,$ at some
point of $\,M$. This third family is still closely related, through a form of
analytic continuation, to B\'erard Bergery's and Page's second family, while,
for $\,m=2$, it consists precisely of the examples already given in
\cite{hwang-simanca} and \cite{tonnesen-friedman}. See \S\ref{thrd}.

Finally, in every odd complex dimension $\,m\ge9$, we exhibit in
\S\ref{frth} a new, fourth family of quadruples with (\ref{zon}). It is
distinguished by a natural notion of {\it duality} (Remark~\ref{invol}):
every quadruple in the first three families is its own dual, but none in the
fourth family is.

The {\it moduli curve\/} mentioned in the title plays a prominent role in our
construction. It is a subset $\,\mathcal{C}\hs$ of an algebraic curve in
$\,\rto\nh$, depending also on the complex dimension $\,m\ge2$. Any quadruple
$\,\mgmt\,$ with (\ref{zon}) or (\ref{zto}) gives rise to a point
$\,(\th,\ha)\in\mathcal{C}$, defined as follows. If $\,g\,$ is locally
reducible as a K\"ahler metric, $\,(\th,\ha)=(0\hs,\nh0)$. Otherwise, we set
$\,\th=\,\mathrm{min}\hskip3pt\vp/\y\,$ and
$\,\ha=\,\mathrm{max}\hskip3pt\vp/\y\hs$, with
$\,\y\in\bbR\smallsetminus\{0\}\,$ characterized by the property that
$\,|\navp|^2$ is a rational function of $\,\vp\,$ having a unique real pole at
$\,\y\,$ (Remark~\ref{detrm}). In this way one obtains not all
$\,(\th,\ha)\in\mathcal{C}$, but only a countably infinite set of points that
we call $\,\px\hs${\it-rational}. See Remark~\ref{prati}.

The quadruple $\,\mgmt\,$ can in turn be explicitly reconstructed (see
\S\ref{stat}) from the corresponding $\,\px\hs$-rational point $\,(\th,\ha)\,$
coupled with some additional data which include a compact K\"ahler-Einstein
manifold $\,(N,h)\,$ with $\,\dimc N=m-1\,$ such that $\,M\,$ is a holomorphic
$\,\bbCP^1$ bundle over $\,N$. The present paper provides a proof of this fact
and a relatively detailed description of the set of $\,\px\hs$-rational
points; our four families arise when one requires $\,(\th,\ha)\,$ to lie in
one of four specific subsets of $\,\mathcal{C}\,$ (see Remark~\ref{fourf}).
The number of $\,\px\hs$-rational points in each of the four subsets is {\it
one} for the first family (\S\ref{redu}) and {\it infinite} for the third
family (\S\ref{thrd} and Theorem~\ref{ctbly}). For the second and fourth
families this number is {\it finite} and varies with $\,m\,$ so that, as
$\,m\to\infty$, it is asymptotic to a positive constant times $\,m^2$
(Theorem~\ref{asymp}).

Although only $\,\px\hs$-rational points are directly used in our
classification of (\ref{zon}) and (\ref{zto}), the other points of
$\,\mathcal{C}\hs$ have a similar geometric interpretation. Namely, the last
two paragraphs are valid even if one replaces $\,\px\hs$-rational points with
arbitrary points of $\,\mathcal{C}$, provided that, instead of compactness of
$\,M\,$ and $\,N$, one only requires $\,M\,$ to be {\it vertically
compact\hs}, which amounts to compactness of the $\,\mathbf{C}\mathrm{P}^1$
fibres, but not necessarily of the base $\,N\nh$. See \S\ref{vert}.

\section{Statement of the main results}\label{stat}
\setcounter{equation}{0}
In this section $\,m\,$ is an integer with $\,m\ge2\,$ and $\,\th,\ha\,$ are
the Cartesian coordinates in $\,\rto\nh$. Most objects discussed here depend
on the choice of $\,m\hs$.

We denote by $\,\mathcal{H}\,$ the hyperbola $\,\th\ha=\th+\ha\,$ in $\,\rto$
and let $\,\mathcal{T}\subset\rto$ be the set given by $\,T(\th,\ha)=0$, where
$\,T\,$ is the symmetric polynomial of degree $\,3(m-2)\,$ described in
Lemmas~\ref{polyt} and~\ref{xpant}. Thus, $\,T\,$ and $\,\mathcal{T}\,$ depend
on $\,m$. See Fig.\ 1.

The {\it moduli curve\/} corresponding to a given value of $\,m\,$ is a subset
$\,\mathcal{C}\hs$ of the half-plane $\,\th\le\ha\,$ in $\,\rto$ such that
$\,\mathcal{C}\subset\hut\,$ for odd $\,m\hs$, while, if $\,m\,$ is even,
$\,\mathcal{C}=\{(\th,\ha)\in\mathcal{H}:(\th,\ha)\ne(1\hs,\nh1)\hskip4.5pt
\mathrm{and}\hskip4.5pt\th\le\ha\}$. Thus, $\,\mathcal{C}\hs$ is the same for
all even $\,m\hs$.
\begin{figure}
\centering
   \psfrag{1}[bl][bl][0.7][0]{$1$}
   \psfrag{2}[bl][bl][0.7][0]{$2$}
   \psfrag{]}[bl][bl][0.7][0]{$-1$}
   \psfrag{3}[bl][bl][0.7][0]{$-2$}
   \psfrag{u}[bl][bl][0.9][0]{$u$}
   \psfrag{v}[bl][bl][0.9][0]{$v$}
   \psfrag{T}[bl][bl][0.8][0]{$\mathcal{T}$}
   \psfrag{d}[bl][bl][0.8][0]{$\mathcal{T}$}
   \psfrag{t}[bl][bl][0.8][0]{$\mathcal{T}$}
   \psfrag{P}[bl][bl][0.8][0]{$\mathcal{H}$}
   \psfrag{Q}[bl][bl][0.8][0]{$\mathcal{H}$}
   \psfrag{|}[bl][bl][0.8][0]{\hskip1pt{\tensans I}}
   \psfrag{!}[bl][bl][0.8][0]{\di}
   \psfrag{X}[bl][bl][0.8][0]{\hskip1pt{\tensans X}}
   \psfrag{M}[bl][bl][0.8][0]{$m=3$}
   \hskip19pt
   \includegraphics[scale = 0.477]{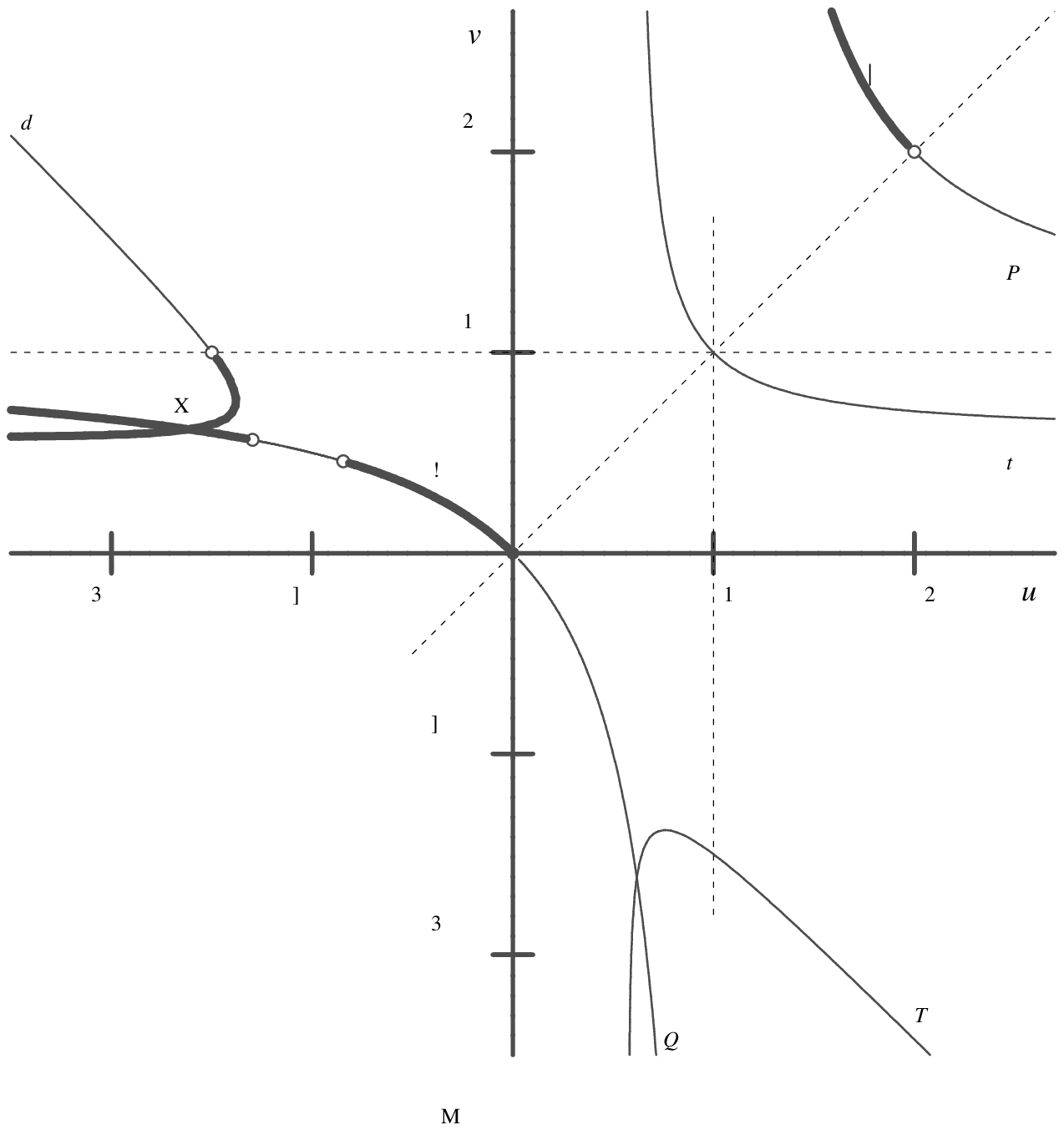}
   \psfrag{1}[bl][bl][0.7][0]{$1$}
   \psfrag{2}[bl][bl][0.7][0]{$2$}
   \psfrag{]}[bl][bl][0.7][0]{$-1$}
   \psfrag{3}[bl][bl][0.7][0]{$-2$}
   \psfrag{u}[bl][bl][0.9][0]{$u$}
   \psfrag{v}[bl][bl][0.9][0]{$v$}
   \psfrag{P}[bl][bl][0.8][0]{$\mathcal{H}$}
   \psfrag{Q}[bl][bl][0.8][0]{$\mathcal{H}$}
   \psfrag{|}[bl][bl][0.8][0]{\hskip1pt{\tensans I}}
   \psfrag{!}[bl][bl][0.8][0]{\di}
   \psfrag{N}[bl][bl][0.8][0]{$m\,$ even}
   \hskip9.2pt
   \includegraphics[scale = 0.477]{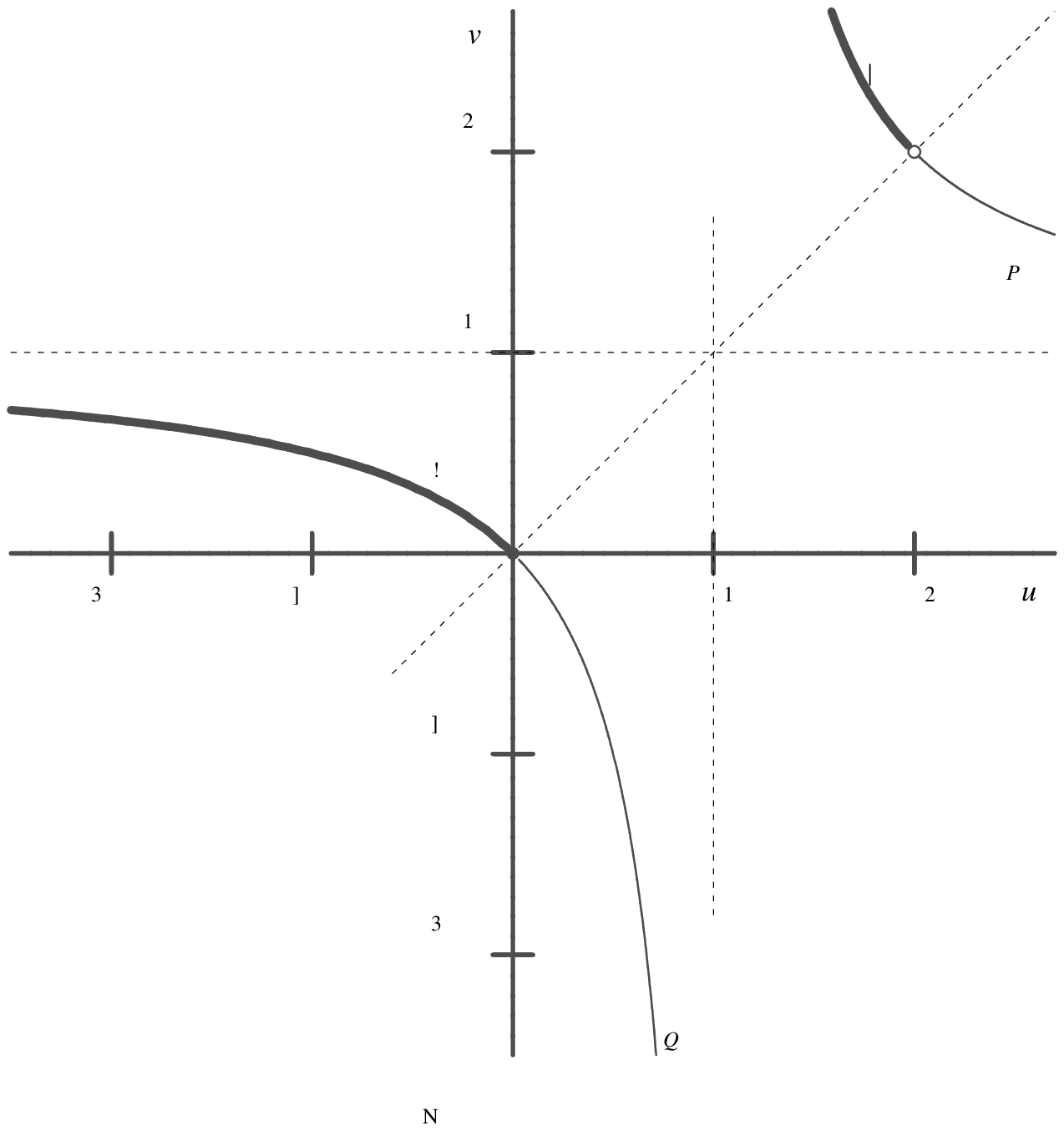}
\caption{The components
$\,\hs\di\hs,\,\mathsf{I}\hs\,$ of $\,\mathcal{C}\hs$ and, for odd
$\,m\hs$, the two beams of $\,\mathsf{X}\hs$, are the heavily marked curve
segments, each contained in $\,\mathcal{H}\,$ except for the
$\hs\mathcal{T}\nh$-beam of $\,\mathsf{X}\hs$, contained in $\hs\mathcal{T}$.
For $\,m=3$, the set $\hs\mathcal{T}$ consists of the isolated point
$\,(0\hs,\nh0)\,$ and three disjoint real-analytic curves in $\hs\rto\nh$,
diffeomorphic to $\hs\bbR\hs$. (See Remark~\ref{graph}.)}
\vskip-18pt
\end{figure}
The following explicit description of $\,\mathcal{C}$, although different from
the  definition given at the end of \S\ref{modu}, is equivalent to it (see
Theorem~\ref{modcu}). Namely, $\,\mathcal{C}\hs$ is the disjoint union of its
connected components (defined below):
\begin{equation}
\mathcal{C}\,=\hskip6pt\di\hskip5pt\cup\hskip5pt\mathsf{I}\hskip16pt{\rm if\ }
\,m\,{\rm\ is\ even,\hskip8ptand\hskip16pt}\mathcal{C}\,=\hskip6pt\mathsf{X}
\hskip4pt\cup\hskip5pt\di\hskip5pt\cup\hskip5pt\mathsf{I}\hskip16pt{\rm if\ }
\,m\,{\rm\ is\ odd.}\label{xii}
\end{equation}
Specifically, the $\,\,\mathsf{X}\,\,$ component exists for odd $\,m\,$
only and is contained in $\,\hut$, so that it is the union of its
$\,\mathcal{H}\nh${\it-beam\hs} $\,\xah\,$ and
$\,\mathcal{T}\hskip-1.3pt${\it-beam\hs} $\,\xat$. The first two of the sets
$\,\hs\mathsf{X}\hs,\,\hs\di\hs\,$ and $\,\hs\mathsf{I}\hs\,$ in $\,\rto$
depend on $\,m\hs$, and the $\,\mathcal{T}\nh$-beam $\,\xat\,$ of
$\,\hs\mathsf{X}\hs\,$ is the intersection of $\,\mathcal{T}\hs$ with
$\,(-\infty,0)\times(0,1)$. Next, $\,\xah\hs$, $\,\di\,$ and
$\,\,\mathsf{I}\,\,$ are subsets of the hyperbola $\,\mathcal{H}$, namely, the
segments of $\,\mathcal{H}\,$ that project onto the following intervals in the
$\,\th\,$ axis: $\,(-\infty,\tz)\,$ (for $\,\xah\,$ and odd $\,m$), or
$\,(\tw,0\hs]\,$ (for $\,\di\,$ and odd $\,m$), or $\,(-\infty,0\hs]\,$ (for
$\,\di\,$ and even $\,m$), or $\,(1,2)\,$ (for $\,\,\mathsf{I}\,\,$ and any
$\,m$), with constants $\,\tz,\tw\,$ such that $\,\tz<\tw<0$, defined in
\S\ref{anth}, and depending on the odd integer $\,m\hs$.

In \S\ref{anth} the symbol $\,\tz\,$ is assigned a meaning also when $\,m\,$
is even; the constant $\,\tz\in(-\infty,0)\,$ depending on an {\it even\hs}
integer $\,m\,$ appears in the following definition of a function
$\,\dv:\mathcal{C}\to\{-\hs1,0,1\}$, which also depends on $\,m\hs$. Namely,
we set $\,\dv=1\,$ both on $\,\,\mathsf{X}\,\,$ for odd $\,m\,$ and on
$\,\,\mathsf{I}\,\,$ for all $\,m\hs$, as well as $\,\dv=\hs-1\,$ on
$\,\,\di\,\,$ for odd $\,m\,$ and $\,\dv=\,\sgn\hs(\tz-\th)\,$ at any
$\,(\th,\ha)\in\,{\di}\hs\,$ for even $\,m\hs$. Thus, when $\,m\,$ is even,
$\,\,\di\,\,$ contains a distinguished point $\,(\tz,\zj)\,$ with
$\,\zj=\tz/(\tz-1)$, at which $\,\dv=0$.

The symbols $\,\,\mathsf{X}\hs,\,\,\di\,\,$ and
$\,\,\mathsf{I}\,\,$ imitate the topology of the sets in question:
$\,\,\mathsf{I}\hs$, $\,\xah$, $\,\xat\,$ and $\,\,\di\,\,$ are
real-analytic submanifolds of $\,\rto$ and $\,\,\di\,\,$ is diffeomorphic to
$\,(0,1\hs]$, while the other three are diffeomorphic to $\,\bbR\hs$, and the
two beams $\,\xah\,$ and $\,\xat\,$ forming $\,\,\mathsf{X}\,\,$ have a
single transverse intersection point. See Proposition~\ref{crvsg}(ii) and
Remark~\ref{trnsv}.

We let $\,\px\,$ stand for a specific rational function of the variables
$\,\th,\ha$, depending on $\,m\hs$, which is defined in \S\ref{trfp}. The
restriction of $\,\px\,$ to $\,\mathcal{C}\hs$ is finite (that is,
well-defined) and nonzero everywhere in
$\,\mathcal{C}\smallsetminus\{(0\hs,\nh0)\}\,$ (for odd $\,m$) or in
$\,\mathcal{C}\smallsetminus\{(0\hs,\nh0),(\tz,\zj)\}\,$ (for even $\,m$). In
addition, $\,\px=0\,$ at $\,(0\hs,\nh0)\,$ for any $\,m\hs$, while, for even
$\,m\,$ only, $\,\px\,$ is undefined at $\,(\tz,\zj)$.
\begin{defn}\label{pratl}
Given a fixed integer $\,m\ge2$, by a
$\,\px\hs${\it-rational point\/} we mean any $\,(\th,\ha)\in\mathcal{C}\hs$ at
which either $\,\dv=0\,$ (that is, $\,m\,$ is even and
$\,(\th,\ha)=(\tz,\zj)$), or $\,\dv=-1\,$ and the value of $\,\px\,$ is
rational, or, finally, $\,\dv=1\,$ and $\,\px\,$ equals $\,\ka/\el\,$ for some
$\,\ka\in\bbZ\,$ and $\,\el\in\{1,\dots,m\}$. In particular,
$\,(0\hs,\nh0)\in\,\di\,\subset\mathcal{C}\,$ is $\,\px\hs$-rational.
\end{defn}
We will now describe how $\,\px\hs$-rational points on $\,\mathcal{C}\hs$ are
related, for any given $\,m\ge2$, to quadruples $\,\mgmt\,$ with (\ref{zon})
or (\ref{zto}). A similar geometric interpretation is valid for {\it
arbitrary} points of $\,\mathcal{C}$, provided that the compactness
requirement in (\ref{zon}) or (\ref{zto}) is suitably relaxed. (See
\S\ref{vert}.)

Every $\,\px\hs$-rational point $\,(\th,\ha)\in\mathcal{C}$, for any fixed
integer $\,m\ge2$, corresponds to a mapping which assigns a quadruple
$\,\mgmt\,$ with (\ref{zon}) or (\ref{zto}) to any suitable input data. Those
data consist of a compact K\"ahler-Einstein manifold $\,(N,h)\,$ with
$\,\dimc N=m-1$, a holomorphic line bundle $\,\mathcal{L}\,$ over $\,N$, a
$\,\hs\mathrm{U}\hs(1)\,$ connection in $\,\mathcal{L}\,$ with a curvature
form $\,\varOmega$, and a constant $\,\y\in\bbR\smallsetminus\{0\}$, which
satisfy the following additional requirements, depending on $\,(\th,\ha)$, and
involving the Ricci form $\,\roh$ of $\,h$, its K\"ahler form $\,\omh$, and
the signum of the Einstein constant $\,\kx\,$ of $\,(N,h)$. Namely,
$\,\hs\sgn\,\kx\,$ must equal the value of $\,\dv:\mathcal{C}\to\{-1,0,1\}\,$
at $\,(\th,\ha)$, and, for the rational function $\,\px\,$ mentioned above,
one of the following two cases has to occur:
\begin{equation}
\begin{array}{rl}\arraycolsep11pt
\mathrm{i)}&\kx\ne0\,\,\mathrm{\ and\ }\,\,\varOmega\,\,\mathrm{\ equals\ }
\,\,\roh\hs\mathrm{\ times\ the\ value\ of\ }\,\,\px\,\,\mathrm{\ at\ }
\,\,(\th,\ha),\\
\mathrm{ii)}&\kx=0\,\,\mathrm{\ and\ }\,\,\varOmega\,\,\mathrm{\ is\ a\
nonzero\ real\ multiple\ of\ }\,\,\omh\hskip-1.2pt.
\end{array}\label{kpo}
\end{equation}
Note that $\,\roh$ and $\,\varOmega$, divided by $\,2\pi$, represent the real
Chern classes $\,c_1(N)\,$ and $\,c_1(\mathcal{L})\,$ in $\,H^2(N,\bbR)$.
Also, $\,\roh=\kx\hskip1.2pt\omh$. Thus, in case (i), $\,c_1(N)\ne0\,$ and
$\,c_1(\mathcal{L})\,$ equals $\,c_1(N)\,$ times the value of $\,\px\,$ at
$\,(\th,\ha)$. Similarly, in case (ii), $\,c_1(N)=0\,$ and
$\,c_1(\mathcal{L})\,$ is a nonzero real multiple of the K\"ahler class of
$\,h$. Moreover, $\,\varOmega=0\,$ in (i) only if $\,(\th,\ha)=(0\hs,\nh0)$,
while (ii) occurs only when $\,m\,$ is even and $\,(\th,\ha)=(\tz,\zj)$.

The quadruple $\,\mgmt\,$ is constructed out of such data as follows. If
$\,(\th,\ha)=(0\hs,\nh0)$, we proceed as described in \S\ref{redu}. Now let
$\,(\th,\ha)\in\mathcal{C}\smallsetminus\{(0\hs,\nh0)\}$. By
Remark~\ref{uvone}, $\,1<\th<\ha\,$ or $\,\th<\ha<1$, that is,
$\,I=[\th,\ha]\,$ is a nontrivial closed interval and $\,1\notin I$. Hence
we may choose $\,\ve\in\{1,-1\}\,$ with $\,\ve\hs\y\hs(\ps-1)>0\,$ for every
$\,\ps\,$ in the interior of $\,I$, where $\,\y\,$ comes from our data. In
terms of $\,\kx\,$ and the value of $\,\px\,$ appearing in (\ref{kpo}), we
define $\,\ah\in\bbR\,$ either by $\,\ah=-\hs\ve\hs\px\kx/2\,$ (case
(\ref{kpo}.i)), or by $\,\varOmega=-\hs2\ve\hskip.2pt\ah\hs\omh$ (case
(\ref{kpo}.ii)). On the other hand, the definition of $\,\mathcal{C}\hs$ at
the end of \S\ref{modu} guarantees the existence of a rational function
$\,Q\,$ satisfying some specific conditions that involve $\,I$. Such $\,Q$,
unique up to a constant factor (Lemma~\ref{onedi}), is now made unique by
requiring (\ref{cea}) to hold for $\,Q\,$ and $\,\y,\ve,\ah$. Our data and
$\,\ve,\ah,\px,\kx\,$ also satisfy (\ref{nho}), as one sees using
Remark~\ref{kemac} and noting that, as $\,\roh=\kx\hskip1.2pt\omh\nh$, our
choice of $\,\ah\,$ gives $\,\varOmega=-\hs2\ve\hskip.2pt\ah\hs\omh$ in {\it
both} cases (\ref{kpo}.i), (\ref{kpo}.ii).

Applying the construction described in \S\ref{main} we now obtain the required
quadruple $\,\mgmt$. In particular, $\,M\,$ is the projective compactification
of the total space of $\,\mathcal{L}\hs$, that is, the bundle of Riemann
spheres associated with $\,\mathcal{L}\hs$.
\begin{prop}\label{quadr}For any $\,\px\hs$-rational point on the moduli
curve $\,\mathcal{C}$, with any integer $\,m\ge2$, additional data as above
exist, and, applying to them the construction just described, we always obtain
a quadruple $\,\mgmt\,$ with {\rm(\ref{zon})} or {\rm(\ref{zto})}.
\end{prop}
Namely, the data exist according to Remark~\ref{prati}, while (\ref{zon}) or (\ref{zto}) is verified in \S\ref{main}.  On the other hand, in \S\ref{more} we prove the
following result:
\begin{thm}\label{cnver}Conversely, every quadruple\/ $\,\mgmt\,$ with\
{\rm(\ref{zon})} or\/ {\rm(\ref{zto}}) is, up to a $\,\vp$-preserving
biholomorphic isometry, obtained from the above construction applied to some
$\,\px\hs$-rational point\/ $\,(\th,\ha)\in\mathcal{C}\hs$ and some additional
data with the properties just listed.
\end{thm}
\begin{rem}\label{fourf}
The four families of quadruples with (\ref{zon}) or
(\ref{zto}), mentioned in \S\ref{intr}, correspond to a decomposition of the
moduli curve $\,\mathcal{C}\hs$ into the following four disjoint subsets, from
which the $\,\px\hs$-rational points then are chosen. For the first family, it
is the subset consisting of $\,(0\hs,\nh0)\,$ alone; for the second, the
$\,\,\mathsf{I}\,\,$ component; for the third,
$\,\,\di\hs\smallsetminus\{(0\hs,\nh0)\}$, augmented (only when $\,m\,$ is
odd) by the $\,\mathcal{H}$-beam $\,\xah\,$ of the $\,\,\mathsf{X}\,\,$
component; and, for the fourth family, the empty set if $\,m\,$ is even, or,
if $\,m\,$ is odd, the $\,\mathcal{T}\nh$-beam $\,\xat\,$ minus its unique
intersection point with the $\,\mathcal{H}\nh$-beam.
\end{rem}
Proposition~\ref{quadr} and Theorem~\ref{cnver} together form a classification
result for (\ref{zon}) and (\ref{zto}). Its meaning, however, remains obscure
unless one addresses the questions of abundance or scarcity, including that of
the very existence, of $\,\px\hs$-rational points corresponding to each of the
four families. The next two results, proved in sections \ref{bbpg} --
\ref{frth} and \ref{valu}, provide information of this kind.
\begin{thm}\label{asymp}For any given integer $\,m\ge2$, the set of\/
$\,\px\hs$-rational points is a countably infinite subset of the moduli curve
$\,\mathcal{C}$. Its intersections with $\,\,\mathsf{I}\hs$, $\,\xah$ and\/
$\,\,\mathsf{X}\,\cap\,(\mathcal{T}\smallsetminus\mathcal{H})\hs$, which we
denote, respectively, by $\,\smi,\,\smh$ and\/ $\,\smt$, are all finite, and
their cardinalities satisfy the asymptotic relations
\[
\begin{array}
{rllll}\arraycolsep9pt
\mathrm{a)}\quad&|\hs\smi|\hs/m^2&\to&3\hs/\pi^2&\mathrm{as}\hskip10ptm\to
\infty\hs,\\
\mathrm{b)}\quad&|\hs\smh|\hs/m^2&\to&\,9\sqrt2\hs(\sqrt2+1)\hs/\pi^2\quad&
\mathrm{as\ \ (odd)}\hskip5ptm\to\infty\hs,\\
\mathrm{c)}\quad&|\hs\smt|\hs/m^2&\to&\,(6-4\sqrt2\hs)/\pi^2&\mathrm{as\ \
(odd)}\hskip5ptm\to\infty\hs,
\end{array}
\]
where $\,\,\mathsf{X}\hs$, $\,\smh$ and $\,\smt$ are defined for odd\/
$\,m\,$ only. The values assumed by $\,\px\,$ on $\,\smh$, or
$\,\smt$, all lie in the interval $\,(-\hs2,0)$, or, respectively,
$\,(-\hs1,0)$. In addition,
\begin{enumerate}\alph{enumi}
\item[d)] $\px\,$ maps $\,\smi$ bijectively onto the set of all rational
numbers in $\,(0,1)\,$ with positive denominators not exceeding $\,m$, cf.\
{\rm Remark~\ref{farey}}.
\item[e)] If\/ $\,m\,$ is odd, $\,\px\,$ assumes the values $\,-\hs1\,$ and\/
$\,-\hs1+1/m\,$ at a unique pair of points of $\,\smh$. Therefore,
$\,|\hs\smh|\ge2$.
\item[f)] $\smt\hskip-1.5pt=$\hskip5pt\emp\hskip7ptif\/ $\,m\in\{3,5,7\}$,
while $\,|\hs\smt|\ge2\,$ if\/ $\,m\,$ is odd and\/ $\,m\ge9$. In particular,
$\,|\hs\mathcal{S}_9^{\mathcal{T}}|=2\,$ and\/ $\,|\hs\smt|\ge2(m-17)\,$ for
any odd\/ $\,m\ge19$.
\end{enumerate}
\end{thm}
The limits in (a) -- (c) are, approximately, 0.304, 3.113 and 0.035. For the
remaining component \ \di\hs, we have the following theorem.
\begin{thm}\label{ctbly}Given an integer $\,m\ge2$, the set of\/
$\,\px\hs$-rational points in the $\,\,\di\,$ component of\/
$\,\mathcal{C}\hs$ is countably infinite.
\begin{enumerate}
\item If\/ $\,m\,$ is odd, the values of\/ $\,\px\,$ on\/ $\,\,\di\,$ all lie
in $\,(0,1/m)$.
\item For even $\,m\hs$, the set of those $\,\px\hs$-rational points in
by $\,\px\,$ bijectively onto the set of rational numbers in
$\,(-\infty,-1)\,$ that have positive denominators not exceeding $\,m\,$ or,
respectively, onto the set of rational numbers in $\,[\hs0,\infty)$.
\end{enumerate}
\end{thm}

\section{The moduli curve}\label{modu}
\setcounter{equation}{0}
Following \cite{local} we set, for any fixed integer $\,m\ge1$,
\begin{equation}
F(\ps)\,=\,{(\ps-2)\ps^{2m-1}\over(\ps-1)^m}\,,\qquad E(\ps)\,
=\,(\ps-1)\,\sum_{j=1}^m{j\over m}{2m-j-1\choose m-1}\ps^{j-1}\hs.\label{fet}
\end{equation}
Thus, $\,F\,$ and $\,E\,$ are rational functions of the real variable
$\,\ps\,$ and
\begin{equation}
E(\ps)\,=\,(\ps-1)\vs(\ps)\,,\qquad\mathrm{where}\qquad
\vs(\ps)\,=\,\sum_{j=1}^m{j\over m}{2m-j-1\choose m-1}\ps^{j-1}\hs.
\label{esg}
\end{equation}
Any real constants $\,\ax,\bx,\cx\,$ now give rise to a rational function
$\,Q\,$ with
\begin{equation}
Q(\ps)\,=\,(\ps-1)\,\lj\hs\ax\,+\,\bx E(\ps)\,+\,\cx F(t)\hs\rj\,,\qquad
\mathrm{for}\hskip5ptE,F\hskip5pt\mathrm{as\ in\ (\ref{fet}).}\label{qef}
\end{equation}
(Cf.\ \cite[formula~(21.5)]{local}.) For a fixed integer $\,m\ge2$, let
\begin{equation}
\jv\,=\,\,\mathrm{Span}\hs\{\ps-1,\,(\ps-1)E,\,(\ps-1)F\}\,,\qquad
\mathrm{with}\hskip9ptE,F\hskip7pt\mathrm{as\ in\ (\ref{fet}),}
\label{spv}
\end{equation}
$\ps\,$ being the identity function. Thus, $\,\jv\,$ is the $\,3$\diml\ real
vector space, depending on $\,m\hs$, of all rational functions of the form
(\ref{qef}) with $\,\ax,\bx,\cx\in\bbR\hs$.

Given a nontrivial closed interval $\,I\subset\bbR\hs$, we may impose on
$\,Q\in\jv\,$ and $\,I\hs$ one or more of the following five conditions (cf.\
\cite[formula~(34.2)]{potentials}):
\begin{equation}
\begin{array}{rl}\arraycolsep13pt
\mathrm{a)}&Q\,\mathrm{\ is\ analytic\ on\ }\,I\mathrm{,\ that\ is,\ }
\,I\hs\mathrm{\ does\ not\ contain\ }\,1\,\mathrm{\ unless\ }\,\cx=0\,
\mathrm{\ in\ (\ref{qef}).}\\
\mathrm{b)}&Q=0\,\mathrm{\ at\ both\ endpoints\ of\ }\,I.\\
\mathrm{c)}&Q\ne0\,\mathrm{\ at\ all\ interior\ points\ of\ }\,I\hs
\mathrm{\ at\ which\ }\,Q\,\mathrm{\ is\ analytic.}\\
\mathrm{d)}&\dot Q=\,dQ/dt\,\mathrm{\ exists\ and\ is\ nonzero\ at\
both\ endpoints\ of\ }\,I.\\
\mathrm{e)}&\mathrm{The\ values\ of\ }\,\dot Q\,\mathrm{\ at\ the\ endpoints\
of\ }\,I\hs\mathrm{\ exist\ and\ are\ mutually\ opposite.}
\end{array}
\label{qdq}
\end{equation}
\begin{rem}\label{endpo}
Conditions (\ref{qdq}.a,\hs c) alone imply
that $\,1\,$ cannot be an interior point of $\,I$. In fact, if $\,1\in I$,
(\ref{qdq}.a) gives $\,\cx=0$, and so $\,Q(1)=0\,$ due to (\ref{qef}). Thus,
by (\ref{qdq}.c), $\,1\,$ is an endpoint of $\,I$.
\end{rem}
Let $\,m\ge2\,$ be a fixed integer. We define the {\it moduli curve\/} to
be the set $\,\mathcal{C}\subset\rto$ consisting of $\,(0\hs,\nh0)\,$ and all
$\,(\th,\ha)\,$ with $\,\th<\ha\,$ for which there exists a function $\,Q\,$
in the space $\,\jv\,$ given by (\ref{spv}), satisfying all of (\ref{qdq}) on
the interval $\,I=[\th,\ha]$.

\section{The main step in the construction}\label{main}
\setcounter{equation}{0}
The simplest examples of quadruples $\,\mgmt\,$ satisfying (\ref{zon}) or
(\ref{zto}) are described in \S\ref{redu}; they involve K\"ahler metrics
$\,g\,$ that are locally reducible. The following construction, which also
appears in \cite{potentials}, leads to quadruples with (\ref{zon}) or
(\ref{zto}) that are {\it not\/} locally reducible (cf.\ Remark~\ref{irred}).

Let a function $\,Q\in\jv$, with $\,\jv\,$ as in (\ref{spv}) for a fixed
integer $\,m\ge2$, satisfy all five conditions (\ref{qdq}) on a nontrivial
closed interval $\,I=[\th,\ha]\,$ with $\,1\notin I$. Replacing $\,Q\,$ by
$\,-\hs Q$, if necessary, we also require (cf.\ (\ref{qdq}.c)) that $\,Q>0\,$
on the open interval $\,(\th,\ha)$. We now choose
$\,\y,\ve,\ah\in\bbR\,$ with
\begin{equation}
\begin{array}{l}
\ve\in\{1,-1\}\,\mathrm{\ and\ }\,\,\ve\hs\y\hs(\ps-1)>0\,\,\mathrm{\
for\ every\ }\,\ps\,\mathrm{\ in\ the\ open\ interval\ }\,(\th,\ha),\\
\mathrm{while\ \ }\dot Q(\th)=-\hs2\ah\y\mathrm{\ \ and\ \ }\dot Q(\ha)
=2\ah\y\mathrm{,\ \ with\ \ }\dot Q=\hs dQ/d\ps.
\end{array}
\label{cea}
\end{equation}
(Such $\,\y,\ve,\ah\,$ must exist by (\ref{qdq}.a\hs,\hs c\hs,\hs e).) Next,
let there be given
\begin{equation}
\begin{array}{l}
\mathrm{a\ \ compact\ \ K}\ddot\mathrm{a}\mathrm{hler}\hyp\mathrm{Einstein\ \
manifold\ \ }\,\hs(N,h)\hs\,\mathrm{\ \ of\ \ complex\ \ dimension\ \ }\,
\hs m-1\\
\mathrm{ having\ the\ Ricci\ form\ \ }\roh\nh=\kx\hskip.9pt\omh\mathrm{\ \
for\ \ }\kx=\ve\hs m\ax/\y\hs\mathrm{,\hs\ a\ complex\ line\ bundle\ \ }
\mathcal{L}\\
\mathrm{over\ }\,N\nh\mathrm{,\ and\ a\ }\,\hs\mathrm{U}\hs(1)\,\mathrm{\
connection\ in\ }\,\mathcal{L}\,\mathrm{\ with\ the\ curvature\ form\ }\,
\varOmega=-\hs2\ve\hskip.2pt\ah\hs\omh\hskip-1.3pt,
\end{array}
\label{nho}
\end{equation}
where $\,\ax\,$ is determined by $\,Q\,$ via (\ref{qef}) and $\,\omh$ denotes
the K\"ahler form of $\,(N,h)$. The question whether such objects exist is
discussed in \S\ref{rati}.

With $\,m,Q,I,\y,\ve,\ah\,$ and the objects (\ref{nho}) fixed as above, let us
also choose a positive function $\,r\,$ of the variable $\,\ps\,$ restricted
to the interior of $\,I$, such that $\,dr/dt=\,\ah\y\hs r/Q$. By (\ref{qdq}.b)
-- (\ref{qdq}.d), $\,\log\hs r\,$ and $\,r\,$ have the ranges
$\,(-\infty,\infty)\,$ and $\,(0,\infty)$. We may thus treat $\,\ps\,$ along
with $\,\vp=\y\hs\ps\,$ and $\,Q$, restricted to the interior of $\,I$, as
functions of a new variable $\,r\in(0,\infty)$, so that
$\,\th=\,\mathrm{inf}\hskip3pt\ps\,$ and $\,\ha=\,\mathrm{sup}\hskip3pt\ps\,$
for $\,\ps:(0,\infty)\to\bbR\hs$.

The total space of the line bundle in (\ref{nho}) is denoted by the same
symbol $\,\mathcal{L}$, and $\,r\,$ also stands for the function
$\,\mathcal{L}\to(0,\infty)\,$ which, restricted to each fibre, is the norm
corresponding to the $\,\hs\mathrm{U}\hs(1)\,$ structure. Being functions of
$\,r>0$, both $\,\ps\,$ and $\,\vp=\y\hs\ps$, as well as $\,Q$, become
functions on $\,\mathcal{L}\smallsetminus N\nh$, where
$\,N\subset\mathcal{L}\,$ is the zero section.

We now define a metric $\,g\,$ on the complex manifold
$\,\mathcal{L}\smallsetminus N\,$ by letting $\,g\,$ on each fibre of
$\,\mathcal{L}\,$ coincide with $\,Q/(\ah r)^2\,$ times the standard Euclidean
metric, declaring the horizontal distribution of the connection in
$\,\mathcal{L}\,$ to be $\,g$-normal to the fibres, and requiring that $\,g\,$
restricted to the horizontal distribution equal $\,2|\vp-\y\hs|\,$ times the
pullback of $\,h\,$ under the projection $\,\mathcal{L}\to N$.

Finally, let the compact complex manifold $\,M\,$ be the projective
compactification of $\,\mathcal{L}$, that is, the Riemann sphere bundle
obtained when the total spaces of $\,\mathcal{L}\,$ and its dual
$\,\mathcal{L}^*$ are glued together by the biholomorphism
$\,\mathcal{L}\smallsetminus N\to\mathcal{L}^*\smallsetminus N\,$ which sends
each $\,\ph\in\mathcal{L}_y\smallsetminus\{0\}$, $\,y\in N$, to the unique
$\,\chi\in\mathcal{L}^*_y$ with $\,\chi(\ph)=1$.

Theorem~34.3 in \cite{potentials} now shows that $\,g\,$ and
$\,\vp=\y\hs\ps\,$ have $\,C^\infty$ extensions to $\,M$ such that the
resulting quadruple $\,\mgmt\,$ satisfies (\ref{zon}) or (\ref{zto}).

This yields Proposition~\ref{quadr}, except for the existence assertion,
proved later in Remark~\ref{prati}. Also, the above inf\hs/sup relations give
$\,\th=\,\mathrm{min}\hskip3pt\ps$, $\,\hs\ha=\,\mathrm{max}\hskip3pt\ps\,$ on
$\,M$.
\begin{rem}\label{rescl}
Equation $\,dr/dt=\,\ah\y\hs r/Q\,$ determines $\,r\hs$ only up to a positive
constant factor. A different choice of $\,r\hs$ thus amounts to rescaling the
Hermitian fibre metric in $\,\mathcal{L}$, and the resulting quadruple is
equivalent to the original one under an obvious biholomorphic isometry.
\end{rem}
\begin{rem}\label{detrm}
The constants $\,\y,\th,\ha\,$ used in the above
construction are in turn uniquely determined by the biholomorphic-isometry
type of the resulting quadruple $\,\mgmt$. Specifically, one easily sees that
the $\,g$-gradient $\,\navp\,$ of $\,\vp\,$ equals $\,\ah\,$ times the
``identity'' vertical field on $\,\mathcal{L}$, which in turn gives
$\,g(\navp,\navp)=Q(\vp/\y)$. Thus, $\,\y\,$ is the unique real pole of
$\,|\navp|^2$ treated as a rational function of the variable $\,\vp\,$ (cf.\
Remark~\ref{apole}), while $\,\th=\,\mathrm{min}\hskip3pt\vp/\y\,$ and
$\,\ha=\,\mathrm{max}\hskip3pt\vp/\y\hs$, since $\,\vp/\y=\ps$.

A purely local definition of $\,\y\,$ is also possible, even in a much more
general situation; see \cite[Lemma~12.5 (and Corollary~9.3)]{local}. Finally,
if $\,m\,$ is fixed, $\,\th\,$ and $\,\ha\,$ depend only on the homothety
class of the Riemannian metric on the $\,2$-sphere obtained by restricting
$\,g\,$ to some, or any, fibre of the $\,\bbCP^1$ bundle $\,M$. (We will not
use this fact, which follows since $\,\vp\,$ is a Killing potential, while the
fibre geometry determines the corresponding Killing field uniquely up to a
factor.)
\end{rem}

\section{Rationality conditions}\label{rati}
\setcounter{equation}{0}
\setcounter{equation}{0}
We will use a result of Kobayashi and Ochiai \cite{kobayashi-ochiai}, as
quoted in subsection 9.124 of \cite{besse}: if $\,N\hs$ is a compact complex
manifold with $\,\dimc N=m-1\,$ such that $\,c_1(N)\in H^2(N,\bbZ)\,$ is
positive and divisible by an integer $\,\el\ge1$, then $\,\el\le m\hs$, with
equality only if $\,N\,$ is biholomorphic to $\,\bbCP^{m-1}\hskip-1.2pt$.

For $\hs m,Q,I\nh,\th,\y,\ve,\ax\hs$ as in (\ref{cea}) -- (\ref{nho}), let
$\hs\px\in\bbR$ and $\hs\dv\nh\in\nh\{-1,0,1\}$ be given by
\begin{equation}
\quad\mathrm{i)}\hskip12ptm\px\,=\,\dot Q(\th)/\ax\hskip11pt
\mathrm{(only\ if}\hskip4pt\ax\ne0\mathrm{),\hskip14ptii)}\hskip12pt\dv
=\,\sgn\,\kx\hs,\hskip9pt\mathrm{where}\hskip9pt\kx=\ve\hs m\ax/\y\hs,
\label{pde}
\end{equation}
$\hs\sgn\,$ being the usual signum function with $\,\,\sgn\,0=0\,$ and
$\,\hs\sgn\,\xh=\xh/|\xh|\,$ for $\,\xh\in\bbR\smallsetminus\{0\}$.

The invariant $\,\dv\,$ in (\ref{pde}.ii) depends just on the original
$\,m,Q,I\nh$, and not on $\,\y\,$ or $\,\ve$. Both $\,\px$, defined only if
$\,\ax\ne0$, and $\,\dv\,$ remain unaffected when $\,Q\,$ is multiplied by a
positive constant. In fact, $\,I\hs$ and $\,Q\,$ determine
$\,\hs\sgn(\ve\hs\y)\,$ via (a), while rescaling $\,Q\,$ leads to
multiplication of $\,\ax\,$ in (\ref{qef}) and $\,\dot Q\,$ by the same
positive factor.

Unlike $\,\y,\ve,\ah\,$ in (\ref{cea}), the objects (\ref{nho}) need not
exist. Moreover, whether they exist or not depends just on $\,m,Q,I\nh$, and
not on how we chose $\,\y,\ve,\ah$. Namely, let $\,\px,\dv\,$ be
determined by $\,m,Q\,$ and $\,I\hs$ as in (\ref{pde}). Then (\ref{nho}) holds
for some $\,N,h,\mathcal{L}\,$ and a $\,\hs\mathrm{U}\hs(1)\,$ connection in
$\,\mathcal{L}\hs$, if and only if
\begin{equation}
\begin{array}{l}
\mathrm{either\ }\,\dv=1\,\mathrm{\ and\ }\,\px=\ka/\el\,\mathrm{\ for\ some\ }
\,\ka\in\bbZ\,\mathrm{\ and\ }\,\el\in\{1,\dots,m\},\\
\mathrm{or\ }\,\,\dv=0\hs,\hskip8pt\mathrm{ or,\hskip8ptfinally,\ }
\,\,\dv=-\hs1\,\,\mathrm{\ and\ }\,\,\px\,\,\mathrm{\ is\ rational.}
\end{array}
\label{edo}
\end{equation}
In fact, given (\ref{nho}) with $\,\kx=\ax=0$, we have $\,\dv=0\,$
by (\ref{pde}.ii), and (\ref{edo}) follows. Also, if $\,\kx\ne0\,$ in
(\ref{nho}) (so that $\,\ax\ne0$), then (\ref{pde}) and the Kobayashi-Ochiai theorem
mentioned above give (\ref{edo}) with $\,\dv=\pm\hs1\,$ (cf.\ Remark~\ref{intgr}
below).

Conversely, let (\ref{edo}) (and (\ref{cea})) be satisfied. If $\,\ax=0$, (\ref{nho}) is easily realized by choosing $\,(N,h)\,$ to be a compact
Ricci-flat K\"ahler manifold whose K\"ahler class equals $\,-\hs\ve\pi/\ah\,$
times an integral class (for instance, a suitable flat complex torus). We then
select $\,\mathcal{L}\,$ so that the latter class is $\,c_1(\mathcal{L})$, and hence
$\,\varOmega\,=\hs-\hs2\hs\ve\hskip.2pt\ah\hs\omh$ is the curvature form of
some $\,\hs\mathrm{U}\hs(1)\,$ connection in $\,\mathcal{L}$.

Finally, let us assume (\ref{edo}) with $\,\ax\ne0$, that is, $\,\dv=\pm\hs1$.
Thus, $\,\px=\ka/\el\,$ with relatively prime integers $\,\ka\,$ and
$\,\el\ge1$. For any integer $\,s\ge1$, let $\,N\hs$ be the {\it Fermat
hypersurface\/} of degree $\,s\,$ in $\,\bbCP^m\nh$, given by
$\,z_0^s+z_1^s+\ldots+z_m^s=0\,$ in homogeneous coordinates $\,z_0,\dots,z_m$.
The adjunction formula \cite[p.\ 147]{griffiths-harris} implies that
$\,c_1(N)=(m+1-s)\hs[e]\,$ in $\,H^2(N,\bbR)$, where $\,[e]\,$ is the
restriction to $\,N\hs$ of the positive generator of $\,H^2(\bbCP^m\nh,\bbR)$.
This has three consequences. First, if $\,s>m+1$, or $\,s\in\{1,\dots,m\}$,
then $\,c_1(N)<0\,$ or, respectively, $\,c_1(N)>0$, and so a K\"ahler-Einstein
metric $\,h\,$ on $\,N\hs$ exists in view of the Aubin and Yau solution to
Calabi's conjecture \cite{aubin}, \cite{yau} or, respectively, by a recent
result of Tian \cite{tian}. Next, if $\,s=m+1-\dv\el\ge1$, then
$\,\kx\hs,\hs\dv,\hs m+1-s,\hs c_1(N)\,$ and the Ricci form $\,\roh$ of
$\,h\,$ all have the same sign; hence, rescaling $\,h$, we can always ensure
that $\,\roh=\kx\hskip1.2pt\omh\nh$. Third, if $\,\mathcal{E}\hs$ is the
restriction to $\,N\hs$ of the dual of the tautological bundle over
$\,\bbCP^m$ and $\,\mathcal{L}=\mathcal{E}^{\otimes\ka}\nh$, then
$\,c_1(\mathcal{L})=\px\hs c_1(N)\,$ in $\,H^2(N,\bbR)$, as both sides equal
$\,\ka\hs[e]$. Now (\ref{nho}) follows, with a connection in $\,\mathcal{L}\,$
chosen as in Remark~\ref{curfo} below.
\begin{rem}\label{intgr}
We clearly have $\,\varOmega=\px\hs\roh$ and
$\,c_1(\mathcal{L})=\px\hs c_1(N)\,$ in $\,H^2(N,\bbR)$, for $\,\px\,$ as in (\ref{pde}), whenever (\ref{cea}) -- (\ref{nho}) hold with $\,\ax\ne0$. Since
both Chern classes are integral, $\,\px\,$ must be a rational number and its
irreducible denominator divides $\,c_1(N)\,$ in $\,H^2(N,\bbZ)$.
\end{rem}

\begin{rem}\label{curfo}
If $\,\mathcal{L}\,$ is a complex line bundle over a compact
K\"ahler manifold $\,(N,h)\,$ such that $\,c_1(\mathcal{L})=\px\hs c_1(N)\,$ in
$\,H^2(N,\bbR)\,$ for some (necessarily rational) number $\,\px$, then
$\,\varOmega=\px\hs\roh$ is the curvature form of some
$\,\hs\mathrm{U}\hs(1)\,$ connection in $\,\mathcal{L}$. In fact,
$\,2\pi\hs\varOmega\,$ represents $\,c_1(\mathcal{L})\,$ in $\,H^2(N,\bbR)$, since
$\,\roh$ represents $\,c_1(N)$.
\end{rem}

\begin{rem}\label{speca}
The importance of the objects $\,N,h,\mathcal{L}\,$ in (\ref{nho}) is due to their role as building blocks for the construction in
\S\ref{main}, leading to quadruples $\,\mgmt\,$ with (\ref{zon}) or (\ref{zto}). As the
case where such $\,N,h,\mathcal{L}\,$ exist is completely characterized by (\ref{edo}), the next natural question concerns the extent of freedom in
choosing $\,N,h,\mathcal{L}$, for any fixed $\,m,Q,I\hs$ having the properties
listed immediately before (\ref{cea}), along with (\ref{edo}). When
$\,\ax\ne0\,$ in (\ref{qef}), we can make the following comments.

In view of Remark~\ref{intgr} and \cite[Remark~2.4]{potentials}, once the
K\"ahler
manifold $\,(N,h)\,$ is selected, the choices of $\,\mathcal{L}\,$ become quite
limited: up to tensoring by holomorphic line bundles with flat
$\,\hs\mathrm{U}\hs(1)\,$ connections, $\,\mathcal{L}=\mathcal{E}^{\otimes p}$ for the
anticanonical bundle $\,\mathcal{E}=[TN]^{\wedge(m-1)}$ of $\,N\nh$. (A connection
required in (\ref{nho}) exists by Remark~\ref{curfo}.)

As for selecting $\,(N,h)$, there are two interesting special cases. First,
(\ref{edo}) obviously holds if the rational number $\,\px\,$ (cf.\
Remark~\ref{intgr}) corresponding to the given $\,m,Q,I\hs$ with $\,\ax\ne0\,$
is an integer; our $\,(N,h)\,$ then can be {\it any\/} compact
K\"ahler-Einstein manifold $\,(N,h)\,$ with $\,\dimc N=m-1\,$ that has
the correct value of the Einstein constant $\,\kx$, with
$\,\mathcal{L}\,$ as in the last paragraph.

An opposite extreme occurs, when (\ref{pde}), for our $\,m,Q,I\hs$ with
$\,\ax\ne0$, gives $\,\dv=1\,$ and $\,\px=\ka/m\,$ for an integer $\,\ka\,$
such that $\,\ka,m\,$ are relatively prime (which clearly implies
(\ref{edo})). The objects realizing (\ref{nho}) then are {\it essentially
unique\hs}: $\,(N,h)\,$ must be biholomorphically isometric to $\,\bbCP^{m-1}$
with a constant multiple of the Fubini-Study metric, in such a way that
$\,\mathcal{L}\,$ becomes the $\,\ka\hs$th tensor power of the dual
tautological bundle. This is due to the equality clause in the
Kobayashi-Ochiai theorem (see the beginning of this section), since a
holomorphic line bundle $\,\mathcal{L}\,$ over $\,N=\bbCP^{m-1}$ is uniquely
determined by $\,c_1(\mathcal{L})\in H^2(N,\bbR)$, while a K\"ahler-Einstein
metric on $\,\bbCP^{m-1}$ is essentially unique
(\cite[pp.\ 144--145]{griffiths-harris}, and \cite{bando-mabuchi}).
\end{rem}

\section{Some functional relations}\label{func}
\setcounter{equation}{0}
Throughout this section, $\,F,E\,$ and $\,\vs\,$ are the functions with
(\ref{fet}) -- (\ref{esg}) for a fixed integer $\,m\ge2$. By (\ref{fet}), the
derivative $\,\dot F=\,dF/d\ps\,$ is given by
\begin{equation}
\mathrm{a)}\hskip9pt
\dot F(\ps)\,=\,{\displaystyle{\ps^{2m-2}\over(\ps-1)^{m+1}}}\,
\varLambda(\ps)\,,\hskip20pt\mathrm{where}\hskip14pt
\mathrm{b)}\hskip9pt
\varLambda(\ps)\,=\,m\ps^2\,-\,2(2m-1)(\ps-1)\,>\,0\hs.
\label{dtf}
\end{equation}
The dependence of $\,F,E,\vs\,$ on $\,m\,$ will usually be suppressed in our
notation. Right now, however, we make it explicit by writing $\,F_m,E_m,\vs_m$
rather than $\,F,E,\vs$. For $\,m\ge2\,$ one then has, as in
\cite[the paragraph preceding (21.3)]{local},
\begin{equation}
\mathrm{i)}\hskip7pt
F_m(\ps)\,=\,{\displaystyle{\ps^2\over \ps-1}\,F_{m-1}(\ps)}\hs,\hskip13pt
\mathrm{ii)}\hskip7pt
E_m(\ps)\,=\,{\displaystyle{\ps^2\over \ps-1}\,E_{m-1}(\ps)\,
-\,{1\over m}{2m-2\choose m-1}}\hs,\hskip11ptE_1(\ps)\,=\,\ps-1\,.
\label{fmt}
\end{equation}
Here (\ref{fmt}.i) follows from (\ref{fet}) and (\ref{fmt}.ii) is obtained by
expanding the difference of the two sides into powers of $\,\ps\,$ via
(\ref{fet}) -- (\ref{esg}). Thus,
\begin{equation}
{E\over F}\,=\,{E_{m-1}\over F_{m-1}}\,
-\,{2m-2\choose m-1}{1\over mF}\,,\qquad\mathrm{where}\quad
E=E_m\hskip5pt\mathrm{and}\hskip7ptF=F_m\hs.\label{efm}
\end{equation}
Consequently, for every $\,\ps\in\bbR\smallsetminus\{0,1,2\}$, induction on
$\,m\ge2\,$ gives
\begin{equation}
{d\over d\ps}\,[E/F]\,=\,
-{2m\choose m}{(\ps-1)^m\over(\ps-2)^2\ps^{2m}}\,.\label{def}
\end{equation}
Specifically, the inductive step comes from (\ref{efm}), where one
differentiates $\,1/F\,$ using (\ref{dtf}.a), then replaces $\,F\,$ with the
expression in (\ref{fet}), and uses (\ref{dtf}.b). Thus,
\begin{equation}
\begin{array}{rl}\arraycolsep9pt
\mathrm{i)}&
\ps(\ps-1)(\ps-2)\hs\dot F(\ps)\,=\,\varLambda(\ps)\hs F(\ps)\qquad
\mathrm{for\ all}\hskip6pt\ps\in\bbR\smallsetminus\{1\}\hs,\\
\mathrm{ii)}&
\ps(\ps-1)(\ps-2)\dot E(\ps)\,=\,\varLambda(\ps)\hs E(\ps)\,
-\,2(2m-1)(\ps-1)\hs\vs(0)\hskip7pt\mathrm{for}\hskip5pt\ps\in\bbR\hs,
\hskip10pt\end{array}
\label{dte}
\end{equation}
for $\,\varLambda\,$ as in (\ref{dtf}.b). Namely, (\ref{dte}.i) is obvious
from (\ref{dtf}.a) and (\ref{fet}), while (\ref{dte}.ii) follows if one
rewrites (\ref{def}) multiplied by $\,F\,$ using the quotient rule for
derivatives, (\ref{dte}.i), and the definitions of $\,F,\vs\,$ in
(\ref{fet}), (\ref{esg}). Also, by (\ref{fet}) -- (\ref{esg}),
\begin{equation}
\mathrm{i)}\hskip9pt
{\displaystyle\vs(0)\,=\,\dot\vs(0)\,=\,{1\over m}{2m-2\choose m-1}},\hskip12pt
\mathrm{ii)}\hskip9ptE(0)\,=\,-\,\vs(0)\,,\hskip12pt
\mathrm{iii)}\hskip9pt\dot E(0)\,=\,0\,.
\label{sto}
\end{equation}
\begin{rem}\label{roots}
The coefficients of the polynomial $\,\vs\,$ given by (\ref{esg}), for any fixed integer $\,m\ge1$, are all positive, and so
$\,\vs(\ps)>0\,$ whenever $\,\ps\ge0$. Therefore, all real roots of $\,\vs\,$
are negative. The same applies to the derivative $\,\dot\vs\,$ when $\,m\ge2$.
Thus, by (\ref{esg}), $\,E(\th)<0<E(\ha)\,$ whenever $\,0\le\th<1<\ha$.
\end{rem}

If $\,Q\in\jv$, with $\,\jv\,$ as in (\ref{spv}) for a fixed integer $\,m\ge2$,
and $\,\ax,\bx,\cx\in\bbR\,$ represent $\,Q\,$ in (\ref{qef}), then
$\,\ps(\ps-1)(\ps-2)\dot Q(\ps)=[\hs\ps(\ps-2)+\varLambda(\ps)\hs]\hs Q(\ps)
-(\ps-1)\hs\ax\varLambda(\ps)-2(2m-1)(\ps-1)^2\bx\vs(0)$, as one sees using
(\ref{qef}), (\ref{dte}) and, again, (\ref{qef}). Thus,
$\,\th(\th-2)\dot Q(\th)=-\hs\ax\varLambda(\th)-2(2m-1)(\th-1)\hs\bx\vs(0)\,$
whenever $\,\th\in\bbR\smallsetminus\{1\}\,$ and $\,Q(\th)=0$. Hence, if
$\,\ax\ne0$, (\ref{dtf}.b) yields, for $\,\px\,$ as in (\ref{pde}.i),
\begin{equation}
\px\,=\,{2(\th-1)\lx-\th^2\over(\th-2)\th}
\quad\mathrm{whenever}\hskip6pt\th\in\bbR\smallsetminus\{0,1,2\}\hskip6pt
\mathrm{and}\hskip6ptQ(\th)=0\hs,\label{pum}
\end{equation}
where, for any given $\,\th\in\bbR\smallsetminus\{1\}\,$ with $\,Q(\th)=0$, we
define $\,\lx\in\bbR\,$ by
\begin{equation}
\lx\,=\,(2\,-\,1/m)\hs(1\,-\,\bx\vs(0)/\ax)\,.\label{lba}
\end{equation}

\section{Monotonicity intervals}\label{mono}
\setcounter{equation}{0}
Let $\,F,E\,$ be defined by (\ref{fet}) with an integer $\,m\ge2$. By
(\ref{fet}) -- (\ref{esg}), (\ref{dtf}) and (\ref{def}), $\,F\,$ (or, $\,E/F$)
has a nonzero derivative everywhere in $\,\bbR\smallsetminus\{0,1\}\,$ (or,
respectively, in
$\,\bbR\smallsetminus\{0,1,2\}$). One also easily sees that the rational
functions $\,F\,$ and $\,E/F\,$ of the real variable $\,\ps\,$ have the
values/limits at $\,\pm\infty,0,1,2\,$ listed below. One-sided limits, if
different, are separated by vertical arrows indicating the direction of the
jumps. The slanted arrows show which kind of monotonicity the given function
has on each of the four intervals forming $\,\bbR\smallsetminus\{0,1,2\}$.
\begin{equation}
\begin{array}{lccccccccc}\arraycolsep5pt
\mathrm{value\ or\ limit\ at\hskip-3.2pt:}&-\infty&&0&&1&&2&&+\infty\\
[1pt]
\hline
&&&&&&&&&\\[-9.7pt]
\mathrm{for}\hskip4.5ptF\hskip4pt\mathrm{(}m\hskip4pt\mathrm{even)\hskip-3.2pt:
}&+\infty&\sea&0&\sea&-\infty&\nea&0&\nea&+\infty\\
\mathrm{for}\hskip4.5ptF\hskip4pt\mathrm{(}m\hskip4pt\mathrm{odd)\hskip-3.2pt:}
&-\infty&\nea&0&\nea&+\infty\vd{-\infty}&\nea&0&\nea&+\infty\\
\mathrm{for}\hskip4.5ptE/F\hskip4pt\mathrm{(}m\hskip4pt\mathrm{even)
\hskip-3.2pt:}\hskip0pt&1&\sea&-\infty\vu{+\infty}&\sea&0&\sea&-\infty\vu
{+\infty}&\sea&1\\
\mathrm{for}\hskip4.5ptE/F\hskip4pt\mathrm{(}m\hskip4pt\mathrm{odd)\hskip-3.2pt
:}\hskip0pt&1&\nea&+\infty\vd{-\infty}&\nea&0&\sea&-\infty\vu{+\infty}&\sea&1
\end{array}
\label{mon}
\end{equation}

\section{Some inequalities}\label{smiq}
\setcounter{equation}{0}
For any fixed integer $\,m\ge2\,$ we have, with $\,F,E\,$ as in (\ref{fet}),
\begin{equation}
\begin{array}{rl}\arraycolsep9pt
\mathrm{a)}&E\,<\,F\hskip5pt\mathrm{on}\hskip5pt(-\infty,0\hs]\hs,\hskip10pt
\mathrm{b)}\hskip5ptE\,>\,0\hskip6pt\mathrm{on}\hskip5pt(1,\infty)\hs,
\hskip10pt
\mathrm{c)}\hskip5ptE\,>\,F\hskip5pt\mathrm{on}\hskip5pt(1,\infty)\hs,\\
\mathrm{d)}&\mathrm{If}\hskip5ptm\hskip5pt\mathrm{is\ odd,}
\hskip5ptE\,<\,0\hskip5pt\mathrm{on}\hskip5pt(-\infty,1)\hs.\end{array}
\label{elf}
\end{equation}
In fact, (\ref{elf}.b) is clear from Remark~\ref{roots}. Next, (\ref{elf}.a) and (\ref{elf}.c) follow since (\ref{mon}) gives $\,E/F<1\,$ and $\,F>0\,$
on $\,(-\infty,0)\,$ for even $\,m\hs$, and $\,E/F>1\,$ and $\,F<0\,$ on
$\,(-\infty,0)\,$ for odd $\,m\hs$, while, for all $\,m\hs$, it yields $\,E/F<0\,$
and $\,F<0\,$ on $\,(1,2)\,$ (so that $\,F<0<E\,$ there), as well as
$\,E/F>1\,$ and $\,F>0\,$ on $\,(2,\infty)$. (Also, $\,E(0)<0<E(2)\,$ and
$\,F(0)=F(2)=0\,$ by Remark~\ref{roots} and (\ref{fet}).) Finally, for odd $\,m\hs$, we
have $\,E<0\,$ both on $\,(-\infty,0\hs]\,$ (from the inequalities just
listed) and on $\,(0,1)\,$ (from Remark~\ref{roots}), which implies (\ref{elf}.d).

\section{A linear-independence property}\label{indp}
\setcounter{equation}{0}
With $\,E(\ps)\,$ as in (\ref{fet}) for a fixed integer $\,m\ge2$, we define a
polynomial function $\,\,\mathbf{r}:\bbR\to\rtr$ by
\begin{equation}
\mathbf{r}(\ps)\,=\,(\ps-1)^m\hskip1pt\mathbf{i}\,\,
+\,(\ps-1)^mE(\ps)\hskip1pt\mathbf{j}\,\,
+\,(\ps-2)\ps^{2m-1}\hskip1pt\mathbf{k}\,,\label{blr}
\end{equation}
where $\,\hs\mathbf{i}\hs,\,\mathbf{j}\hs,\,\mathbf{k}\hs\,$ form the
standard basis of $\,\rtr\hskip-1.2pt$. Note that, as $\,\dot E(0)=0\,$ by
(\ref{sto}.iii),
\begin{equation}
\mathbf{r}(0)\,=-\,\dot{\mathbf{r}}(0)/m\,
=\,(-1)^m\,\lj\,\mathbf{i}\,\,+\,E(0)\hskip1pt\mathbf{j}\,\rj\,,\hskip15pt
\mathbf{r}(1)\,=\,-\,\mathbf{k}\,.
\label{rzr}
\end{equation}
\begin{prop}\label{indep}For any integer\/ $\,m\ge2\,$ and\/
$\,\th,\ha\in\bbR\,$ with $\,\th\ne\ha$, the vectors\/
$\,\,\mathbf{r}(\th),\,\mathbf{r}(\ha)\,$ defined as in {\rm(\ref{blr})} are
linearly independent.
\end{prop}
\begin{proof}Writing $\,\,\mathbf{r}=(\xh,\yh,\zh)\,$ we have
$\,\xh(\ps)\ne0\,$ for $\,\ps\ne1$, while $\,\,\mathbf{r}(1)=(0,0,-1)\,$ by
(\ref{rzr}). Therefore, our assertion follows if one of $\,\th,\ha\,$ equals
$\,1$. On the other hand, $\,\zh(\ps)/\xh(\ps)=F(\ps)\,$ for $\,\ps\ne1$, with
$\,F(\ps)\,$ as in (\ref{fet}). This proves our assertion in the case where
$\,\th\,$ and $\,\ha\,$ are both greater than $\,1\,$ or both less than
$\,1\,$ since, according to (\ref{mon}), $\,F\,$ is injective both on
$\,(-\infty,1)\,$ and on $\,(1,\infty)$.

Therefore, switching $\,\th\,$ and $\,\ha\,$ if necessary, we may assume that
$\,\th<1<\ha$. Contrary to our assertion, let $\,\,\mathbf{r}(\th)\,$ and
$\,\,\mathbf{r}(\ha)\,$ be linearly {\it dependent}. Now (\ref{blr}) gives
$\,(\th-1)^{-m}\hs\mathbf{r}(\th)=(\ha-1)^{-m}\hs\mathbf{r}(\ha)$, and
hence, by (\ref{fet}), $\,F(\th)=F(\ha)\,$ and $\,E(\th)=E(\ha)$. However, the
last equality contradicts the relation $\,E(\th)<E(\ha)$, which is immediate
both when $\,0\le\th<1<\ha\,$ (see Remark~\ref{roots}) and in the case where
$\,\th<0<1<\ha\,$ (since (\ref{elf}.a,\hs c) then yield
$\,E(\th)<F(\th)=F(\ha)<E(\ha)$). This contradiction completes the
proof.
\end{proof}

\section{Condition \hs{\rm(\ref{qdq}.b)}\hs\ alone}\label{lone}
\setcounter{equation}{0}
For $\,Q\in\jv\,$ for $\,\jv\,$ as in (\ref{spv}) with a fixed integer
$\,m\ge2$, we have
\begin{equation}
\begin{array}{rl}\arraycolsep9pt
\mathrm{i)}&
(\ps-1)^{m-1}\hs Q(\ps)\,
=\,\,\mathbf{p}\cdot\mathbf{r}(\ps)\,,\quad\mathrm{and}\\
\mathrm{ii)}&
(\ps-1)^m\hs\dot Q(\ps)\,=\,\,(\ps-1)\,\mathbf{p}\cdot
\dot{\mathbf{r}}(\ps)\,-\,(m-1)\,\mathbf{p}\cdot\mathbf{r}(\ps)\,,
\quad\mathrm{with}\hskip15pt\\
\mathrm{iii)}&
\mathbf{p}\,\,=\,\ax\hskip1pt\mathbf{i}\,\,+\,\bx\hskip1pt\mathbf{j}\,\,
+\,\cx\hskip1pt\mathbf{k}\,\in\,\rtr,\end{array}
\label{qpr}
\end{equation}
where $\,\hs\mathbf{r}(\ps),\,\mathbf{i}\hs,\,\mathbf{j}\hs,
\,\mathbf{k}\,\hs$ are as in (\ref{blr}), while $\,\hs\cdot\hs\,$ denotes the
inner product of $\,\rtr\nh$, and $\,\ax,\bx,\cx\,$ correspond to $\,Q\,$ via (\ref{qef}). (In fact, (\ref{qef}) gives (\ref{qpr}.i), and differentiation then
leads to (\ref{qpr}.ii).) Thus, given such $\,Q,\hs\jv\,$ and $\,m\hs$, we have
\begin{equation}
\begin{array}{rl}\arraycolsep9pt
\mathrm{a)}&
\dot Q(0)\,=\,0\quad\mathrm{whenever}\quad Q\in\jv
\quad\mathrm{and}\quad Q(0)\,=\,0\hs,
\hskip43pt\\
\mathrm{b)}&
\mathrm{If}\hskip7ptQ\in\jv\hs\mathrm{, \ then}\hskip8ptQ(\th)\,=\,0
\hskip4pt\mathrm{\ if \ and \ only \ if}\hskip8pt\mathbf{p}\cdot
\mathbf{r}(\th)\,=\,0\hs,
\end{array}
\label{dqo}
\end{equation}
for any $\,\th\in\bbR\hs$, where $\,\,\mathbf{p}\in\rtr$ corresponds to
$\,Q\,$ as in (\ref{qpr}.iii). Namely, from (\ref{rzr}) we have
$\,\,\mathbf{p}\cdot\mathbf{r}(0)=-\,\mathbf{p}\cdot
\dot{\mathbf{r}}(0)/m=(-1)^m\,\lj\hs\ax+\bx\hs E(0)\hs\rj\,$ which, combined
with (\ref{qpr}), gives $\,\dot Q(0)=-\hs Q(0)=\ax+\bx\hs E(0)$, and hence
(\ref{dqo}.a). Next, when $\,\th\ne1$, (\ref{dqo}.b) is clear from
(\ref{qpr}.i), while, if $\,\th=1$, condition $\,Q(1)=0\,$ is, by (\ref{qef})
with (\ref{fet}), equivalent to $\,\cx=0$, that is (cf.\ (\ref{qpr}.iii) and
(\ref{rzr})), to $\,0=\,\mathbf{p}\cdot\mathbf{k}\,
=\,-\,\mathbf{p}\cdot\mathbf{r}(1)$.
\begin{lem}\label{onedi}Given an integer\/ $\,m\ge2\,$ and a nontrivial closed
interval\/ $\,I\subset\bbR$, those $\,Q\in\jv\,$ which vanish at both
endpoints of\/ $\,I\hs$ form a one\diml\ vector subspace of the space\/
$\,\jv\,$ with {\rm(\ref{spv})}. Thus, $\,Q\in\jv\smallsetminus\{0\}\,$
satisfying {\rm(\ref{qdq}.b)} on $\,I$ exists, for any such\/ $\,m,I$, and
is unique up to a nonzero constant factor.

Explicitly, up to a factor,
$\,(\ps-1)^{m-1}\hs Q(\ps)=\lj\hskip1pt\mathbf{r}(\th)\times\hskip.3pt
\mathbf{r}(\ha)\hs\rj\cdot\,\mathbf{r}(\ps)$, where\/ $\,\th,\ha\,$ denote
the endpoints of\/ $\,I\hs$ and\/ $\,\,\mathbf{r}\hs\,$ is given by
{\rm(\ref{blr})}, while\/ $\,\cdot\,$ and\/ $\,\times\,$ denote the inner
product and vector product in $\,\rtr$.
\end{lem}
In fact, up to a factor, $\,\,\mathbf{p}\,\,$ in (\ref{qpr}.i) equals
$\,\,\mathbf{r}(\th)\times\hskip.3pt\mathbf{r}(\ha)\ne\,\mathbf{0}$,
since $\,\,\mathbf{r}(\th),\,\mathbf{r}(\ha)\,$ are linearly independent
(Proposition~\ref{indep}) and orthogonal to $\,\hs\mathbf{p}\hs\,$ (by (\ref{dqo}.b)).

\section{A determinant formula for $\,Q(\ps)$}\label{detr}
\setcounter{equation}{0}
Given an integer $\,m\ge2\,$ and a nontrivial closed interval $\,I$, let
$\,Q\in\jv\smallsetminus\{0\}$, with $\,\jv\,$ as in (\ref{spv}), be a function
satisfying condition (\ref{qdq}.b) on $\,I\nh$. Such $\,Q$, which exists and is
unique up to a factor, for any $\,m\,$ and $\,I$, is explicitly described in
Lemma~\ref{onedi}. Let us now also assume that $\,\th\ne1\ne\ha$, where
$\,\th\,$ and $\,\ha\,$ are the endpoints of $\,I$. With $\,F,E\,$ as in (\ref{fet}) for our fixed $\,m\ge2$, (\ref{blr}) and (\ref{fet}) give
$\,\,\mathbf{r}(\ps)=(\ps-1)^m\hs\mathbf{w}(\ps)\,$ for $\,\ps\ne1$, with
$\,\,\mathbf{w}(\ps)=\,\mathbf{i}\,+E(\ps)\hskip1pt\mathbf{j}\,
+F(\ps)\hskip1pt\mathbf{k}$. Thus,
$\,(\ps-1)^{-1}Q(\ps)=(\th-1)^m(\ha-1)^m\hs[\mathbf{w}(\th)\times\hskip.3pt
\mathbf{w}(\ha)]\cdot\,\mathbf{w}(\ps)\,$ by Lemma~\ref{onedi}, that is, up to
another factor,
\begin{equation}
Q(\ps)\,=\,(\ps-1)\hh(\ps)\hskip17pt\mathrm{for\ all}\hskip10pt
\ps\in\bbR\smallsetminus\{1\}\hs,\hskip22pt\mathrm{where}\quad\hh(\ps)\,
=\,\,\det\left[\begin{array}{ccc}
1&E(\th)&F(\th)\\1&E(\ha)&F(\ha)\\1&E(\ps)&F(\ps)
\end{array}\right]\nh.\label{qth}
\end{equation}
Writing $\,\hh,\,F,\,E\,$ for $\,\hh(\ps),\,F(\ps),\,E(\ps)$, and
$\,(\,\,)\dot{\,}=\,d/d\ps$, we have, by (\ref{qth}),
\begin{equation}
\begin{array}{rl}\arraycolsep9pt
\mathrm{i)}&
\hh\,=\,E_0F_1\,-\,F_0E_1\,+\,(F_0\,-\,F_1)E\,+\,(E_1\,-\,E_0)F\,,\\
\mathrm{ii)}&
\dot\hh\,=\,(F_0\,-\,F_1)\dot E\,+\,(E_1\,-\,E_0)\dot F\,,\qquad
\mathrm{where}\\
\mathrm{iii)}&
E_0\,=\,E(\th)\,,\quad E_1\,=\,E(\ha)\,,\quad
F_0\,=\,F(\th)\,,\quad F_1\,=\,F(\ha)\,.\hskip20pt\end{array}
\label{hef}
\end{equation}
\begin{lem}\label{abceq}Let\/ $\,\jv\,$ be the space {\rm(\ref{spv})} for a
fixed integer\/ $\,m\ge2$, and let a nontrivial closed interval\/
$\,I\subset\bbR\,$ have the endpoints\/ $\,\th,\ha\,$ with $\,\th\ne1\ne\ha$.
If a function $\,Q\in\jv\,$ satisfies condition {\rm(\ref{qdq}.b)} on $\,I$,
that is, $\,Q(\th)=Q(\ha)=0$, while $\,(\ax,\bx,\cx)\in\rtr$ corresponds to
$\,Q\,$ as in {\rm(\ref{qef})}, then, up to a nonzero overall factor,
\begin{equation}
(\ax,\bx,\cx)\,=\,(E(\th)F(\ha)\,-\,F(\th)E(\ha)\hs,\,F(\th)-F(\ha)\hs,\,
E(\ha)-E(\th))\,.\label{abc}
\end{equation}
\end{lem}
This is clear from (\ref{qef}) and (\ref{hef}.i\hs,\hs iii), as
$\,Q(\ps)=(\ps-1)\hh(\ps)\,$ by (\ref{qth}).

\section{A convexity lemma}\label{cnvx}
\setcounter{equation}{0}
According to \S\ref{mono}, the variable $\,\ps\,$ may, on suitable intervals,
be diffeomorphically replaced with $\,F$. As shown next, this makes $\,E\,$ or
$\,-\hs E\,$ a convex function of $\,F$.
\begin{lem}\label{cnvex}Let\/ $\,F,E,\varLambda\,$ be as in {\rm(\ref{fet})}
and\/ {\rm(\ref{dtf})} for an integer\/ $\,m\ge2$. Then
\begin{equation}
{d\over d\ps}\,[\dot E/\dot F\hs]\,=\,{2m\choose m}
{m(\ps-1)^m\over[\varLambda(\ps)]^2\hskip1pt\ps^{2m-2}}\hskip11pt\mathrm{at\
every\ real}\hskip9pt\ps\ne0\hs,\hskip6pt\mathrm{with}\hskip6pt
(\,\,)\dot{\,}=\,d/d\ps\hs.\label{cnv}
\end{equation}
\end{lem}
In fact, let $\,\gy={2m\choose m}$. Writing $\,F\,$ for $\,F(\ps)$, etc., let
us differentiate (\ref{def}) multiplied by $\,F^2$ and then multiply the
result by $\,\dot F\ns/F$, obtaining
$\,\dot F\ddot E-E\dot F\ddot F\ns/F=-\,\gy\hs F^{-1}\dot F\,
d\hs[\ps^{-1}(\ps-2)^{-1}F\hs]/d\ps$. Also, multiplying (\ref{def}) by
$\,F\ddot F\,$ we get
$\,\dot E\ddot F-E\dot F\ddot F\ns/F=-\,\gy\hs\ps^{-1}(\ps-2)^{-1}\ddot F$.
Subtracting the last two relations, we see that
$\,(\dot F\ddot E-\dot E\ddot F)/\gy\,$ coincides with
$\,L^{-1}\ddot F-F^{-1}\dot F\hs d\hs[F/L]/d\ps\,$ for $\,L=\ps(\ps-2)$,
which, for any given $\,C^2$ functions $\,F,L\,$ of the real variable $\,\ps$,
obviously equals $\,L^{-2}F\hs d\hs[L\dot F/F\hs]/d\ps\,$ wherever $\,FL\ne0$.
Hence
$\,\dot F\ddot E-\dot E\ddot F=\gy\hs L^{-2}F\hs d\hs[L\dot F\ns/F\hs]/d\ps\,$
with $\,L=\ps(\ps-2)$. Let us now divide both sides by $\,(\dot F)^2$ and
successively replace: $\,d\hs[L\dot F\ns/F\hs]/d\ps\,$ by
$\,m\ps(\ps-2)/(\ps-1)^2$ (noting that
$\,L\dot F\ns/F=\varLambda(\ps)/(\ps-1)$, cf.\ (\ref{dte}.i), and using (\ref{dtf}.b)), then $\,\dot F\,$ by the expression provided by (\ref{dte}.i),
$\,F\,$ by its description in (\ref{fet}), $\,L\,$ by $\,\ps(\ps-2)\,$ and,
finally, $\,\gy\,$ by $\,{2m\choose m}$. This yields (\ref{cnv}) (also at
$\,\ps\in\{0,2\}$, as both sides are rational functions of $\,\ps$).
\begin{rem}\label{dotef}
Let $\,F,E\,$ be as in (\ref{fet}) for an integer
$\,m\ge2$, and let $\,(\,\,)\dot{\,}=\,d/d\ps$. By (\ref{cnv}) and
(\ref{dtf}.b), $\,\dot E/\dot F\,$ then has a nonzero derivative everywhere in
$\,\bbR\smallsetminus\{0,1\}$. The values/limits of $\,\dot E/\dot F\,$ at
$\,\pm\infty,0,1\,$ are listed below. We use the notations of \S\ref{mono},
with slanted arrows indicating, again, the monotonicity type of
$\,\dot E/\dot F$.
\begin{equation}
\begin{array}{lccccccc}\arraycolsep9pt
\mathrm{value\ or\ limit\ at\hskip-3.2pt:}&-\infty&&0&&1&&+\infty\\
[1pt]
\hline
&&&&&&&\\[-9.7pt]
\mathrm{for}\hskip4.5pt\dot E/\dot F\hskip4pt\mathrm{(}m\hskip4pt
\mathrm{even)\hskip-3.2pt:}&1&\nea&+\infty\vd{-\infty}&\nea&0&\nea&1\\
\mathrm{for}\hskip4.5pt\dot E/\dot F\hskip4pt\mathrm{(}m\hskip4pt
\mathrm{odd)\hskip-3.2pt:}&1&\sea&-\infty\vu{+\infty}&\sea&0&\nea&1
\end{array}
\end{equation}
\noindent In fact, the rational function $\,\dot E/\dot F\,$ must have some
limits at $\,\pm\infty$. By l'Hospital's rule, they coincide with those of
$\,E/F\,$ in \S\ref{mono}. (Both $\,E,F\,$ have infinite limits at $\,\pm\infty$, cf.\
\S\ref{mono}.) The limits of $\,\dot E/\dot F\,$ at $\,1\,$ and $\,0\,$ are easily
found using (\ref{fet}) and (\ref{dtf}): $\,\dot F\,$ has a pole at $\,1$, and
$\,\dot E\,$ does not, while $\,\dot F\,$ has at $\,0\,$ a zero of order
$\,2m-2$, greater than the order of a zero at $\,0\,$ for the degree $\,m-1\,$
polynomial $\,\dot E$.
\end{rem}
\begin{prop}\label{btacd}Given an integer\/ $\,m\ge2$, let\/
$\,Q\in\jv\smallsetminus\{0\}$, with\/ $\,\jv\,$ as in {\rm(\ref{spv})}, satisfy
condition {\rm(\ref{qdq}.b)} on a nontrivial closed interval\/ $\,I\hs$ with\/
$\,0\notin I\hs$ and\/ $\,1\notin I$. Then $\,Q\,$ must also satisfy
conditions {\rm(\ref{qdq}.a)}, {\rm(\ref{qdq}.c)}, {\rm(\ref{qdq}.d)}.
\end{prop}
\begin{proof}Since $\,1\notin I$, (\ref{qdq}.a) follows. Let us now suppose
that all the assumptions hold, yet, contrary to our claim, one of conditions
(\ref{qdq}.c), (\ref{qdq}.d) fails. In view of (\ref{qdq}.b), the function
$\,\,\hh\,$ with (\ref{qth}) then not only vanishes at both endpoints
$\,\th,\ha\,$ of $\,I$, but, in addition, its derivative
$\,\dot\hh=\,d\hh/d\ps\,$ is zero at one of the endpoints (if (\ref{qdq}.d) fails),
or $\,\hh=0\,$ at some interior point of $\,I\hs$ (if (\ref{qdq}.c) fails; note
that, to evaluate $\,\dot\hh\,$ at an endpoint, e.g., $\,\th$, we may treat
the $\,(\ps-1)\,$ factor in (\ref{qth}) like a nonzero constant, since
$\,\th\ne1\,$ and $\,Q(\th)=0$). In either case, Rolle's theorem gives
$\,\dot\hh=0\,$ at two distinct points of $\,I$. On the other hand, by (\ref{hef}) and Lemma~\ref{abceq}, $\,\dot\hh=0\,$ at precisely those $\,\ps\,$ at
which $\,\dot E/\dot F=-\hs\cx/\bx$. Note that $\,I\hs$ is contained in one of
the intervals $\,(-\infty,0),\,(0,1),\,(1,\infty)$, and hence, according to
(\ref{mon}), $\,F\,$ is strictly monotone on $\,I$, so that $\,\bx\ne0\,$
by Lemma~\ref{abceq}; however, for a similar reason, $\,\dot E/\dot F\,$ is strictly
monotone on $\,I\hs$ (see Remark~\ref{dotef}), and so it cannot assume the value
$\,-\hs\cx/\bx\,$ twice. This contradiction completes the proof.
\end{proof}

\section{Conditions \hs{\rm(\ref{qdq}.c)}, {\rm(\ref{qdq}.d)}\hs\ with an
endpoint at $\,1$}\label{cdcd}
\setcounter{equation}{0}
The following lemma lists some obvious facts that will help us understand
which integers $\,m\ge2\,$ and nontrivial closed intervals $\,I\hs$ containing
$\,1\,$ as an endpoint have the property that a function
$\,Q\in\jv\smallsetminus\{0\}\,$ with (\ref{qdq}.b) on $\,I\hs$ also satisfies
conditions (\ref{qdq}.c), (\ref{qdq}.d). Cf.\ also Lemma~\ref{onedi}.

For most of our discussion, the symbol $\,E\,$ has stood for the function
appearing in (\ref{fet}) with a fixed integer $\,m\ge2$. The following obvious
lemma, however, is an exception, as we allow $\,E\,$ to be much more general.
\begin{lem}\label{czero}Let\/ $\,E:\bbR\to\bbR\,$ be any function. For\/
$\,\ax,\bx\in\bbR\,$ with\/ $\bx\ne0$, and $\,\ps\in\bbR\hs$, let us set
\begin{equation}
Q(\ps)\,=\,(\ps-1)\,[\hs\ax+\bx E(\ps)\hs]\,.\label{qta}
\end{equation}
\begin{enumerate}
\item If\/ $\,\th\in\bbR\smallsetminus\{1\}$, condition\/ $\,Q(\th)=0\,$ holds
if and only if\/ $\,E(\th)=-\hs\ax/\bx$.
\item Given\/ $\,\ps,\th\in\bbR\,$ with\/ $\,E(\th)=-\hs\ax/\bx$, we have
$\,Q(\ps)=0\,$ if and only if\/ $\,\ps=1\,$ or $\,E(\ps)=E(\th)$.
\item If\/ $\,E\,$ is of class\/ $\,C^1$ and\/ $\,E(1)=0$, while
$\,\th\in\bbR\smallsetminus\{1\}\,$ and\/ $\,E(\th)=-\hs\ax/\bx$, then
\begin{enumerate}
\item[a)] $\dot Q(1)=0\,$ if and only if\/ $\,E(\th)=0$. Here and in
{\rm(b),\hs(c)}, $\,(\,\,)\dot{\,}\,$ stands for\/ $\,d/d\ps$.
\item[b)] $\dot Q(\th)=0\,$ if and only if\/ $\,\dot E(\th)=0$.
\item[c)] $\dot Q(1)+\dot Q(\th)=0\,$ if and only if\/
$\,\dot\vs(\th)=0$, where $\,\vs(\ps)=(\ps-1)^{-1}E(\ps)\,$ for\/ $\,\ps\ne1$.
\end{enumerate}
\end{enumerate}
\end{lem}
\begin{rem}\label{ptone}
Let $\,\jv\,$ be the space (\ref{spv}) with a fixed
integer $\,m\ge2$. For any nontrivial closed interval $\,I$, a function
$\,Q\in\jv\smallsetminus\{0\}\,$ satisfying (\ref{qdq}.b) on $\,I\hs$ exists and is
unique up to a factor (see Lemma~\ref{onedi}). In the case where the endpoints
of $\,I\hs$ are $\,\ha=1\,$ and $\,\th\ne1$, this $\,Q\,$ is given by
(\ref{qta}) with any constants $\,\bx\ne0\,$ and $\,\ax\,$ chosen so
that $\,E(\th)=-\hs\ax/\bx$, for $\,E\,$ as in (\ref{fet}).

In fact, (\ref{qta}) implies (\ref{qef}), while, by Lemma~\ref{czero}(i), such
a choice of $\,\ax,\bx\,$ gives (\ref{qdq}.b).
\end{rem}

The next section comprises facts we need in order to apply Lemma~\ref{czero} to
$\,E(\ps)\,$ and $\,Q(\ps)\,$ given by (\ref{fet}) and (\ref{qef}). The eventual
conclusions about functions $\,Q\in\jv\,$ satisfying conditions (\ref{qdq}.b), (\ref{qdq}.c), (\ref{qdq}.d) on intervals $\,I\hs$ with an endpoint at $\,1\,$ will be
presented later; see (iii) in \S\ref{posi}.

\section{Monotonicity properties of the function $\,E\,$ with
\hs{\rm(\ref{fet})}}\label{mopr}
\setcounter{equation}{0}
Let $\,\dot E=dE/d\ps\,$ with $\,E(\ps)\,$ as in (\ref{fet}) for an integer
$\,m\ge2$. By (\ref{sto}.i) -- (\ref{sto}.ii), $\,E(0)<0\,$ and
$\,\dot E(0)=E(1)=0$. The following claims will be verified in \S\ref{prfs}:

If $\,m\,$ is even, then $\,\dot E(\ps)\ne0\,$ for every
$\,\ps\in\bbR\smallsetminus\{0\}\,$ and there exists a unique
$\,\zu\in\bbR\smallsetminus\{1\}\,$ with $\,E(\zu)=0$. This unique $\,\zu\,$
is negative. With the notations and conventions of \S\ref{mono}, we then have
\begin{equation}
\begin{array}{lccccccccc}\arraycolsep9pt
\mathrm{value\ or\ limit\ at\hskip-3.2pt:}\quad&-\infty&&\zu&&0&&1&&+\infty\\
[1pt]
\hline
&&&&&&&&&\\[-9.7pt]
\mathrm{for}\hskip4.5ptE\hskip4pt\mathrm{(}m\hskip4pt\mathrm{even)\hskip-3.2pt:}
&+\infty&\sea&0&\sea&E(0)&\nea&0&\nea&{+\infty}\\
\end{array}
\end{equation}
If $\,m\,$ is odd, there exist unique numbers
$\,\zu,\wu\in\bbR\smallsetminus\{0\}\,$ with $\,E(\zu)=E(0)\,$ and
$\,\dot E(\wu)=0$. They satisfy the relations $\,\zu<\wu<0\,$ and
$\,E(0)<E(\wu)<0$. Also, $\,\dot E\ne0\,$ everywhere in
$\,\bbR\smallsetminus\{\zu,\wu,0,1\}$, and, in the notations of \S\ref{mono},
\begin{equation}
\begin{array}{lccccccccccc}\arraycolsep9pt
\mathrm{value\ or\ limit\ at\hskip-3.2pt:}\quad&-\infty&&\zu&&\wu&&0&&1&&
+\infty\\
[1pt]
\hline
&&&&&&&&&&&\\[-9.7pt]
\mathrm{for}\hskip4.5ptE\hskip4pt\mathrm{(}m\hskip4pt\mathrm{odd)\hskip-3.2pt:}
&-\infty&\nea&E(0)&\nea&E(\wu)&\sea&E(0)&\nea&0&\nea&{+\infty}
\end{array}
\label{emo}
\end{equation}

\section{Proofs of the claims made in \S\ref{mopr}}\label{prfs}
\setcounter{equation}{0}
According to Remark~\ref{dotef}, if $\,m\,$ is even, $\,\dot E/\dot F\,$ is
positive on $\,(-\infty,0)\cup(1,\infty)\,$ and negative on $\,(0,1)$, while,
if $\,m\,$ is odd, $\,\dot E/\dot F\,$ vanishes at a unique
$\,\wu\in(-\infty,0)$, is positive on
$\,(-\infty,\wu)\cup(0,1)\cup(1,\infty)$, and is negative on $\,(\wu,0)$.
Combined with the signs of $\,\dot F\,$ on the individual intervals (cf.\ the
slanted arrows for $\,F\,$ in \S\ref{mono}), this gives the required signs of
$\,\dot E$, that is, slanted arrows for $\,E\,$ in \S\ref{mopr}. (As
(\ref{esg}) gives $\,\dot E(\ps)=\vs(\ps)+(\ps-1)\dot\vs(\ps)$, we have
$\,\dot E(1)=\vs(1)>0\,$ by Remark~\ref{roots}.)

Since $\,E\,$ is a nonconstant polynomial, its limits at $\,\pm\infty\,$ are
infinite, with the signs required in \S\ref{mopr} (see the slanted arrows). This
proves all statements except for those involving $\,\zu\,$ and the relation
$\,E(0)<E(\wu)<0\,$ for odd $\,m\hs$. However, $\,E(0)<E(\wu)\,$ as $\,E\,$ is
decreasing on $\,[\wu,0\hs]$, so that the already-established monotonicity
pattern of $\,E\,$ gives, for all $\,m\hs$, the existence and uniqueness of
$\,\zu\,$ along with $\,\zu<0\,$ ($m\,$ even) or $\,\zu<\wu\,$ ($m\,$ odd).
Finally, if $\,m\,$ is odd, $\,E(\wu)<0\,$ by (\ref{elf}.d).

\section{Conditions equivalent to \hs{\rm(\ref{qdq}.c)} --
{\rm(\ref{qdq}.d)}}\label{cdeq}
\setcounter{equation}{0}
\begin{lem}\label{hdhdf}Given an integer\/ $\,m\ge2$, let the functions\/
$\,E,F,\hh\,$ be given by {\rm(\ref{fet})} and {\rm(\ref{hef})} for any fixed\/
$\,\th,\ha\in\bbR\,$ with\/ $\,\th\ne1\ne\ha\ne\th$. If the restrictions of\/
$\,E\,$ and\/ $\,\hh\,$ to the interval\/ $\,(-\infty,1)\,$ of the variable\/
$\,\ps\,$ are treated as continuous functions of the new variable\/
$\,F\in\bbR\hs$, differentiable 
on $\,\bbR\smallsetminus\{0\}$, cf.\ {\rm\S\ref{mono}}, then
\begin{equation}
\mathrm{i)}\quad\hh\,=\,\ax\,+\,\bx E\,+\,\cx F\,,\qquad\quad
\mathrm{ii)}\quad d\hh/dF\,=\,\cx\,+\,\bx\,dE/dF\label{hab}
\end{equation}
wherever\/ $\,F\ne0$, with\/ $\,\ax,\bx,\cx\,$ defined by {\rm(\ref{abc})}. In
addition, $\,\dsq E/dF^2<0\,$ wherever $\,F\ne0$, which, if\/ $\,\bx<0$,
amounts to $\,\dsq\hh/dF^2>0$.
\end{lem}
In fact, (\ref{hef}.i) and (\ref{abc}) give (\ref{hab}.i), and hence (\ref{hab}.ii),
while $\,\dsq E/dF^2<0\,$ by (\ref{dtf}) and (\ref{cnv}), so that (\ref{hab}.ii)
with $\,\bx<0\,$ yields $\,\dsq\hh/dF^2>0$.

\begin{prop}\label{dtefu}Given an integer\/ $\,m\ge2\,$ and a nontrivial
closed interval\/ $\,I\subset\bbR\,$ with the endpoints\/ $\,\th,\ha$, such
that\/ $\,0\,$ is an interior point of\/ $\,I$ and\/ $\,1\notin I$, let\/
$\,F,E,\jv\,$ be as in {\rm(\ref{fet})}, (\ref{spv}), and let\/
$\,Q\in\jv\smallsetminus\{0\}\,$ satisfy on $\,I\hs$ condition {\rm(\ref{qdq}.b)}.

With\/ $\,\ax,\bx,\cx\,$ depending on $\,\th,\ha\,$ as in {\rm(\ref{abc})} and
with\/ $\,\dot F=dF/d\ps$, $\,\dot E=dE/d\ps$, we then have\/
$\,\bx\dot F(\th)\dot F(\ha)\ne0$. In addition, $\,Q\,$ and\/ $\,I\hs$
satisfy {\rm(\ref{qdq}.a)}, (\ref{qdq}.c) and {\rm(\ref{qdq}.d)} if and only if one of the
following two conditions holds\/{\rm:}
\begin{enumerate}
\item $\dot E(\th)/\dot F(\th)\,<\,-\hs\cx/\bx\,<\,\dot E(\ha)/\dot F(\ha)$,
with the endpoints $\,\th,\ha$ of\/ $\,I\hs$ switched, if necessary, to ensure
that\/ $\,(-1)^m(u-v)>0$, or
\item $\ax/\bx\,>\,-\hs E(0)$.
\end{enumerate}
\end{prop}
\begin{proof}Our assumption on $\,I\hs$ states that $\,\th<0<\ha<1\,$ or
$\,\ha<0<\th<1$. As in (\ref{hef}.iii), let us set $\,F_0=F(\th)$,
$\,F_1=F(\ha)$. As $\,(-1)^m\dot F<0\,$ on $\,(-\infty,0)\cup(0,1)\,$ and
$\,F(0)=0$, cf.\ (\ref{mon}), the additional assumption $\,(-1)^m(u-v)>0\,$
made in (i) now gives $\,F_0<0<F_1$, while, by (\ref{abc}), $\,\bx<0\,$ and
$\,\bx\dot F(\th)\dot F(\ha)\ne0$.

Let $\,\hh\,$ be the function defined in (\ref{qth}) (that is, (\ref{hef}))
for our $\,\th,\ha$. Using the coordinate change $\,\ps\mapsto F=F(\ps)$, as
in Lemma~\ref{hdhdf}, we may treat $\,E\,$ and $\,\hh\,$ not as functions of
$\,\ps\in(-\infty,1)$, but rather as functions of the variable $\,F\in\bbR$,
continuous at $\,F=0\,$ and of class $\,C^\infty$ everywhere else. Then, as
$\,\dot E/\dot F=\hs dE/dF$,
\begin{equation}
\mathrm{Condition\ (i)\ means\ that}\hskip4ptd\hh/dF\hskip4pt
\mathrm{is\ positive\ at}\hskip4ptF=F_0\hskip2.6pt\mathrm{\ and\ negative\ at}
\hskip4ptF=F_1\hs,\label{cim}
\end{equation}
by (\ref{hab}.ii) with $\,\bx<0$, while, by (\ref{hab}.i) with $\,F(0)=0\,$ (cf.\ (\ref{fet})) and $\,\bx<0$,
\begin{equation}
\mathrm{Inequality\ (ii)\ states\ that}\hskip6pt\hh<0\hskip5pt
\mathrm{at}\hskip5ptF=0\,.\label{inq}
\end{equation}
Also, since $\,\hh\,$ and $\,Q\,$ are, up to a nonzero factor, related by (\ref{qth}),
\begin{equation}
\begin{array}{l}
\mathrm{Conditions\ (\ref{qdq}.c),\ (\ref{qdq}.d)\ for\ our\ pair}\hskip5ptQ,I
\hskip5pt\mathrm{mean\ that}\hskip6pt\hh\ne0
\hskip5pt\mathrm{at\ all}\hskip5ptF\\
\mathrm{with}\hskip9ptF_0<F<F_1\,,
\hskip6pt\mathrm{while}\hskip9ptd\hh/dF\ne0\hskip8pt\mathrm{at}\hskip7ptF
=F_0\hskip7pt\mathrm{and}\hskip6ptF=F_1\,.\end{array}
\label{ccd}
\end{equation}
In fact, to see if $\,d\hh/dF=(d\hh/d\ps)/(dF/d\ps)\,$ is zero or not, we only
need to apply $\,d/d\ps\,$ to the $\,Q(\ps)\,$ factor in (\ref{qth}), as
$\,Q(\th)=Q(\ha)=0$, while, by (\ref{dtf}), $\,dF/d\ps\ne0\,$ when
$\,0\ne\ps\ne1$.

Let $\,Q\,$ and $\,I\hs$ now satisfy (\ref{qdq}.a), (\ref{qdq}.c) and (\ref{qdq}.d). By (\ref{ccd}), $\,\hh\,$ must be nonzero throughout the
whole interval $\,(F_0,F_1)\,$ of the variable $\,F$, and $\,d\hh/dF\ne0\,$ at
the endpoints $\,F_0,\,F_1$. Thus, $\,\hh\,$ is positive (or, negative) on
$\,(F_0,F_1)$, which yields the clause about $\,\hh\,$ in (\ref{cim}) (or, (\ref{inq}), and hence (i) or, respectively, (ii).

Conversely, let us assume (i) or (ii). In case (i), using (\ref{cim}) and the
inequality $\,\dsq\hh/dF^2>0\,$ whenever $\,F\ne0\,$ (Lemma~\ref{hdhdf}), we
see that $\,d\hh/dF\,$ is positive for all $\,F\,$ with $\,F_0<F<0$, and
negative if $\,0<F<F_1$. As $\,\hh=0\,$ at both $\,F=F_0\,$ and $\,F=F_1$,
this in turn gives $\,\hh>0\,$ at every $\,F\,$ with $\,F_0<F<F_1$. Since, by (\ref{cim}), $\,d\hh/dF\ne0\,$ at the endpoints $\,F_0$, $F_1$, conditions (\ref{qdq}.c) and (\ref{qdq}.d) for $\,Q\,$ and $\,I\hs$ follow in view of
(\ref{ccd}) (while (\ref{qdq}.a) is obvious as $\,1\notin I$).

Finally, let us consider the remaining case (ii). By (\ref{inq}), we
then have $\,\hh<0\,$ at
$\,\,F=0$. It now follows that $\,d\hh/dF<0\,$ at $\,F=F_0$ and $\,\hh\ne0\,$
(so that $\,\hh<0$) everywhere in the interval $\,(F_0,0)$. In fact, if either
of these claims failed, we could find $\,F_2$, $\,F_3$ with
$\,F_0\le F_2\le F_3<0$, $\,d\hh/dF\ge0\,$ at $\,F=F_2$, and $\,\hh=0\,$ at
$\,F=F_3$. (Specifically, we set $\,F_2=F_3=F_0$ if $\,d\hh/dF\ge0\,$ at
$\,F=F_0\hs\mathrm{;}$ while, if $\,\hh=0\,$ at some $\,F_3$ in $\,(F_0,0)$,
we can use Rolle's theorem to select $\,F_2$.) Since $\,d\hh/dF\,$ is a
strictly increasing function of $\,F\in(F_0,0)\,$ (as $\,\dsq\hh/dF^2>0\,$ by
Lemma~\ref{hdhdf}), we have $\,d\hh/dF>0\,$ on $\,(F_3,0)$, which is not
possible as $\,\hh=0\,$ at $\,F=F_3$ and $\,\hh<0\,$ at $\,\,F=0$. A
completely analogous argument shows that, if (ii) holds, we must have
$\,d\hh/dF>0\,$ at $\,F=F_1$ and $\,\hh\ne0\,$ (that is, $\,\hh<0$) everywhere
in the interval $\,(0,F_1)$. In other words, by (\ref{ccd}), assuming (ii) we
obtain (\ref{qdq}.a), (\ref{qdq}.c) and (\ref{qdq}.d) as well. This completes the
proof.
\end{proof}

\section{Positivity}\label{posi}
\setcounter{equation}{0}
Let $\,\jv\,$ be the space (\ref{spv}) for a given integer $\,m\ge2$. According
to Lemma~\ref{onedi}, for every nontrivial closed interval $\,I$, a function
$\,Q\in\jv\smallsetminus\{0\}\,$ satisfying on $\,I\hs$ condition (\ref{qdq}.b)
exists and is unique up to a constant factor. We will say that $\,I\hs$
satisfies the {\it positivity condition\/} if, in addition to (\ref{qdq}.b), we
also have (\ref{qdq}.a), (\ref{qdq}.c) and (\ref{qdq}.d) for some, or any, such $\,Q$.

Now let $\,m\ge2\,$ be fixed. The discussion in the preceding sections has
determined that for any given nontrivial closed interval $\,I\hs$ the
positivity condition holds
\begin{enumerate}
\item Always, if $\,I\hs$ contains neither $\,0\,$ nor $\,1\,$
(Proposition~\ref{btacd}).
\item Never, if $\,I\hs$ contains $\,1\,$ as an interior point or $\,0\,$
as an endpoint. This is clear from Remark~\ref{endpo} or, respectively, the fact
that, by (\ref{dqo}.a), condition (\ref{qdq}.b) for the endpoint $\,0\,$ contradicts
condition (\ref{qdq}.d).
\item When $\,I\hs$ contains $\,1\,$ as an endpoint: if and only if the
other endpoint lies in $\,(-\infty,\zu)\cup(0,\infty)\,$ for the number
$\,\zu<0\,$ defined in \S\ref{mopr}. (See below.)
\item When $\,I\hs$ contains $\,0\,$ as an interior point and does not
contain $\,1\,$ at all: if and only if the endpoints $\,\th,\ha\,$ of $\,I\hs$
satisfy (i) or (ii) in Proposition~\ref{dtefu}.
\end{enumerate}
Only (iii) still requires an explanation. Namely, as $\,1\in I$, condition (\ref{qdq}.a) gives $\,\cx=0$, so that (\ref{qef}) becomes (\ref{qta}), and the remaining endpoint $\,\th\ne1\,$ determines $\,Q\,$ up to a
factor via Lemma~\ref{czero}(i). Hence, by Lemma~\ref{czero}(ii), (iii)a),b),
conditions (\ref{qdq}.d) and (\ref{qdq}.c), in addition to (\ref{qdq}.b) and (\ref{qdq}.a), amount
to requiring that $\,E(\th)\hs\dot E(\th)\ne0\,$ and $\,E(\ps)\ne E(\th)\,$
for all $\,\ps\,$ in the open interval connecting $\,1\,$ and $\,\th$. The
claims made in (iii) now are trivial consequences of the descriptions in
\S\ref{mopr} of the monotonicity intervals for $\,E\,$ and the roots of $\,E\,$ and
$\,\dot E$.

\section{Factorization of \hs{\rm(\ref{qdq}.e)}}\label{fact}
\setcounter{equation}{0}
\begin{lem}\label{dvsbl}Let\/ $\,\varPhi,\vh\,$ be polynomials in two
variables\/ $\,\th,\ha$. Then\/ $\,\varPhi\,$ is divisible by\/ $\,\vh\,$ if
\begin{enumerate}
\item $\mathrm{deg}\hskip2pt\vh=1$, while\/ $\,\varPhi=0\,$ wherever\/
$\,\vh=0$, or
\item $\vh=(\th-\ha)^3$ and\/ $\,\varPhi\hs$ is antisymmetric, while\/
$\,\partial\varPhi/\partial\ha=0\,$ wherever\/ $\,\th=\ha$.
\end{enumerate}
\end{lem}
In fact, (i) is clear if we use new affine coordinates $\,\xh,\yh\,$ with
$\,\vh=\xh$. Next, for $\,\varPhi\,$ as in (ii),
$\,\tilde\varPhi=\varPhi/(\ha-\th)\,$ is, by (i), a symmetric polynomial with
$\,\tilde\varPhi=\,\partial\varPhi/\partial\ha\,$ wherever $\,\th=\ha$. Thus,
for the new coordinates $\,\xh,\yh\,$ given by
\begin{equation}
\th\,=\,\xh+\yh\,,\quad\ha\,=\,\xh-\yh\,,\hskip20pt\mathrm{that\ is,}\hskip15pt
\xh\,=\,(\th+\ha)/2\,,\quad\yh\,=\,(\th-\ha)/2\,,\label{uxe}
\end{equation}
$\tilde\varPhi\,$ is even in $\,\yh\,$ (and so it is a polynomial in $\,\xh\,$
and $\,\yh^2$) and vanishes wherever $\,\yh=0\,$ (due to the assumption on
$\,\partial\varPhi/\partial\ha$). Hence, by (i), $\,\tilde\varPhi\,$ is
divisible by $\,\yh^2\nh$.

\begin{lem}\label{cpiuv}With\/ $\,F,\,E\,$ as in {\rm(\ref{fet})} for an integer\/
$\,m\ge2$, and\/ $\,(\,\,)\dot{\,}=\,d/d\ps$, let
\begin{equation}
\pv(\th,\ha)\,=\,(\ha-1)\left\{\lj F(\ha)-F(\th)\rj\hs\dot
E(\ha)\,-\,\lj E(\ha)-E(\th)\rj\hs\dot F(\ha)\right\}\nh.\label{puv}
\end{equation}
Then, for\/ $\,\ax,\bx\,$ depending on $\,\th,\ha\,$ as in {\rm(\ref{abc})},
$\,\th(\th-2)\ha(\ha-2)\hs\lj\pv(\th,\ha)-\pv(\ha,\th)\rj\,$ equals
\begin{equation}
2\hs(\th\ha-\th-\ha)\hskip7pt\mathrm{times}\hskip7pt
\lj\hs m\hs\th\ha-(2m-1)(\th+\ha-2)\hs\rj\ax\,
+\,(2m-1)\hs(\th+\ha-2)\hs\bx\hs\vs(0)\,.\label{muv}
\end{equation}
\end{lem}
In fact, by (\ref{puv}) and (\ref{dte}),
$\hs\th(\th-2)\ha(\ha-2)\hs\lj\pv(\th,\ha)-\pv(\ha,\th)\rj
=\lj\ha(\ha-2)\varLambda(\th)+\th(\th-2)\varLambda(\ha)\rj\ax
+2(2m-1)\hs\lj\th(\th-2)(\ha-1)+\ha(\ha-2)(\th-1)\rj\bx\vs(0)$, with
$\,\varLambda\,$ as in (\ref{dtf}). This proves our claim, since
$\,\ha(\ha-2)\varLambda(\th)+\th(\th-2)\varLambda(\ha)
=2m\hs(\th\ha-\th-\ha)\hs\lj\th\ha-(2-1/m)(\th+\ha-2)\rj\,$ and
\begin{equation}
\th(\th-2)(\ha-1)\,+\,\ha(\ha-2)(\th-1)\,=\,(\th\ha-\th-\ha)(\th+\ha-2)\,.
\label{uuv}
\end{equation}
\begin{lem}\label{polyt}For any integer\/ $\,m\ge2$, there exists a unique
symmetric polynomial\/ $\,T\,$ in the variables\/ $\,\th,\ha\,$ such that,
for\/ $\,\pv\,$ as in {\rm(\ref{puv})},
\begin{equation}
(\th-1)^m(\ha-1)^m\,\lj\pv(\th,\ha)-\pv(\ha,\th)\rj\,
=\,2(\ha-\th)^3(\th\ha-\th-\ha)\,T(\th,\ha)\,.\label{umv}
\end{equation}
\end{lem}
\begin{proof}Multiplication by $\,(\th-1)^m(\ha-1)^m$ turns (\ref{muv}), as well
as $\,\pv(\th,\ha)\,$ and $\,\pv(\ha,\th)$, into polynomials in $\,\th,\ha\,$
(cf.\ (\ref{puv}) and (\ref{fet}), (\ref{dtf})). Lemmas~\ref{dvsbl} --~\ref{cpiuv} now show
that the left-hand side of (\ref{umv}) is a polynomial divisible by
$\,\th\ha-\th-\ha$.

On the other hand, $\,\pv\,$ in (\ref{puv}) clearly vanishes whenever
$\,\th=\ha$. Also, the $\,\,\partial\varPhi/\partial\ha\,$ clause in Lemma~\ref{dvsbl}(ii) is satisfied both by $\,\varPhi=\pv\,$ and
$\,\varPhi(\th,\ha)=\pv(\ha,\th)$. (One verifies this without evaluating
$\,\ddot E,\ddot F$, since, in the Leibniz-rule expression for the partial
derivative, only the factors $\,F(\ha)-F(\th)\,$ and $\,E(\ha)-E(\th)\,$ need
to be differentiated, as they vanish when $\,\th=\ha$.) In view of Lemma~\ref{dvsbl}(ii) and the last paragraph, the polynomial on the left-hand side of
(\ref{umv}) divided by $\,\th\ha-\th-\ha\,$ is still divisible by
$\,(\ha-\th)^3$. This completes the proof.
\end{proof}
\begin{lem}\label{dotso}Given an integer\/ $\,m\ge2\,$ and\/ a
nontrivial closed interval\/ $\,I\subset\bbR\,$ with the endpoints\/
$\,\ha=1\,$ and\/ $\,\th\ne1$, let\/ $\,\dot\vs=\,d\vs/d\ps\,$ for\/
$\,\vs\,$ defined by {\rm(\ref{esg})}. Also, let\/ $\,\jv\,$ be the space\/
{\rm(\ref{spv})} and let\/ $\,Q\in\jv\smallsetminus\{0\}\,$ be a function,
unique up to a factor, which satisfies {\rm(\ref{qdq}.b)} on $\,I$, cf.\
{\rm Lemma~\ref{onedi}}. Condition\/ $\,\dot\vs(\th)=0\,$ then is necessary
and sufficient for $\,Q\,$ and\/ $\,I\hs$ to satisfy {\rm(\ref{qdq}.e)}.
\end{lem}
In fact, as $\,1\,$ is an endpoint of $\,I$, (\ref{qdq}.b) gives (\ref{qdq}.a) and we
have (\ref{qef}) with $\,\cx=0$, that is, (\ref{qta}) for some $\,\bx\ne0\,$
and
$\,\ax\,$ with $\,E(\th)=-\hs\ax/\bx\,$ (see Lemma~\ref{czero}(i)). Conditions
$\,\dot\vs(\th)=0\,$ and (\ref{qdq}.e) now are equivalent by Lemma~\ref{czero}(iii)c).
\begin{thm}\label{facto}Given an integer\/ $\,m\ge2\,$ and a nontrivial
closed interval\/ $\,I$ with endpoints\/ $\,\th,\ha$, let a function
$\,Q\in\jv\smallsetminus\{0\}$, with\/ $\,\jv\,$ as in {\rm(\ref{spv})}, satisfy
on $\,I$ condition {\rm(\ref{qdq}.b)}, that is, $\,Q(\th)=Q(\ha)=0$. By\/ {\rm Lemma~\ref{onedi}}, such\/ $\,Q\,$ exists and is unique up to a constant factor. Also,
let\/ $\,(\,\,)\dot{\,}=\,d/d\ps$.

Then $\,Q\,$ satisfies {\rm(\ref{qdq}.e)}, that is,
$\,\dot Q(\th)+\dot Q(\ha)=0$, if
and only if, for the polynomials\/ $\,\vs,T\,$ defined in {\rm(\ref{esg})} and\/
{\rm Lemma~\ref{polyt}}, one of the following three cases occurs\/{\rm:}
\begin{equation}
\begin{array}{rl}\arraycolsep9pt
\mathrm{a)}&
\th=1\hskip8pt\mathrm{and}\hskip7pt\dot\vs(\ha)=0\,,\qquad\mathrm{or}
\qquad\ha=1\hskip8pt\mathrm{and}\hskip7pt\dot\vs(\th)=0\,.\hskip43pt\\
\mathrm{b)}&
\th\ha\,=\,\th\,+\,\ha\,.\\
\mathrm{c)}&
T(\th,\ha)\,=\,0\qquad\quad\mathrm{and}\qquad\th\ne1\ne\ha\hs.\end{array}
\label{trc}
\end{equation}
\end{thm}
\begin{proof}If one of $\,\th,\ha\,$ equals $\,1$, our assertion, stating
in this case that condition (\ref{qdq}.e) is equivalent to (\ref{trc}.a), is
nothing
else than Lemma~\ref{dotso}. Let us therefore assume that $\,\th\ne1\ne\ha$,
and consider a function $\,Q\in\jv\smallsetminus\{0\}\,$ satisfying on $\,I\hs$
condition (\ref{qdq}.b). Using the subscript convention (\ref{hef}.iii) (also for
functions of $\,\ps\,$ other than $\,E$, $\,F$), we can rewrite (\ref{qdq}.e) as
$\,0=\dot Q_0+\dot Q_1=(\th-1)\dot\hh_0+(\ha-1)\dot\hh_1\,$ (see (\ref{qth}); by (\ref{qdq}.b), $\,\hh_0=\hh_1=0$). Thus, (\ref{qdq}.e) states that
\begin{equation}
(\th-1)\lj(F_1-F_0)\dot E_0-(E_1-E_0)\dot F_0\rj
+(\ha-1)\lj(F_1-F_0)\dot E_1-(E_1-E_0)\dot F_1\rj=0\hs,
\end{equation}
cf.\ (\ref{hef}.ii), and so $\,\pv(\ha,\th)=\pv(\th,\ha)$, for $\,\pv\,$ given
by (\ref{puv}). As $\,\th\ne1\ne\ha\ne\th$, (\ref{umv}) shows that, when $\,\th\ne1\ne\ha$, (\ref{qdq}.e) holds if and only if we have (\ref{trc}.b) or
(\ref{trc}.c). Since, conversely, (\ref{trc}.b) implies that
$\,\th\ne1\ne\ha$, this completes the proof.
\end{proof}
Conditions (\ref{trc}.b), (\ref{trc}.c) are not mutually exclusive: when
$\,m\,$ is odd, they may occur simultaneously, with $\,\th\ne\ha$. See
Lemma~\ref{xxast}(b) in \S\ref{sgnt}.

\section{First subcase of \hs{\rm(\ref{qdq}.e)}: condition
\hs{\rm(\ref{trc}.a)}}\label{sbon}
\setcounter{equation}{0}
Let $\,m\ge2\,$ be a fixed integer. We will now find all pairs $\,Q,I\hs$
formed by a function $\,Q\in\jv\smallsetminus\{0\}$, with $\,\jv\,$ as in
(\ref{spv}), and a nontrivial closed interval $\,I\hs$ such that $\,1\in I\hs$
and $\,Q\,$ satisfies on $\,I\hs$ the boundary conditions (\ref{qdq}.b), (\ref{qdq}.e).
As we show in Proposition~\ref{onein} below, such $\,Q,I\hs$ exist if and only if
$\,m\,$ is odd, and then they are unique up to multiplication of $\,Q\,$ by a
nonzero constant.

With $\,(\,\,)\dot{\,}=\,d/d\ps$, (\ref{esg}) gives
$\,(\ps-1)\dot\vs(\ps)=\dot E(\ps)-\vs(\ps)\,$ for all real $\,\ps$.
Multiplying this by $\,\ps(\ps-1)(\ps-2)\,$ and using (\ref{dte}.ii), (\ref{esg}), (\ref{dtf}.b), we easily obtain
\begin{equation}
\ps(\ps-1)(\ps-2)\dot\vs(\ps)\,
=\,[(m-1)(\ps-2)^2+2]\hs\vs(\ps)\,-\,2(2m-1)\hs\vs(0)\hs.\label{tds}
\end{equation}
\begin{rem}\label{zerpo}We will repeatedly use the following obvious fact.
If the derivative of a $\,C^1$ function $\,\varPhi\,$ on an interval
$\,\mathcal{I}\,$ is positive wherever $\,\varPhi=0\,$ in $\,\mathcal{I}$,
then there exists at most one $\,t\in\mathcal{I}\,$ with $\,\varPhi(t)=0$. If
such $\,t\,$ exists, then $\,\varPhi<0\,$ on $\,(-\infty,t)\cap\mathcal{I}\,$
and $\,\varPhi>0\,$ on $\,(t,\infty)\cap\mathcal{I}$.
\end{rem}
\begin{prop}\label{sgrts}Given an integer\/ $\,m\ge2$, let\/ $\,\vs\,$ be
the polynomial in {\rm(\ref{esg})}, and let\/
$\,X(\ps)=2(2m-1)\hs\vs(0)\hs[(m-1)(\ps-2)^2+2]^{-1}\hskip-1.3pt$.
\begin{enumerate}
\item[a)] If\/ $\,m\,$ is even, then $\,\dot\vs=\hs d\vs/d\ps\,$ has no real
roots.
\item[b)] If\/ $\,m\,$ is odd, $\,\dot\vs\,$ has exactly one real root\/
$\,\sx$, and that root is negative.
\item[c)] $X(\ps)>0\,$ for all\/ $\,\ps\in\bbR\hs$. For any $\,\ps<0$, we have
$\,\dot\vs(\ps)=0\,$ if and only if\/ $\,\vs(\ps)=X(\ps)$.
\end{enumerate}
\end{prop}
\begin{proof}The function $\,X\,$ is positive since $\,\vs(0)>0\,$ (cf.\
(\ref{sto}.i)), and so (\ref{tds}) yields (c). Next, one and only one of the
following two conditions must hold:
\begin{enumerate}
\item $\dot\vs>0\,$ and $\,\vs<X\,$ everywhere in $\,(-\infty,0)$.
\item There exists a unique $\,\sx\in(-\infty,0)\,$ with
$\,\dot\vs(\sx)=0$, and then $\,\vs>X\,$ on $\,(-\infty,\sx)$, while
$\,\vs(\sx)=X(\sx)\,$ and $\,\vs<X\,$ on $\,(\sx,0)$.
\end{enumerate}
In fact, the definition of $\,X\,$ easily gives
\begin{enumerate}
\item[a)] $\dot X=\hs dX/d\ps\,$ is positive everywhere in $\,(-\infty,0\hs]$.
\hskip15ptb)\hskip9pt$X(0)=\vs(0)\,$ and $\,\dot X(0)<\dot\vs(0)$,
\end{enumerate}
the last inequality being clear as $\,\dot X(0)=(2m-2)\vs(0)/(2m-1)$, while
$\,\dot\vs(0)=\vs(0)>0\,$ (see (\ref{sto}.i)). Now let $\,\varPhi=X-\vs$. By (d)
and (c), $\,\dot\varPhi>0\,$ wherever $\,\varPhi=0\,$ in $\,(-\infty,0)$.
Thus, according to Remark~\ref{zerpo}, if some $\,\sx\in(-\infty,0)\,$ has
$\,\dot\vs(\sx)=0\,$ (that is, $\,\varPhi(\sx)=0$, cf.\ (c)), then such $\,\sx\,$
is unique, while $\,\varPhi<0\,$ on $\,(-\infty,\sx)\,$ and $\,\varPhi>0\,$ on
$\,(\sx,0)$. This yields (ii). However, if no such $\,\sx\,$ exists,
we have $\,\varPhi\dot\vs\ne0\,$ everywhere in $\,(-\infty,0)$. Hence
$\,\dot\vs>0\,$ on $\,(-\infty,0\hs]$, as $\,\dot\vs(0)>\dot X(0)>0\,$ by (d)
and (e). Also, by (e), $\,\dot\varPhi(0)<0=\varPhi(0)$, so that
$\,\varPhi(\ps)>0\,$ for all $\,\ps<0\,$ close to $\,0$. Thus, $\,\varPhi>0\,$
on $\,(-\infty,0)$, and so $\,\vs<X\,$ on $\,(-\infty,0)$, which gives (i).

In case (i), Remark~\ref{roots} shows that $\,\dot\vs\,$ has no real roots, and so
$\,m\,$ is even (as $\,\,\deg\,\dot\vs=m-2$). In the remaining case (ii),
$\,\dot\vs\,$ has exactly one real root $\,\sx\,$ (cf.\ Remark~\ref{roots}), and so
$\,m\,$ must be odd, which is clear since $\,\vs\,$ is a polynomial of degree
$\,m-1\,$ with a positive leading coefficient, cf.\ (\ref{esg}), while
$\,\vs>X>0\,$ on $\,(-\infty,\sx)\,$ by (ii) and (c). This yields (a) and (b),
completing the proof.
\end{proof}
\begin{prop}\label{onein}For a fixed integer\/ $\,m\ge2$, let\/
$\,F,E,\vs,\jv\,$ be given by {\rm(\ref{fet})}, {\rm(\ref{esg})},
{\rm(\ref{spv})}. The following two conditions are equivalent\/{\rm:}
\begin{enumerate}
\item There exist a nontrivial closed interval\/ $\,I$ containing\/
$\,1\hs$ and a function $\,Q\in\jv\smallsetminus\{0\}\,$ satisfying on
$\,I$ conditions {\rm(\ref{qdq}.b)}, {\rm(\ref{qdq}.e)}.
\item $m\,$ is odd.
\end{enumerate}
If\/ $\,m\,$ is odd, $\,Q,I$ in {\rm(i)} are unique up to multiplication of\/
$\,Q\,$ by a nonzero constant. Specifically, the endpoints of\/ $\,I\hs$ then
are\/ $\,1\,$ and the unique\/ $\,\sx\in(-\infty,0)\,$ with\/
$\,\dot\vs(\sx)=0$, cf.\ {\rm Proposition~\ref{sgrts}(b)}, while\/ $\,Q\,$ is
given by {\rm(\ref{qta})} for any constants $\,\bx\ne0\,$ and\/ $\,\ax\,$ with
$\,E(\sx)=-\hs\ax/\bx$.
\end{prop}
In fact, if $\,m\,$ is even, $\,\dot\vs(\th)\ne0\,$ for all $\,\th\in\bbR\,$
(Proposition~\ref{sgrts}(a)); thus, by Remark~\ref{endpo} and Lemma~\ref{dotso}, no pair
$\,Q,I\hs$ has the properties listed in (i).

Conversely, if $\,m\,$ is odd, Lemma~\ref{dotso}, Proposition~\ref{sgrts}(b) and
Remark~\ref{endpo} imply both the existence of a pair $\,Q,I\hs$ with (i), and
uniqueness of the interval $\,I$, the endpoints of which must be those
described in our assertion. Finally, uniqueness of $\,Q\,$ up to a factor and
its required form are immediate from Remark~\ref{ptone}.

\section{Failure of positivity in case \hs{\rm(\ref{trc}.a)}}\label{fail}
\setcounter{equation}{0}
In this section we show (Proposition~\ref{noone}) that case (\ref{trc}.a) of
condition (\ref{qdq}.e) cannot occur simultaneously with the positivity
condition introduced in \S\ref{posi}.
\begin{lem}\label{sleta}Given an odd integer\/ $\,m\ge3$, let\/ $\,\sx\,$ be
the unique real number with\/ $\,\dot\vs(\sx)=0\,$ for\/ $\,\vs\,$ as in
{\rm(\ref{esg})}, cf.\ {\rm Proposition~\ref{sgrts}(b)}. Then, with\/
$\,\zu,\wu\,$ as in {\rm\S\ref{mopr}},
\begin{equation}
\begin{array}{rll}\arraycolsep9pt
\mathrm{i)}&
\sx\,\le\,-\hs\eta\,,&\mathrm{where\ }\,\,\eta\,=\,(2m-3)/[3(m-2)]\,,\\
\mathrm{ii)}&
\vs(\ps)\,\ge\,\vs(\sx)\,>\,0\quad&\mathrm{for\ all\ }\,\,\ps\in\bbR\,,\\
\mathrm{iii)}&
\zu\,\le\,\sx\,<\,\wu\,<\,0\,,&\mathrm{and\ }\,\,
\zu\,<\,\sx\hskip6pt\mathrm{if}\hskip6ptm>3\hs.\end{array}
\label{zsw}
\end{equation}
\end{lem}
\begin{proof}Differentiating the second equality in (\ref{esg}) we obtain
$\,\dot\vs(\ps)=\sum_{k=2}^ma_k\ps^{k-2}$ with coefficients $\,a_k>0\,$ such
that $\,a_k/a_{k+1}=[1-2/(k+1)][1+(m-1)/(m-k)]$, where both factors in square
brackets clearly are increasing positive functions of $\,k=2,\dots,m-1$. Thus,
$\,a_k/a_{k+1}\ge a_2/a_3=\eta\,$ (cf.\ (\ref{zsw}.i)), that is,
$\,a_k\ge a_{k+1}\eta$, and so
$\,\dot\vs(-\hskip.5pt\eta)=(a_2-a_3\eta)+(a_4-a_5\eta)\eta^2
+\ldots+(a_{m-1}-a_m\eta)\eta^{m-3}\ge0$. As $\,\dot\vs\,$ is an odd-degree
polynomial with a positive leading coefficient and with a unique root at
$\,\sx$, it must be negative on $\,(-\infty,\sx)\,$ and positive on
$\,(\sx,\infty)$, so that the last inequality proves (\ref{zsw}.i), and
$\,\vs(\ps)\ge\vs(\sx)\,$ for all $\,\ps\in\bbR\hs$. Furthermore,
Proposition~\ref{sgrts}(c), with $\,\ps=\sx$, gives $\,\vs(\sx)=X(\sx)>0\,$ for $\,X(\ps)\,$
as in Proposition~\ref{sgrts}, which yields the remaining inequality in (\ref{zsw}.ii).

In the remainder of this proof, all inequalities are strict if $\,m>3$. First,
$\,0\le(2m-5)(m-3)=(2m-3)(m-1)-6(m-2)\,$ (since $\,m\ge3$), so
that $\,\eta\ge2/(m-1)\,$ for $\,\eta\,$ as in (\ref{zsw}.i). Now (\ref{zsw}.i)
yields $\,\sx\le-\,2/(m-1)$. This gives $\,[(m-1)\sx+2]\sx\ge0\,$ and
$\,[(m-1)\sx+2]\sx\hs\vs(0)\,\ge\,0\,$ (since $\,\vs(0)>0\,$ by (\ref{sto}.i)).
As $\,\dot\vs(\sx)=0$, the right-hand side of (\ref{tds}) with $\,\ps=\sx\,$
vanishes, which, multiplied by $\,\sx-1$, reads
$\,[(m-1)(\sx-2)^2+2]\hs E(\sx)-2(2m-1)(\sx-1)\hs\vs(0)=0\,$ (see (\ref{esg})).
Adding this side-by-side to the last inequality we get
$\,[(m-1)(\sx-2)^2+2]\hs[E(\sx)+\vs(0)]\,\ge\,0$. Hence
$\,E(\sx)\ge-\,\vs(0)=E(0)\,$ (cf.\ (\ref{esg})). Also, as
$\,\dot E(\ps)=\vs(\ps)+(\ps-1)\dot\vs(\ps)\,$ (by (\ref{esg})) and
$\,\dot\vs(\sx)=0$, (\ref{zsw}.ii) gives $\,\dot E(\sx)=\vs(\sx)>0$. Since
$\,\wu<0\,$ (see \S\ref{mopr}) and $\,\sx<0$, while $\,E(\sx)\ge E(0)\,$ and
$\,\dot E(\sx)>0\,$ (see above), (\ref{emo}) yields (\ref{zsw}.iii). This
completes the proof.
\end{proof}
\begin{prop}\label{noone}Given an integer\/ $\,m\ge2$, there exists no
function $\,Q\in\jv$, for\/ $\,\jv\,$ as in {\rm(\ref{spv})}, satisfying\/
{\rm(\ref{qdq})} on any nontrivial closed interval\/ $\,I$ that contains
$\,1\nh$.
\end{prop}
In fact, let (\ref{qdq}.a,\hs b,\hs e) hold for $\,Q\,$ and $\,I$. By
Proposition~\ref{onein}, $\,m\,$ is odd, and the endpoint of $\,I$ other than
$\,1\,$ is the unique $\,\sx<0\,$ with $\,\dot\vs(\sx)=0$. However, by
(\ref{zsw}.iii), $\,\sx\,$ fails the ``positivity test'' in (iii) of
\S\ref{posi}, that is, $\,Q\,$ and $\,I\hs$ cannot satisfy all four conditions
(\ref{qdq}.a) -- (\ref{qdq}.d).
\begin{rem}\label{uvone}
For the moduli curve $\,\mathcal{C}\hs$ defined at the end
of \S\ref{modu} and any point $\,(\th,\ha)\in\mathcal{C}\smallsetminus\{(0\hs,\nh0)\}$,
Proposition~\ref{noone} clearly implies that $\,1<\th<\ha\,$ or $\,\th<\ha<1$.
\end{rem}

\begin{example}\label{zswtr}For $\,m=3$, we have $\,\zu=\sx=-\hs1\,$ and
$\,\wu=-\hs2/3\,$ in (\ref{zsw}.iii). In fact, $\,-\hs1,-\hs1,-\hs2/3\,$ are the
unique negative roots of $\,E-E(0)$, $\,\dot\vs\,$ and $\,\dot E$, since, by (\ref{esg}), $\,\vs(\ps)=\ps^2+2\ps+2\,$ for $\,m=3$, and so
$\,E(\ps)=\ps^3+\ps^2-2$, $\,E(\ps)-E(0)=\ps^2(\ps+1)$,
$\,\dot\vs(\ps)=2(\ps+1)\,$ and $\,\dot E(\ps)=\ps(3\ps+2)$. Cf.\ Example~\ref{xymtr} and the end of \S\ref{anth}.
\end{example}

\section{Second subcase of \hs{\rm(\ref{qdq}.e)}: condition
\hs{\rm(\ref{trc}.b)}}\label{sbto}
\setcounter{equation}{0}
In this section we explicitly classify all pairs $\,Q,I\hs$ that satisfy the
boundary conditions (\ref{qdq}.b), (\ref{qdq}.e) and are of type (\ref{trc}.b).

Given an integer $\,m\ge2$, let $\,F,E\,$ be as in (\ref{fet}). Then, for all
$\,\,\ps\in\bbR\smallsetminus\{1\}$,
\begin{equation}
\mathrm{i)}\quad F(\hatt)\hs=\hs-\hs F(\ps)\hs,\qquad\quad\mathrm{ii)}\quad
E(\hatt)\hs=\hs E(\ps)\hs-\hs F(\ps),\label{fta}
\end{equation}
where $\,\hatt\,$ is, for any $\,\ps\ne1$, given by either of the two
equivalent relations
\begin{equation}
\mathrm{a)}\quad\hatt\,=\,\ps/(\ps-1)\,,\hskip30pt\mathrm{b)}\quad
\hatt-1\,=\,1/(\ps-1)\,.\label{ast}
\end{equation}
In fact, (\ref{fta}.i) is immediate from (\ref{fet}), and (\ref{fta}.ii) is easily
verified by induction on $\,m\hs$, using (\ref{fmt}). Since (\ref{ast}.a) is
equivalent to (\ref{ast}.b), we obtain
\begin{lem}\label{uuast}The assignment\/ $\,\th\mapsto\hu$ with\/
$\,\hu=\th/(\th-1)\,$ is an involution of\/ $\,\bbR\smallsetminus\{1\}$,
decreasing both on $\,(-\infty,1)\,$ and on $\,(1,\infty)$. It interchanges
the interval\/ $\,(1,2\hs]\,$ with $\,[\hs2,\infty)$, and\/
$\,(-\infty,0\hs]\,$ with $\,[\hs0,1)$, while its fixed points are $\,0\,$
and\/ $\,2$.
\end{lem}

\begin{lem}\label{uastq}Given an integer\/ $\,m\ge2\,$ and a real number\/
$\,\th\notin\{0,1,2\}$, let\/ $\,\ha=\hu$, where\/ $\,\hu=\th/(\th-1)$.
Formula {\rm(\ref{qth})} then defines a function $\,Q\in\jv\smallsetminus\{0\}$,
for $\,\jv\,$ as in {\rm(\ref{spv})}, for which {\rm(\ref{qdq}.a)},
{\rm(\ref{qdq}.b)} and\/ {\rm(\ref{qdq}.e)}
hold on the nontrivial closed interval\/ $\,I$ with the endpoints\/
$\,\th,\ha$. Moreover, $\,\th,\ha\,$ satisfy condition {\rm(\ref{trc}.b)}.
\end{lem}
In fact, $\,\th\ne\ha\,$ since $\,\th\notin\{0,1,2\}\,$ (cf.\ Lemma~\ref{uuast}).
Thus, in view of (\ref{qth}), $\,Q\in\jv\smallsetminus\{0\}\,$ and (\ref{qdq}.b)
holds. Also, by Lemma~\ref{uuast}, $\,1\notin I$, which implies (\ref{qdq}.a). Finally,
relation $\,\ha=\hu$ is nothing else than (\ref{trc}.b). Theorem~\ref{facto} now
gives (\ref{qdq}.e) for $\,Q\,$ and $\,I$.

\begin{thm}\label{qihyp}Let\/ $\,m\ge2\,$ be a fixed integer. Assigning to
each real number\/ $\,\th\,$ with\/ $\,\th<0\,$ or\/ $\,1<\th<2\,$ the
interval\/ $\,I=[\th,\ha]\,$ with\/ $\,\ha=\th/(\th-1)\,$ and the function
$\,Q\,$ defined up to a factor by {\rm(\ref{qth})}, with\/ $\,E,F\,$ as in
{\rm(\ref{fet})}, we obtain a bijective correspondence between
\begin{enumerate}
\item The subset\/ $\,(-\infty,0)\cup(1,2)\,$ of\/ $\,\bbR$, and
\item The set of all equivalence classes, modulo multiplication of\/
$\,Q\,$ by nonzero scalars, of pairs\/ $\,Q,I\hs$ formed by a nontrivial
closed interval\/ $\,I\subset\bbR\,$ and a function
$\,Q\in\jv\smallsetminus\{0\}$, with\/ $\,\jv\,$ as in {\rm(\ref{spv})}, such
that\/ $\,Q\,$ and\/ $\,I\hs$ satisfy conditions {\rm(\ref{qdq}.b)} and\/ {\rm(\ref{trc}.b)},
i.e., $\,Q(\th)=Q(\ha)=0\,$ and\/ $\,\th\ha=\th+\ha$, where\/ $\,\th,\ha\,$
are the endpoints of\/ $\,I$.
\end{enumerate}
\end{thm}
\begin{proof}That $\,I\hs$ and $\,Q\,$ described above have the properties
listed in (ii) is clear from Lemma~\ref{uastq}. Injectivity of our assignment is
obvious since $\,Q,I\hs$ obtained from the given $\,\th\,$ determine $\,\th\,$
uniquely (as the lower endpoint of $\,I$).

To show that the assignment is surjective, let us fix $\,Q,I\hs$ as in (ii),
and let $\,I=[\th,\ha]$. Since $\,\th,\ha\,$ satisfy (\ref{trc}.b), we have
$\,(\th-1)\ha=\th$, so that $\,\th\ne1\,$ and $\,\ha=\hu$ (notation of
(\ref{ast})). As $\,\th<\ha$, this gives $\,\th<0\,$ or $\,1<\th<2\,$ (cf.\
Lemma~\ref{uuast}). Lemma~\ref{abceq} now shows that $\,Q\,$ is, up to a factor, given by (\ref{qef}) with (\ref{abc}), that is, by (\ref{qth}). Thus, the equivalence class
of $\,Q,I\hs$ is the image of $\,\th\,$ under our assignment, which completes
the proof.
\end{proof}

\section{Another rational function}\label{anth}
\setcounter{equation}{0}
Given an integer $\,m\ge2$, we define a rational function $\,G\,$ of the
variable $\,\ps\,$ by
\begin{equation}
G\,=\,E\,-\,F/2\,,\qquad\mathrm{with}\quad E,F\quad\mathrm{as\ in\
(\ref{fet}).}\label{gef}
\end{equation}
For any $\,\ps\in\bbR\smallsetminus\{1\}\,$ we then have, with
$\,\hatt=\ps/(\ps-1)\,$ as in (\ref{ast}) and $\,(\,\,)\dot{\,}=\,d/d\ps$,
\begin{equation}
\begin{array}{rl}\arraycolsep9pt
\mathrm{a)}&
G(\hatt)\,=\,G(\ps)\,,\qquad\mathrm{b)}\quad
\dot G(\hatt)/\dot F(\hatt)\,=\,-\hs\dot G(\ps)/\dot F(\ps)\quad\mathrm{if}
\hskip5pt\ps\hs\ne\hs0\hs,\\
\mathrm{c)}&
\ps(\ps-1)(\ps-2)\dot G(\ps)\,=\,\varLambda(\ps)\hs G(\ps)\,
-\,2(2m-1)(\ps-1)\hs\vs(0)\hs,\\
\mathrm{d)}&
G(0)\,=\,-\hs\vs(0)\,<\,0\hs,\hskip8pt\dot G(0)\,=\,0\hs,\quad\mathrm{e)}
\hskip9pt(2m-3)\hs\ddot G(0)\,=\,m\hs\vs(0)\,>\,0\hs,\hskip4pt\\
\mathrm{f)}&
G(2)\,=\,(2m-1)\hs\vs(0)\hs,\qquad\dot G(2)\,=\,0\hs,\end{array}
\label{gta}
\end{equation}
where $\,\varLambda\,$ is given by (\ref{dtf}.b). In fact, (\ref{gta}.a),
which is immediate from (\ref{gef}) and (\ref{fta}), states that $\,G\,$
restricted to either of the intervals $\,(-\infty,1)$, $\,(1,\infty)\,$ is an
even function of the new variable $\,F\in(-\infty,\infty)$, cf.\ \S\ref{mono}
and (\ref{fta}.i). As the derivative $\,d\hskip.4ptG/dF=\dot G/\dot F\,$ of
the even function $\,G\,$ is odd, (\ref{gta}.b) follows. (Also,
$\,\dot F\ne0\,$ on $\,\bbR\smallsetminus\{0,1\}\,$ by (\ref{dtf}).) Next,
(\ref{dte}) and (\ref{gef}) yield (\ref{gta}.c), while (\ref{fet}) gives
$\,F(0)=\dot F(0)=0$, and so (\ref{gta}.d) is immediate from (\ref{sto}).
Finally, (\ref{gta}.e) and (\ref{gta}.f) are easily obtained from
(\ref{gta}.c) by evaluating it at $\,\ps=2$, or differentiating it once/twice
at $\,\ps=2\,$ and, respectively, $\,\ps=0$.

Note that (\ref{fta}) and (\ref{gef}) lead to the following special case of
(\ref{abc}):
\begin{equation}
(\ax,\bx,\cx)=(-\hs2F(\th)G(\th)\hs,\,2F(\th)\hs,\,-F(\th))\hskip9pt
\mathrm{whenever}\hskip8.5pt\th\ne1\hskip4.5pt\mathrm{and}\hskip4.5pt\ha
=\hu\nh.
\label{fug}
\end{equation}
If $\,m\,$ is even, $\,\dot G=\hs d\hskip.4ptG/d\ps\,$ is negative on
$\,(-\infty,0)\,$ and there exists a unique negative real number $\,\tz\,$
with $\,G(\tz)=0$. The following diagram uses slanted arrows to describes the
monotonicity type of $\,G\,$ (that is, the sign of $\,\dot G$) on some specific
intervals, and lists the limits of $\,G\,$ at their endpoints.
\begin{equation}
\begin{array}{lccccc}\arraycolsep9pt
\mathrm{value\ or\ limit\ at\hskip-3.2pt:}\quad&-\infty&&\tz&&0\\
[1pt]
\hline
&&&&&\\[-9.7pt]
\mathrm{for}\hskip4.5ptG\hskip4pt\mathrm{(}m\hskip4pt\mathrm{even)\hskip-3.2pt:
}&+\infty&\seam&0&\seam&E(0)
\end{array}
\label{gev}
\end{equation}
For odd $\,m\hs$, there exists unique numbers $\,\tz,\tw\in(-\infty,0)\,$ with
$\,G(\tz)=E(0)\,$ and $\,\dot G(\tw)=0$. Moreover, $\,\tz<\tw<0\,$ and
$\,E(0)=G(\tz)<G(\tw)<0$, while $\,\dot G>0\,$ on $\,(-\infty,\tw)\,$ and
$\,\dot G<0\,$ on $\,(\tw,0)$. With the same notations as in (\ref{gev}),
\begin{equation}
\begin{array}{lccccccc}\arraycolsep9pt
\mathrm{value\ or\ limit\ at\hskip-3.2pt:}\quad&-\infty&&\tz&&\tw&&0\\
[1pt]
\hline
&&&&&&&\\[-9.7pt]
\mathrm{for}\hskip4.5ptG\hskip4pt\mathrm{(}m\hskip4pt\mathrm{odd)\hskip-3.2pt:}
&-\infty&\neam&E(0)&\neam&G(\tw)&\seam&E(0)
\end{array}
\label{god}
\end{equation}
\begin{rem}\label{nalgy}
Note the analogy between $\,\tz,\tw\,$ for $\,G\,$ and $\,\zu,\wu\,$ for
$\,E\,$ in \S\ref{mopr}. Also, if $\,m\,$ is odd, $\,\tz<\zu\,$ by
(\ref{god}), since $\,F<0\,$ on $\,(-\infty,0)\,$ (see (\ref{mon})), and so
at $\,\zu\,$ we have $\,G=E-F/2>E=E(0)$.
\end{rem}

The above claims, justified in the next section, lead to a general conclusion
about the rational function $\,G\,$ given by (\ref{gef}) with a fixed integer
$\,m\ge2$, which, by (\ref{fet}), has just one pole, at $\,1$. Namely,
$\,G\ne0\,$ everywhere in $\,\bbR\smallsetminus\{1\}\,$ for odd $\,m\hs$,
while, if $\,m\,$ is even, $\,G\,$ has two zeros in
$\,\bbR\smallsetminus\{1\}$, located at $\,\tz\,$ and $\,\zj=\tz/(\tz-1)$. In
fact, as $\,G(\tw)<0$, (\ref{gev}) -- (\ref{god}) show that $\,G\ne0\,$ everywhere
in $\,(-\infty,0\hs]\,$ except at $\,\tz\,$ for even $\,m\hs$, while (\ref{elf}.c) and (\ref{fet}) give $\,E>F\ge0\,$ on $\,[\hs2,\infty)$, so that
$\,E>F/2\,$ and $\,G>0\,$ on $\,[\hs2,\infty)$. A similar assertion about the
remaining intervals now follows since, by (\ref{gta}.a), $\,G\,$ is invariant
under the involution in Lemma~\ref{uuast}, which sends $\,[\hs0,1)\,$ and
$\,(1,2\hs]\,$ onto $\,(-\infty,0\hs]\,$ and, respectively, $\,[\hs2,\infty)$.

For $\,m=3\,$ we have $\,2\tz=\hs1-\sqrt3-(2\sqrt{3\,})^{1/2}$,
$\,\,6\tw=\hs2-\sqrt{10\,}-(8\sqrt{10\,}-10)^{1/2}$. Thus,
$\,\tz\approx-\hs1.2966\,$ and $\,\tw\approx-\hs0.8456$. In fact,
$\,2(\ps-1)^3\,\lj\hs G(\ps)-E(0)\hs\rj=\ps^2\lj\ps^2
-(1+\sqrt{3\,})(\ps-1)\hs\rj\hs\lj\ps^2-(1-\sqrt{3\,})(\ps-1)\hs\rj\,$
by (\ref{gef}) and (\ref{fet}), and, similarly, since
$\,E(\ps)=\ps^3+\ps^2-2\,$ (cf.\ Example~\ref{zswtr}), we obtain
$\,6(\ps-1)^4\,\dot G(\ps)\,=\,
\ps(\ps-2)\hs\lj3\ps^2-(2+\sqrt{10\,})(\ps-1)\hs\rj\hs
\lj3\ps^2-(2-\sqrt{10\,})(\ps-1)\hs\rj$.

\section{Proofs of \hs{\rm(\ref{gev})} -- {\rm(\ref{god})}}\label{prfg}
\setcounter{equation}{0}
According to (\ref{mon}), $\,(-1)^mF(\ps)\to\infty\,$ and
$\,2\hs G(\ps)/F(\ps)\to1\,$ as $\,\ps\to-\infty$. (Note that
$\,2\hs G/F=-\hs1+2\hs E/F$, cf.\ (\ref{gef}).) Therefore, by (\ref{dtf}) and
(\ref{gef}) with $\,(\,\,)\dot{\,}=\,d/d\ps$,
\begin{equation}
\begin{array}{rl}\arraycolsep9pt
\mathrm{i)}&
(-1)^mG(\ps)\,\to\,\infty\qquad\mathrm{as}\quad\ps\,\to\,-\infty\,,\\
\mathrm{ii)}&
2\hs\dot G/\dot F\,=\,-\hs1\,+\,2\hs\dot E/\dot F\qquad\mathrm{and}\quad
(-1)^m\dot F\,<\,0\quad\mathrm{on}\quad(-\infty,0)\,.\hskip20pt
\end{array}
\label{mgt}
\end{equation}
If $\,m\,$ is even, Remark~\ref{dotef} and (\ref{mgt}.ii) give
$\,2\hs\dot G/\dot F>1$, and hence $\,2\hs\dot G<\dot F<0$, on
$\,(-\infty,0)$. As $\,G(0)=E(0)<0\,$ by (\ref{sto}), this and (\ref{mgt}.i)
show that $\,\tz\,$ exists, is unique, and (\ref{gev}) holds.

Now let $\,m\,$ be odd. By (\ref{mgt}.ii) and Remark~\ref{dotef},
$\,2\hs\dot G/\dot F\,$ decreases on $\,(-\infty,0)\,$ from $\,1\,$ to
$\,-\infty$, assuming the value $\,0\,$ at a unique $\,\tw<0$. Since
$\,\dot F>0\,$ on $\,(-\infty,0)\,$ (see (\ref{mgt}.ii)), we have $\,\dot G>0\,$
on $\,(-\infty,\tw)\,$ and $\,\dot G<0\,$ on $\,(\tw,0)$. Thus $\,G(\tw)\,$ is
the maximum value of $\,G\,$ in $\,(-\infty,0\hs]$, and so
$\,G(\tw)>G(0)=E(0)$. As $\,G\,$ increases on $\,(-\infty,\tw)\,$ from
$\,-\infty\,$ (see (\ref{mgt}.i)) to $\,G(\tw)$, it assumes the intermediate
value $\,E(0)\,$ exactly once in $\,(-\infty,\tw)$, and not at all in
$\,[\tw,0)\,$ (where it decreases to the limit $\,E(0)$), which gives the
existence and uniqueness of $\,\tz\,$ and the relation $\,\tz<\tw$. Finally,
on $\,(-\infty,0)\,$ we have $\,F<0\,$ and $\,E/F>1/2\,$ by (\ref{mon}), and
hence $\,G=E-F/2<0$, so that $\,G(\tw)<0$.

\section{Positivity and type \hs{\rm(\ref{trc}.b)}}\label{ptpb}
\setcounter{equation}{0}
In \S\ref{sbto} we classified all those pairs $\,Q,I\hs$ with (\ref{qdq}.b) and (\ref{qdq}.e)
that are of type (\ref{trc}.b). We will now determine which of them also satisfy
the positivity condition of \S\ref{posi} or, equivalently, all the remaining
conditions in {\rm(\ref{qdq})}.
\begin{prop}\label{pyuzw}Given
$\,\th\in\bbR\smallsetminus\{0,1,2\}$, let\/ $\,\ha=\hu$ with\/
$\,\hu=\th/(\th-1)$, cf.\ {\rm(\ref{ast})}. Also, let\/ $\,\jv\,$ be the space
{\rm(\ref{spv})} for a fixed integer\/ $\,m\ge2$, and let
$\,Q\in\jv\smallsetminus\{0\}\,$ be a function with {\rm(\ref{qdq}.b)} on the closed
interval\/ $\,I\hs$ with the endpoints\/ $\,\th,\ha$. By\/ {\rm Lemma~\ref{onedi}}, such\/ $\,Q\,$ exists and is unique up to a factor. Then, the pair\/
$\,Q,I\hs$ satisfies all five conditions listed in {\rm(\ref{qdq})} if and only if
\begin{enumerate}
\item[a)] $m\,$ is even, or
\item[b)] $m\,$ is odd, while\/ $\,\th\notin[\tz,\tw]\,$ and\/
$\,\ha\notin[\tz,\tw]$, with\/ $\,\tz,\tw\,$ as in {\rm\S\ref{anth}}.
\end{enumerate}
\end{prop}
\begin{proof}Relation $\,\ha=\hu$ amounts to (\ref{trc}.b), and hence
implies (\ref{qdq}.e) (Theorem~\ref{facto}). Therefore we just need to show that (a) or (b)
holds if and only if $\,Q,I\hs$ satisfy (\ref{qdq}.a), (\ref{qdq}.c), (\ref{qdq}.d), that is,
the positivity condition of \S\ref{posi}.

First, let one of $\,\th,\ha\,$ lie in $\,(1,2)\cup(2,\infty)$. Then so does
the other (Lemma~\ref{uuast}), and the positivity condition for $\,Q,I\hs$
follows from (i) in \S\ref{posi}, while, at the same time, we have (a) or (b), as
$\,[\tz,\tw]\subset(-\infty,0)\,$ when $\,m\,$ is odd (see \S\ref{anth}).

Thus, it suffices to prove our assertion in the remaining case where
$\,\th\,$ and $\,\ha=\hu$ both lie in $\,(-\infty,0)\cup(0,1)$, that is (cf.\
Lemma~\ref{uuast}), one of them is in $\,(-\infty,0)\,$ and the other in
$\,(0,1)$. The hypotheses of Proposition~\ref{dtefu} now are satisfied; also,
$\,F(\th)\ne0\,$ by (\ref{fet}), so that (\ref{fug}) gives $\,\bx\ne0$,
$\,\cx/\bx=-\hs1/2\,$ and $\,\ax/\bx=-\hs G(\th)$. Using (\ref{fug}), we may
now rephrase conditions (i), (ii) of Proposition~\ref{dtefu} as follows:
\begin{enumerate}
\item $\,\dot G(\th)/\dot F(\th)<0<\dot G(\ha)/\dot F(\ha)$, provided that
$\,\th\,$ and $\,\ha=\hu$ have been switched, if necessary, so that
$\,(-1)^m(u-v)>0$.
\item $\,G(\th)<E(0)$.
\end{enumerate}
We will show that (a) or (b) holds if and only if one of these conditions (i),
(ii) is satisfied, thus obtaining our assertion as a direct consequence of
Proposition~\ref{dtefu}.

Let us first suppose that neither (a) nor (b) holds, that is, $\,m\,$ is odd
and one of $\,\th,\ha\,$ lies in $\,[\tz,\tw]$. As $\,\tw<0$, ordering
$\,\th,\ha\,$ as in (i) (so that $\,\th<\ha$), we obtain
$\,\th\in[\tz,\tw]\,$ and $\,\ha\in(0,1)$. Since $\,G(\th)\ge E(0)\,$ (see
(\ref{god})), condition (ii) cannot be satisfied. Also, $\,\dot G(\th)\ge0\,$
(see (\ref{god})); as $\,\dot F>0\,$ on $\,(-\infty,0)\,$ for odd $\,m\,$ by (\ref{dtf}), this gives $\,\dot G(\th)/\dot F(\th)\ge0$, and (i) fails as well.
Thus, if (i) or (ii) is satisfied, we must have (a) or (b).

To prove the converse, let us now assume (a) or (b) and order $\,\th,\ha\,$ as
in (i). If $\,m\,$ is even, we thus have $\,\hu=\ha<0<\th<1\,$ and
$\,\dot G(\ha)<0\,$ (cf.\ (\ref{gev})). Hence, as (\ref{dtf}) yields
$\,\dot F<0\,$ on$\,(-\infty,0)\cup(0,1)$, (\ref{gta}.b) implies the double
inequality in (i). On the other hand, if $\,m\,$ is odd, we have
$\,\th<0<\ha<1\,$ and, as $\,\tw<0$, condition (b) leads to two possible
cases: $\,\th<\tz$, or $\,\tw<\th<0$. If $\,\th<\tz$, (\ref{god}) gives
$\,G(\th)<E(0)$, that is, (ii); while, if $\,\tw<\th<0$, (\ref{god}) shows that
$\,\dot G(\th)<0$, which, combined with (\ref{gta}.b) and the fact that
$\,\dot F>0\,$ on $\,(-\infty,0)\,$ by (\ref{dtf}), again gives the double
inequality in (i). This completes the proof.
\end{proof}

\section{More on the polynomial $\,T$}\label{mort}
\setcounter{equation}{0}
Given an integer $\,m\ge2$, let $\,\av,\bv,\fj,\fy\,$ be the polynomials in
$\,\th,\ha\,$ with
\begin{equation}
\begin{array}{rl}\arraycolsep9pt
\mathrm{a)}&
\av\,=\,(\th-1)^m(\ha-2)\ha^{2m-1}E(\th)\,
-\,(\ha-1)^m(\th-2)\th^{2m-1}E(\ha)\hs,\\
\mathrm{b)}&
\bv\,=(\ha-1)^m(\th-2)\th^{2m-1}\,-\,(\th-1)^m(\ha-2)\ha^{2m-1}\hs,\\
\mathrm{c)}&
\fj\,=\,\th+\ha-2\,,\qquad\quad\fy\,=\,m\hs\th\ha-(2m-1)\hs\fj\hs.\end{array}
\label{aph}
\end{equation}
for $\,E\,$ as in (\ref{fet}). Thus, with $\,\ax,\bx\,$ depending on
$\,\th,\ha\in\bbR\smallsetminus\{1\}\,$ as in (\ref{abc}),
\begin{equation}
\av\,=\,(\th-1)^m(\ha-1)^m\ax\,,\qquad\bv\,=\,(\th-1)^m(\ha-1)^m\bx\,.
\label{auv}
\end{equation}
Lemmas~\ref{polyt} and~\ref{cpiuv} imply that the product of
$\,\th(\th-2)\ha(\ha-2)\,$ and the right-hand side of (\ref{umv}) equals
$\,\,(\th-1)^m(\ha-1)^m$ times $\,2\hs(\th\ha-\th-\ha)\,$ times the expression
(\ref{muv}). Dividing by $\,2\hs(\th\ha-\th-\ha)\,$ and using (\ref{umv}),
(\ref{aph}), (\ref{auv}), we obtain, for the polynomial $\,T\,$ given by
(\ref{umv}),
\begin{equation}
\begin{array}{rl}\arraycolsep9pt
\mathrm{i)}&
(\ha-\th)^3\th(\th-2)\ha(\ha-2)\hs T(\th,\ha)\,
=\,\fy\av\,+\,\cj\hs\bv\fj\hs,\qquad\mathrm{where}\hskip30pt\\
\mathrm{ii)}&
\cj\,=\,(2m-1)\hs\vs(0)\,>\,0\hs,\hskip14pt\mathrm{with}\hskip10pt\vs(0)
\hskip10pt\mathrm{as\ in\ (\ref{sto}.i)}.\end{array}
\label{vuu}
\end{equation}
(This $\,\cj\,$ is not related to the symbol $\,\y\,$ in \S\ref{intr} --
\S\ref{rati}.) For $\,\ha=1$, (\ref{aph}) gives $\,\av=-\hs(\th-1)^mE(\th)$,
$\,\bv=(\th-1)^m\,$ and $\,\fj=\th-1$, $\,\fy=(1-m)\th+2m-1$, so that
(\ref{vuu}.i) with $\,\ha=1\,$ becomes
\begin{equation}
\begin{array}{rl}\arraycolsep9pt
\mathrm{a)}&
\th(\th-2)\hs T(\th,1)\,=\,(\th-1)^{m-2}\pj(\th)\qquad\mathrm{for\ all}
\quad\th\in\bbR\hs,\quad\mathrm{where}\hskip10pt\\
\mathrm{b)}&
\pj(\th)\,=\,[(m-1)\th-2m+1]\hs\vs(\th)\,+\,(2m-1)\hs\vs(0)\quad\mathrm{for}
\hskip7pt\th\in\bbR\hs.\end{array}
\label{uut}
\end{equation}
Also, with $\,F,\vs,G,T\hs$ given by (\ref{fet}), (\ref{esg}), (\ref{gef}), (\ref{umv}) for an integer $\,m\ge2$,
\begin{equation}
\begin{array}{rl}\arraycolsep9pt
\mathrm{a)}&
\th^5(\th-2)^5\,T(\th,\hu)\,=\,2(\th-1)^2\hs F(\th)\hs\xo(\th)\hskip11pt
\mathrm{for\ all}\hskip7pt\th\ne1\hs,\hskip7pt\mathrm{where}\\
\mathrm{b)}&
\xo(\th)\nh=\nh[(m\nh-\nh1)\th^2\hskip-2.2pt-\nh2(2m\nh-\nh1)(\th\nh-\nh1)]\hs
G(\th)+(2m\nh-\nh1)(\th^2\hskip-2pt-\nh2\th\nh+\nh2)\vs(0)\end{array}
\label{tuu}
\end{equation}
and $\,\hu=\th/(\th-1)\,$ for $\,\th\ne1$.
In fact, setting $\,\ha=\hu$ in (\ref{aph}.c) we get
$\,(\th-1)\hs\fj=(\th-1)^2+1\,$ and $\,(1-\th)\fy=(m-1)\th^2-2(2m-1)(\th-1)$.
Since $\,(\th-1)(\ha-\th)=(\th-1)^2\ha(\ha-2)=\th(2-\th)\,$ whenever
$\,\ha=\hu$, we obtain (\ref{tuu}.a) by setting $\,\ha=\hu$ in (\ref{vuu}.i) and
noting that (\ref{auv}) with $\,\ha=\hu$ reads $\,\av=\ax\,$ and $\,\bv=\bx\,$
(cf.\ (\ref{ast}.b)), while $\,\ax,\bx\,$ with $\,\ha=\hu$ are in turn given
by (\ref{fug}).

Finally, for any integer $\,m\ge2$, applying $\,\,\partial/\partial\th\,$ to
(\ref{vuu}.i) at $\,\th=0\,$ and using (\ref{aph}) with $\,\dot E(0)=0\,$ (see
(\ref{sto}.iii)), we obtain
\begin{equation}
2\hs T(0,\ha)=(-1)^mm\hs\vs(0)\hs\ha^{2m-4}\hskip6pt\mathrm{for\ all}
\hskip6pt\ha\in\bbR\hs,\hskip4.5pt\mathrm{with}\hskip5.5pt\vs(0)\hskip5.5pt
\mathrm{as\ in\ (\ref{sto}.i).}\label{tov}
\end{equation}

\section{The sign of $\,T\,$ on a specific curve}\label{sgnt}
\setcounter{equation}{0}
Let $\,T\,$ be the polynomial in $\,\th,\ha\,$ given by (\ref{umv}) for any
fixed integer $\,m\ge2$. In this section we describe the behavior of
$\,\,\sgn\,T\,$ on the intersection of the half-plane $\,\th<0\,$ with the
hyperbola $\,\mathcal{H}\,$ given by $\,\th\ha=\th+\ha\,$ (that is,
$\,\ha=\hu$, cf.\ (\ref{ast})). This will allow us, later in \S\ref{tzer}, to
draw important conclusions about zeros of $\,T\,$ in the region
$\,\th<0<\ha\le1$.

We need the next result to prove Lemmas~\ref{xxast} and~\ref{fdpsi}. See also Example~\ref{xymtr}.
\begin{lem}\label{cfpsi}Given\/ $\,\cj\in\bbR\smallsetminus\{0\}\,$ and\/
$\,C^1$ functions\/ $\,\av,\fy,\ly,\bda,\my,\vi\,$ on an open interval, let\/
$\,(\,\,)'$ be the derivative operator followed by multiplication by
some fixed\/ $\,C^1$ function. If\/ $\,\av'=\bda\av-\cj\hs\my\,$ and\/
$\,\vi=\fy\av+\cj\hs\ly$, then, at every point where\/ $\,\vi=0$,
\begin{enumerate}
\item[a)] $\,\cj^{-1}\fy\vi'=(\ly\hs'-\bda\ly-\fy\my)\hs\fy-\ly\fy'$,
\item[b)] $\,\cj^{-1}\fy\vi'=(\fy\fj'-\fj\fy')\hs\bv-(\fj\ny+\fy\my)\hs\fy\,$
provided that, in addition, $\,\ly=\bv\fj\,$ and\/ $\,\bv\hs'=\bv\bda-\ny\,$
for some\/ $\,C^1$ functions\/ $\,\fj,\bv,\ny$.
\end{enumerate}
\end{lem}
In fact, any point with $\,\vi=0\,$ has $\,\fy\av=-\hs\cj\hs\ly$, and hence
$\,\fy\av'=\bda\fy\av-\cj\fy\my=-\hs\cj\hs(\bda\ly+\fy\my)$, so that (a) is
immediate as $\,\vi=\fy\av+\cj\hs\ly\,$ and
$\,\fy(\fy\av)'=(\fy\av')\fy+(\fy\av)\fy'$. To obtain (b), it now suffices to
replace $\,\ly\,$ and $\,\ly\hs'$ in (a) by $\,\bv\fj\,$ and, respectively,
$\,\bv\fj'+\bv\hs'\nh\fj=\fj\hs'\nh\bv+(\bv\bda-\ny)\fj$.
\begin{lem}\label{xxast}Given an integer\/ $\,m\ge2$, let\/ $\,T,\xo\,$ be the
polynomial and rational function with {\rm(\ref{umv})} and\/
{\rm(\ref{tuu}.b)}. Also, let\/ $\,\hu=\th/(\th-1)\,$ for\/ $\,\th\ne1$, cf.\
{\rm(\ref{ast})}.
\begin{enumerate}
\item[a)] If\/ $\,m\,$ is even, $\,T(\th,\hu)>0\,$ for every\/ $\,\th<0$.
\item[b)] If\/ $\,m\,$ is odd, equation $\,T(\xj,\xj^*)=0\,$ has a unique
negative real solution $\,\xj$. With this\/ $\,\xj\hs$ we have\/
$\,T(\th,\hu)>0\,$ for all\/ $\,\th\in(-\infty,\xj)\,$ and\/
$\,T(\th,\hu)<0\,$ for all\/ $\,\th\in(\xj,0)$, as well as\/ $\,\xj<\tz$,
where\/ $\,\tz\,$ is defined as in {\rm\S\ref{anth}}.
\item[c)] If\/ $\,m\,$ is odd, $\,\,\sgn\,(\th-\xj)=\hs\sgn\,\xo(\th)\,$
for\/ $\,\xj\,$ as in {\rm(b)}, $\,\hs\sgn\hs\,$ as in {\rm\S\ref{rati}}, and
any $\th\in(-\infty,0\hs]$.
\item[d)] $d\hs[T(\th,\hu)]/d\th<0\,$ at every\/ $\,\th<0\,$ with\/
$\,T(\th,\hu)=0$.
\end{enumerate}
\end{lem}
\begin{proof}By (\ref{gta}.d), $\,\xo(0)=\dot\xo(0)=0\,$ and
$\,(2m-3)\hs\ddot\xo(0)=8m(m-1)\hs\vs(0)>0$, where
$\,(\,\,)\dot{\,}=\,d/d\ps$. Thus, $\,\xo(\ps)>0\,$ for $\,\ps<0\,$ near
$\,0$. Also, $\,(-1)^m\xo(\ps)\to\infty\,$ as $\,\ps\to-\infty$, since, by (\ref{mgt}.i), $\,(-1)^m\xo(\ps)/\ps^2\to\infty\,$ as $\,\ps\to-\infty$.

In view of (\ref{gta}.c), the assumptions of Lemma~\ref{cfpsi} are satisfied by
$\,\av=G(\ps)\,$ with (\ref{gef}), $\,\fy=(m-1)\ps^2-2(2m-1)(\ps-1)$,
$\,\ly=\ps^2-2\ps+2$, $\,\bda=\varLambda(\ps)\,$ with (\ref{dtf}.b),
$\,\my=2(\ps-1)$, $\,\vi=\xo(\ps)$, the constant (\ref{vuu}.ii), and
$\,(\,\,)'=\ps(\ps-1)(\ps-2)\,d/d\ps$. Also,
$\,\fy=\bda-\ps^2$ (by (\ref{dtf}.b)) and $\,\ly=\ps^2-\my$, so that
$\,\bda\ly+\fy\my=\bda\ps^2-\bda\my+\bda\my-\my\ps^2
=(\bda-\my)\ps^2=m\ps^2(\ps-2)^2$ and
$\,\ly\hs'-\bda\ly-\fy\my=-\hs\ps(\ps-2)\hs[(m-2)\ps(\ps-2)-2]$. Moreover,
$\,\ly\fy'=2\ps(\ps-1)(\ps-2)(\ps^2-2\ps+2)[(m-1)\ps-2m+1]$. The right-hand
side in Lemma~\ref{cfpsi}(a) thus equals $\,-\hs m(m-1)\ps^2(\ps-2)^4$. Using
Lemma~\ref{cfpsi}(a), with both sides divided by $\,\ps(\ps-2)$, we see that
$\,(\ps-1)\fy\hs\dot\xo(\ps)=-\hs m(m-1)\cj\hs\ps(\ps-2)^3$ at every
$\,\ps\in\bbR\smallsetminus\{0,1,2\}\,$ at which $\,\xo(\ps)=0$. This
gives $\,\dot\xo(\ps)>0\,$ for every $\,\ps<0\,$ with $\,\xo(\ps)=0$, as the
quadratic polynomial $\,\fy$, having a positive value and a negative
derivative at $\,\ps=0$, must be positive for all $\,\ps<0$. Remark~\ref{zerpo}
for $\,\varPhi=\xo\,$ and $\,\mathcal{I}=(-\infty,0)$, combined with the signs of
$\,\xo\,$ near $\,-\infty\,$ and $\,0\,$ determined in the last paragraph, now
gives (c) with some (unique) $\,\xj<0$, as well as $\,\xo>0\,$ on
$\,(-\infty,0)\,$ if $\,m\,$ is even. Since $\,\xo(\xj)=0$, the definition of
$\,\xo\,$ leads to a rational expression for $\,G(\xj)\,$ in terms of $\,\xj$,
which easily shows that $\,G(\xj)<E(0)=-\hs\vs(0)\,$ (cf.\ (\ref{sto}.ii)), so
that $\,\xj<\tz\,$ by (\ref{god}). Also, as $\,\,\sgn\,F=(-1)^m$ on
$\,(-\infty,0)$, cf.\ (\ref{mon}), relation (\ref{tuu}.a) yields
$\,\hs\sgn\,T(\th,\hu)=(-1)^m\hskip2pt\sgn\,\xo(\th)\,$ for $\,\th<0$, with
$\,\hs\sgn\,$ as in \S\ref{rati}. Hence (a), (b) and (d) follow, completing the
proof.
\end{proof}

\section{The differential of $\,T\,$ at points where
$\,T\hskip-1.2pt=0$}\label{dfto}
\setcounter{equation}{0}
\begin{lem}\label{fdpsi}Given an integer\/ $\,m\ge2$, let\/
$\,T,\bv,\fj,\fy,\vi,\zy,\tj\,$ be the polynomials in $\,\th,\ha\,$
given, respectively, by {\rm(\ref{umv})}, {\rm(\ref{aph}.b)},
{\rm(\ref{aph}.c)}, formula
\begin{equation}
\vi\,=\,(\ha-\th)^3\th(\th-2)\ha(\ha-2)\hs T(\th,\ha)\,,\label{psv}
\end{equation}
and\/ $\,\zy=(\th-1)^{m-1}(\ha-2)\ha^{2m-1}$,
$\,\tj=(\ha-1)^{m-1}(\th-2)\th^{2m-1}$. 
For the constant\/ $\,\cj\,$ with {\rm(\ref{vuu}.ii)} we then have, at every
point\/ $\,(\th,\ha)\,$ at which\/ $\,\vi=0$,
\begin{equation}
\begin{array}{ccc}\arraycolsep9pt
(m\cj)^{-1}\fy\,\partial\vi/\partial\th&=&(\th-\ha)\hs\fy\zy\,
-\,\ha(\ha-2)\hs\bv\hskip9pt\mathrm{unless}\quad\th(\th-2)=0\,,\\
(m\cj)^{-1}\fy\,\partial\vi/\partial\ha&=&(\th-\ha)\hs\fy\tj\,
-\,\th(\th-2)\hs\bv\hskip9pt\mathrm{unless}\quad\ha(\ha-2)=0\,.
\end{array}
\label{mcf}
\end{equation}
\end{lem}
\begin{proof}The assumptions of Lemma~\ref{cfpsi}(b) hold for
$\,\bv,\fj,\fy,\vi,\cj\,$ chosen as above, $\,\av\,$ as in (\ref{aph}.a),
$\,\ly=\bv\fj$, $\,\bda=2\hs(m\th-2m+1)$, $\,\my=2(\th-1)\zy$,
$\,\ny=-\hs\varLambda(\th)\zy$, with $\,\varLambda\,$ given by (\ref{dtf}.b),
and $\,(\,\,)'=\th(\th-2)\,\partial/\partial\th$, where the variable $\,\ha\,$
is fixed.

This is clear as $\,\vi=\fy\av+\cj\hs\fj\bv\,$ by (\ref{psv}) and
(\ref{vuu}.i), while $\,\av'=\bda\av-\cj\hs\my\,$ and
$\,\bv\hs'=\bv\bda-\ny\,$ in view of (\ref{dte}.ii) with
$\,\varLambda(\th)+m\hs\th(\th-2)=2(\th-1)(m\th-2m+1)$, which is immediate
from (\ref{dtf}.b).

Now $\,\fj'=\th(\th-2)$, $\,\fy'=\th(\th-2)\hs\hz(\ha)$, and so
$\,\fy\fj'-\fj\fy'\,$ equals $\,\th(\th-2)\,$ times $\,\fy-\fj\hs\hz(\ha)
=m\hs\th\ha-[2m-1+\hs\hz(\ha)\hs]\hs\fj=m\ha\hs(\th-\fj)
=-\hs m\hs\ha(\ha-2)\,$ with
$\,\hz(\ha)=m\hs \ha-2m+1$, that is,
$\,\fy\fj'-\fj\fy'=-\hs m\th\ha(\th-2)(\ha-2)$.

Finally,
$\,(2m-1)\my-\ny=[\hs2(2m-1)(\th-1)+\varLambda(\th)\hs]\hs\zy=m\th^2\zy\,$
(cf.\ (\ref{dtf}.b)). Therefore,
$\,\fj\ny+\fy\my=\fj\ny+[m\hs\th\ha-(2m-1)\fj]\hs\my=-\hs[\hs(2m-1)\my
-\ny\hs]\,\fj+m\hs\th\ha\my=-\hs m\th^2\fj\hs\zy+m\hs\th\ha\my=
-\hs m\hs\th(\th-2)(\th-\ha)\zy$. Dividing both sides of the equality in
Lemma~\ref{cfpsi}(b) by $\,m\hs\th(\th-2)$, we arrive at the first relation in
(\ref{mcf}). The second one is obtained by evaluating the first at the point
$\,(\ha,\th)\,$ rather than $\,(\th,\ha)$. In fact, switching $\,\th,\ha\,$
causes $\,\zy\,$ to be replaced by $\,\tj$, while $\,\fy\,$ is symmetric in
$\,\th,\ha\,$ and $\,\bv,\vi\,$ are anti-sym\-met\-ric, so that
the values of $\,\partial\hs\vi\nh/\partial\th\,$ at $\,(\ha,\th)\,$ and
$\,\partial\hs\vi\nh/\partial\ha\,$ at $\,(\th,\ha)\,$ are mutually opposite.
This completes the proof.
\end{proof}
\begin{lem}\label{polyz}For any integer\/ $\,m\ge2$, let\/ $\,\zi,\vx\,$ be
the polynomials in real variables $\,\th,\ha\,$ with\/ $\,\zi=1\,$ if\/
$\,m\,$ is even, $\,\zi=\hs\th\ha(\th+\ha)-(\th^2+\ha^2)\,$ if\/ $\,m\,$ is
odd, and $\,\vx=[(\th-1)\ha^2]^{m-1}-[(\ha-1)\th^2]^{m-1}$ for all\/ $\,m\hs$.
Then $\,\vx=(\ha-\th)\hs(\th\ha-\th-\ha)\zi\vy\,$ for some
polynomial\/ $\,\vy\,$ in $\,\th,\ha\,$ such that\/
$\,\vy>0\,$ on $\,\rto\smallsetminus\{(0\hs,\nh0),(1\hs,\nh1)\}$.
\end{lem}
This is clear since
$\,\xi^{2n}-\eta^{2n}=(\xi-\eta)(\xi+\eta)(\xi^{2n-2}+\xi^{2n-4}\eta^2
+\ldots+\eta^{2n-2})\,$ and $\,\xi^{m-1}-\eta^{m-1}
=(\xi-\eta)(\xi^{m-2}+\xi^{m-3}\eta+\ldots+\eta^{m-2})\,$ for an integer
$\,n\ge1\,$ or an even integer $\,m\ge2$, and $\,\xi,\eta\in\bbR\hs$. In
both cases, the last factor is positive unless $\,\xi=\eta=0.$ (In fact, if
$\,m\,$ is even and $\,\eta\ne0$, then $\,\eta\,$ is both a simple root and
the unique real root of the polynomial
$\,\xi\,\mapsto\,\xi^{m-1}\hskip-1pt-\eta^{m-1}\hskip-1.3pt$.) Our assertion
now follows if we use $\,n=(m-1)/2\,$ when $\,m\,$ is odd, and set
$\,\xi=(\th-1)\ha^2$, $\,\eta=(\ha-1)\th^2$, so that
$\,\xi-\eta=(\ha-\th)(\th\ha-\th-\ha)\,$ and
$\,\xi+\eta=\th\ha(\th+\ha)-(\th^2+\ha^2)$.

For $\,\fy,\zi\,$ depending on $\,(\th,\ha)\,$ as in (\ref{aph}.c) and
Lemma~\ref{polyz} with a fixed $\,m\ge2$,
\begin{equation}
\fy\,\ge\,m\,>\,0\hskip12pt\mathrm{and}\hskip12pt(-1)^m\zi\,\ge\,0
\hskip12pt\mathrm{whenever}\hskip12pt\th-1\le0<\ha\le1\hs.\label{fmo}
\end{equation}
In fact, $\,\fy=m\hs\tilde\th\tilde\ha+(m-1)(\tilde\th+\tilde\ha)+m\,$ for
$\,\tilde\th=1-\th\,$ and $\,\tilde\ha=1-\ha$, so that $\,\fy\ge m\,$ as
$\,\tilde\th,\tilde\ha\ge0$. That $\,(-1)^m\zi\ge0\,$ follows if one adds
up the inequalities obtained from multiplying the relations $\,\ha\le1\,$ and
$\,\th\ge1\,$ by $\,\th^2$ and, respectively, $\,\ha^2$.
\begin{lem}\label{sgndt}Given an integer\/ $\,m\ge3$, for\/ $\,\hs\sgn\hs\,$
and\/ $\,T\,$ as in {\rm\S\ref{rati}} and\/ {\rm Lemma~\ref{polyt}} we have,
at every point\/ $\,(\th,\ha)\in\rto$ at which\/ $\,\th<0<\ha\le1\,$ and\/
$\,T(\th,\ha)=0$,
\begin{equation}
(-1)^m\,\sgn\hs\,d_{\hs\mathbf{w}}T\,=\,\,\sgn\hs(\th\ha-\th-\ha)\qquad
\mathrm{and}\qquad d\hskip1.3ptT\,\ne\,0\,,\label{sdt}
\end{equation}
$d_{\hs\mathbf{w}}$ being the directional derivative for the vector field\/
$\,\hs\mathbf{w}\,=(\th-1,\hs\ha-1)\,$ on $\,\rto\nh$.
\end{lem}
\begin{proof}With $\,\cj,\fy,\vi\,$ as in (\ref{vuu}.ii), (\ref{aph}.c) and (\ref{psv}), let $\,\vrp\,$ stand for the value of
$\,(m\cj)^{-1}\fy\hs d_{\hs\mathbf{w}}\vi\,$ at any fixed $\,(\th,\ha)\,$
with $\,\th<0<\ha\le1\,$ and $\,T(\th,\ha)=0$. Thus,
$\,\hs d_{\hs\mathbf{w}}$ applied to (\ref{psv}) gives
$\,\hs\sgn\hs\,d_{\hs\mathbf{w}}T=\,-\,\sgn\hs\,\vrp\,$ at
$\,(\th,\ha)$, as $\,\cj\fy>0\,$ (see (\ref{fmo}). Also,
$\,\fy=(\th-2)(m\hs\ha-2m+1)+\ha$. Hence $\,\vrp$, being, by (\ref{mcf}),
the difference of $\,(\th-\ha)\hs\lj(\th-1)^m(\ha-2)\ha^{2m-1}
+(\ha-1)^m(\th-2)\th^{2m-1}\rj\hs\fy\,$ and
$\,\th(\th-2)(\ha-1)\hs\bv+\ha(\ha-2)(\th-1)\hs\bv\,$ at $\,(\th,\ha)$, must
equal the sum of the two expressions
\[
\begin{array}{l}
(\th-\ha)\hs\lj(\th-2)(m\hs\ha-2m+1)+\ha\rj\hs
\lj(\th-1)^m(\ha-2)\ha^{2m-1}\hskip-1pt+(\ha-1)^m(\th-2)\th^{2m-1}\rj\hs,\\
\lj\th(\th-2)(\ha-1)+\ha(\ha-2)(\th-1)\rj\hs
\lj(\th-1)^m(\ha-2)\ha^{2m-1}\hskip-1pt-(\ha-1)^m(\th-2)\th^{2m-1}\rj\hs.
\end{array}
\]
In both displayed products, either of the two factors enclosed in square
brackets is a sum of two polynomials; of these eight polynomials,
four are manifestly divisible by $\,\th-2$. Direct multiplications in both
displayed lines, performed to evaluate $\,\vrp$, thus give rise to
eight product terms, of which only two,
$\,\ha(\ha-2)(\th-1)\,$ times $\,(\th-1)^m(\ha-2)\ha^{2m-1}$ and
$\,(\th-\ha)\hs\ha\,$ times $\,(\th-1)^m(\ha-2)\ha^{2m-1}\nh$, fail
to explicitly contain $\,\th-2\,$ as a factor; their sum, however, is
$\,(\th-2)(\ha-1)(\th-1)^m(\ha-2)\ha^{2m}$. Thus, the polynomial
in $\,\th,\ha\,$ representing the value $\,\vrp\,$ is divisible by
$\,\th-2$, and the calculation just outlined gives
\begin{equation}
\begin{array}{rcl}\arraycolsep9pt
\vrp/(\th-2)&=&\th(\th-1)^m(\ha-1)(\ha-2)\ha^{2m-1}\\
&-&(\th\ha-\th-\ha)\hs\lj\th\hs+\hs(\ha-2)\rj\hs\th^{2m-1}(\ha-1)^m\\
&+&(\th-\ha)\hs\lj(m\hs\th-2m+1)(\ha-2)\,+\,\th\rj\hs\th^{2m-1}(\ha-1)^m\\
&+&(\ha-2)\hs\lj(m\hs\ha-2m+1)\th-(m-1)(\ha-2)\ha\rj\hs(\th-1)^m\ha^{2m-1}\nh.
\end{array}
\label{yps}
\end{equation}
(The second line of (\ref{yps}) is the result of combining two of the four
terms containing the factor $\,\th^{2m-1}(\ha-1)^m$ with the aid of
(\ref{uuv}).) The right-hand side can, as before, be rewritten as the sum of
polynomial product terms, with each of the three square brackets contributing
to two of them; this time, there are {\it seven\/} such terms, and five of
them are manifestly divisible by $\,\ha-2$, while the other two add up to
$\,\lj(\th+\ha-\th\ha)\hs+\hs(\th-\ha)\rj\hs\th^{2m}(\ha-1)^m$, that is,
$\,\ha-2\,$ times $\,-\hs\th^{2m+1}(\ha-1)^m$. Consequently, as a polynomial
in $\,\th,\ha$, our $\,\vrp/(\th-2)\,$ is divisible by $\,\ha-2\,$ and,
proceeding just as we did above to evaluate $\,\vrp/(\th-2)$, we obtain
\begin{equation}
\begin{array}{rcl}\arraycolsep9pt
\vrp/[(\th-2)(\ha-2)]&=&\lj(m-1)\th^2-(m+1)\th\ha-2(m-1)\th
+2m\hs\ha\rj\hs\th^{2m-1}(\ha-1)^m\\
&-&\lj(m-1)\ha^2-(m+1)\th\ha-2(m-1)\ha+2m\hs\th\rj\hs\ha^{2m-1}(\th-1)^m\hs.
\end{array}
\end{equation}
Relation $\,[(\ha-1)\th^2]^{m-1}=[(\th-1)\ha^2]^{m-1}-\vx$, for $\,\vx\,$ as
in Lemma~\ref{polyz}, allows us to replace the factor
$\,\th^{2m-1}(\ha-1)^m=[(\ha-1)\th^2]^{m-1}(\ha-1)\th\,$ with
$\,\,-\hs(\ha-1)\th\vx\,+\,(\ha-1)\th(\th-1)^{m-1}\ha^{2m-2}$, while,
$\,\vx=(\ha-\th)\hs(\th\ha-\th-\ha)\zi\vy\,$ by Lemma~\ref{polyz}. It now
easily follows that $\,\vrp\,$ is equal to
$\,(\th-2)(\ha-2)(\th-\ha)(\th\ha-\th-\ha)\,$ times
\begin{equation}
\begin{array}{l}\arraycolsep4pt
\lj(m-1)\th^2\,-\,(m+1)\th\ha\,-\,2(m-1)\th\,+\,2m\hs\ha\rj\hs(\ha-1)
\hs\th\hs\zi\vy\\
\phantom{\lj(m-1)\th^2\,-\,(m+1)\th\ha\,}
+\,(m-1)(\th+\ha-2)(\th-1)^{m-1}\ha^{2m-2}\nh.
\end{array}
\label{muq}
\end{equation}
As $\,\th<0<\ha\le1$, all four terms within the square brackets in (\ref{muq})
are positive. The signs of the other four factors
$\,(\ha-1),\hs\th,\hs\zi,\hs\vy\,$ in that line are, respectively, \hs($-\,$
or $\,0$), $\,-\hs$, \hs($(-1)^m$ or $\,0$), $\hs+\hs$, cf.\ (\ref{fmo}) and
Lemma~\ref{polyz}, while those of the four factors in the second line are
$\,+\hs,\hs-\hs,\hs(-1)^{m-1},\hs+\hs$, so that $\,\hs\sgn\hs\,$ of (\ref{muq})
equals $\,(-1)^m$, which clearly yields the first relation in (\ref{sdt}). The
second one then is immediate from the first, since, by Lemma~\ref{xxast}(d),
$\,d\hskip1.3ptT\ne0\,$ at our $\,(\th,\ha)\,$ also in the case where
$\,\ha=\hu$, that is, $\,\th\ha-\th-\ha=0$. This completes the proof.
\end{proof}
\begin{rem}\label{submf}
By (\ref{sdt}), $\,0\,$ is a regular value
of $\,T\,$ in $\,\mathcal{R}=(-\infty,0)\times(0,1\hs]$, that is, in the region
$\,\mathcal{R}\,$ in the $\,\th\ha\hs$-plane $\,\rto$ given by $\,\th<0<\ha\le1$.
Thus, if the set of zeros of $\,T\,$ in $\,\mathcal{R}\,$ is nonempty, its
connected components are one\diml\/ real-analytic submanifolds of $\,\mathcal{R}$,
possibly with boundary. When $\,m\,$ is odd, that set of zeros is nonempty, as
it contains $\,(\xj,\xj^*)\,$ (see Lemma~\ref{xxast}(b)).
\end{rem}

\section{An analytic curve segment with $\,T\hskip-1.2pt=0$}\label{cseg}
\setcounter{equation}{0}
The meaning of the symbol $\,r\hs$ in this section is not related to its use
in \S\ref{main}.
\begin{lem}\label{dydro}Let a real number\/ $\,\qx_*\in[-\hs1,0)\,$ and a\/
$\,C^1$ function $\,Y$ of the real variables\/ $\,\qx,\rx$, defined on
the square\/ $\,\mathcal{S}=(-\hs1,0)\times[\hs0,1]$, satisfy the conditions
\begin{enumerate}
\item $\partial Y/\partial\rx<0\,$ at every interior point of\/
$\,\mathcal{S}\,$ at which\/ $\,Y=0\hs$,
\item $\sgn\,Y(\qx,1)=\,\sgn\hs(\qx-\qx_*)\,$ and\/ $\,Y(\qx,0)>0$
\end{enumerate}
for all\/ $\,\qx\in(-\hs1,0)$, with\/ $\,\hs\sgn\hs\,$ as in
{\rm\S\ref{rati}}. Then
\begin{enumerate}
\item[a)] $Y>0\,$ on $\,(\qx_*,0)\times[\hs0,1]$.
\item[b)] For every\/ $\,\qx\in(-\hs1,\qx_*]\,$ there exists a unique\/
$\,\rx\in[\hs0,1]\,$ with\/ $\,Y(\qx,\rx)=0$.
\item[c)] The coordinate function $\,\qx$, restricted to
the set of all zeros of\/ $\,Y\hs$ in the square\/ $\,\mathcal{S}$,
maps it homeomorphically onto\/ $\,(-\hs1,\qx_*]$. In particular, that set of
zeros is empty when $\,\qx_*=-\hs1$.
\end{enumerate}
\end{lem}
In fact, by (ii), both $\,Y(\qx,0)\,$ and $\,Y(\qx,1)\,$ are positive if
$\,\qx\in(\qx_*,0)$, so that (i) and Remark~\ref{zerpo} with
$\,\varPhi(\rx)=-\hs Y(\qx,\rx)\,$ and $\,\mathcal{I}=[\hs0,1\hs]\,$ yield (a).
Next, $\,\rx\,$ required in (b) exists for $\,\qx\in(-\hs1,\qx_*]$, since (ii)
gives $\,Y(\qx,0)>0\ge Y(\qx,1)$, and it is unique in view of Remark~\ref{zerpo}
(for $\,\varPhi\,$ as above) and (i). Finally, $\,\qx\,$ sends the set of
zeros of $\,Y\hs$ bijectively onto $\,(-\hs1,\qx_*]$, while continuity of its
inverse mapping (that is, of the function $\,\qx\mapsto\rx$) follows from an
obvious subsequence argument.

We define a rational function $\,\qx\,$ of the variables $\,\th,\ha\,$ and a
constant $\,\qx_*$ by
\begin{equation}
\qx\,=\,(\ha-\th)/(\th+\ha-2)\,,\hskip20pt
\qx_*\,=\,(2-\xj)\xj/\lj(\xj-1)^2+1\rj\,\in\,(-\hs1,0)\hs,\label{qvu}
\end{equation}
for $\,\xj\,$ as in Lemma~\ref{xxast}(b). Thus, $\,\xj\,$ and $\,\qx_*$ also
depend on an odd integer $\,m\ge3$.
\begin{rem}\label{xtoze}
We have $\,\qx_*\in(-\hs1,0)\,$ in (\ref{qvu}) since the assignment
$\,\xj\mapsto\qx_*$ is an increasing diffeomorphism
$\,(-\infty,0)\to(-\hs1,0)$. In fact, it is the composite
$\,(-\infty,0)\to(1,\infty)\to(-\hs1,0)\,$ of the decreasing diffeomorphisms
$\,\xj\mapsto\zeta=(\xj-1)^2$ and
$\,\zeta\mapsto\qx_*=(1-\zeta)/(1+\zeta)=-\hs1+2/(1+\zeta)$.
\end{rem}
\begin{thm}\label{cvseg}Let\/ $\,\dz\,$ be the set of all\/
$\,(\th,\ha)\in\rto$ with\/ $\,\th<0<\hu\le\ha\le1\,$ and\/ $\,T(\th,\ha)=0$,
where\/ $\,\hu=\th/(\th-1)\,$ and\/ $\,T\hs$ is the polynomial with\/
{\rm(\ref{umv})} for a given integer $\,m\ge2$. Also, let\/
$\,\qx,\qx_*,\xj\,$ be as in\/ {\rm(\ref{qvu})} and\/
{\rm Lemma~\ref{xxast}(b)}.
\begin{enumerate}
\item[a)] If\/ $\,m\,$ is even, $\,\dz\,$ is empty.
\item[b)] If\/ $\,m\,$ is odd, $\,\dz\,$ is a real-analytic compact curve
segment embedded in $\,\rto\hskip-1.2pt$.
\end{enumerate}
For\/ odd $\,m\hs$, there exists a unique negative real number\/ $\,\yj\,$
with\/ $\,T(\yj,1)=0$. The endpoints of the curve segment\/ $\,\dz\,$ then
are\/ $\,(\xj,\xj^*)\,$ and\/ $\,(\yj,1)$, while the restriction of\/
$\,\qx\,$ to\/ $\,\dz\,$ is a homeomorphism\/ $\,\qx:\dz\to[-\hs1,\qx_*]\,$
sending\/ $\,(\yj,1)\,$ and\/ $\,(\xj,\xj^*)\,$ onto\/ $\,-\hs1\,$ and,
respectively, $\,\qx_*$.
\end{thm}
\begin{proof}At any $\,(\th,\ha)\,$ with $\,\th<0\,$ and $\,\ha=1$,
(\ref{sdt}) reads: $\,(-1)^m\hs d\hs[T(\th,1)]/d\th>0$, provided that
$\,T(\th,1)=0$. The assumptions of Remark~\ref{zerpo} thus are satisfied by
$\,\varPhi(\th)=(-1)^m\hs T(\th,1)\,$ on the interval
$\,\mathcal{I}=(-\infty,0)$. Also, $\,\varPhi(0)>0\,$ by (\ref{tov}), while
$\,(-1)^m\hs\varPhi(\th)\to\infty\,$ as $\,\th\to-\infty\,$ (since $\,\pj\hs$
in (\ref{uut}) is a degree $\,m\,$ polynomial with leading coefficient
$\,m-1>0$, cf.\ (\ref{esg})). Remark~\ref{zerpo} now shows that equation
$\,T(\yj,1)=0\,$ has no negative real solutions $\,\yj\,$ when $\,m\,$ is
even, and has exactly one such solution when $\,m\,$ is odd.

Formula $\,\varUpsilon(\th,\ha)=(\qx,\rx)\,$ with $\,\qx\,$ as in (\ref{qvu})
and $\,\rx=\th(\ha-1)/\ha\,$ defines a diffeomorphism
$\,\varUpsilon:\mathcal{K}\to\mathcal{S}\,$ of the set $\,\mathcal{K}\,$ in
the $\,\th\ha\hs$-plane, formed by all $\,(\th,\ha)\,$ with $\,\th\le0<\ha<1\,$
and $\,\hu\le\ha$, onto the square $\,\mathcal{S}=(-\hs1,0)\times[\hs0,1]\,$ in the
$\,\qx\rx\hs$-plane. This is an easy exercise; for instance, $\,\th+\ha-2<0\,$ on
$\,\mathcal{K}$, as $\,\th\le0\,$ and $\,\ha<2$, so that $\,-1<\qx<0$, while
$\,\rx=\th/\hv$, and so $\,0\le\rx\le1\,$ as $\,\hv\le\th\le0$, that is,
$\,0\le\hu\le\ha$. (Cf.\ Lemma~\ref{uuast}.) Also, solving
$\,\ha-\th=(\th+\ha-2)\qx\,$ for $\,\ha$, we can rewrite condition
$\,\rx=\th(\ha-1)/\ha$, for any $\,(\qx,\rx)\in\mathcal{S}$, as the equation
$\,(\qx+1)(\th-\rx-1)\th+2\qx\rx=0$, which is quadratic in $\,\th\,$ and
has a positive leading coefficient, while its left-hand side is nonpositive at
both $\,\th=0\,$ and $\,\th=1$, so that its only nonpositive real root is
simple. Solving it for $\,\th\,$ and using our expression for $\,\ha\,$ in
terms of $\,\th,\qx$, we now get an explicit description of the inverse
$\,\varUpsilon^{-1}$.

The diffeomorphism $\,\varUpsilon:\mathcal{K}\to\mathcal{S}\,$ sends
$\,(-1)^m\,T:\mathcal{K}\to\bbR\,$ and the vector field
$\,\,\mathbf{w}\,=(\th-1,\hs\ha-1)\,$ on $\,\mathcal{K}\,$ onto the function
$\,Y=(-1)^m\,T\circ\varUpsilon^{-1}:\mathcal{S}\to\bbR\,$ and a vector field on
$\,\mathcal{S}$. The latter equals a positive function times the coordinate vector
field in the direction of $\,\rx$, which is clear since
$\,\,d_{\hs\mathbf{w}}\qx=0\,$ for $\,\qx\,$ treated as a function of
$\,\th,\ha$, while
$\,\ha^2\hs d_{\hs\mathbf{w}}\rx=(\ha-1)[(\th-1)\ha+\th]>0\,$ on $\,\mathcal{K}$, as
$\,\ha>0\,$ and $\,\ha-1$, $\,\th-1$, $\,\th\,$ are all negative. Thus, by (\ref{sdt}), our $\,Y\hs$ satisfies (i) in Lemma~\ref{dydro}, since the inequality
$\,\th\ha\le\th+\ha\,$ (that is, $\,\hu\le\ha$) gives $\,\th\ha-\th-\ha<0\,$
on the interior of $\,\mathcal{K}$.

Moreover, $\,\varUpsilon\,$ also maps the boundary curves of $\,\mathcal{K}$,
parameterized by $\,\th\mapsto(\th,\hu)\,$ and, respectively,
$\,\ha\mapsto(0,\ha)\,$ with $\,\th\in(-\infty,0)\,$ and $\,\ha\in(0,1)$, onto
the boundary curves for $\,\mathcal{S}$, given by $\,\qx\mapsto(\qx,1)\,$ and
$\,\qx\mapsto(\qx,0)\,$ with $\,\qx\in(0,1)$, in such a way that the curve
parameter $\,\qx\,$ is an increasing (or, respectively, decreasing) function
of $\,\th\,$ (or, $\,\ha$). Lemma~\ref{xxast}(a),(b) and (\ref{tuu}) now show that
$\,Y\hs$ satisfies condition (ii) in Lemma~\ref{dydro} as well, provided that we
set $\,\qx_*=-\hs1\,$ when $\,m\,$ is even, and define $\,\qx_*$ as in (\ref{qvu}) when $\,m\,$ is odd. Note that, for odd $\,m\hs$, the function on
$\,\mathcal{K}\,$ corresponding under $\,\varUpsilon\,$ to the coordinate function
$\,\qx\,$ on $\,\mathcal{S}\,$ is, obviously, $\,\qx:\mathcal{K}\to\bbR\,$ given by (\ref{qvu}), the value of which at $\,(\xj,\xj^*)\,$ is $\,\qx_*$.

Lemma~\ref{dydro} and Remark~\ref{submf}, combined with our initial conclusion about
the equation $\,T(\yj,1)=0$, now yield both (a) and (b), completing the
proof.
\end{proof}
\begin{example}\label{xymtr}If $\,m=3\,$ we have
$\,\xj=-\hs(\sqrt 5+1)/2\,$ and $\,\yj=-\hs3/2\,$ for $\,\xj,\yj\,$ as in
Lemma~\ref{xxast}(b) and Theorem~\ref{cvseg}. This is clear from
Lemma~\ref{xxast}(c) and (\ref{uut}) since,
using (\ref{fet}) and the explicit expressions for $\,\vs\,$ and $\,E\,$ in
Example~\ref{zswtr}, we obtain
$\,(\ps-1)^3\hs\xo(\ps)=\ps^2(\ps-2)^4(\ps^2+\ps-1)\,$ and
$\,\pj(\ps)=\ps(\ps-2)(2\ps+3)$. Now $\,\xj^2=1-\xj$, so
that $\,(\xj-1)^2+1=3(1-\xj)\,$ and hence $\,\qx_*=-\hs\sqrt5/3$.
\end{example}

\section{A further symmetry}\label{fusy}
\setcounter{equation}{0}
In view of Lemma~\ref{uuast}, for $\,\hu\nh,\hv\,$ as in (\ref{ast}) the assignment
\begin{equation}
\reg\ni(\th,\ha)\,\mapsto\,(\hv\nh,\hu)\in\,\reg\,,\quad\mathrm{where}\quad
\reg\,=\,(\bbR\smallsetminus\{1\})\times(\bbR\smallsetminus\{1\})\,,
\label{fsy}
\end{equation}
is an involution $\,\hs\reg\to\hs\reg\,$ and its fixed-point set is the
hyperbola $\,\mathcal{H}\,$ given by $\,\ha=\hu\nh$. The importance of (\ref{fsy})
for our discussion is due to Proposition~\ref{tildt} below.

Both $\,\mathcal{H}\,$ and (\ref{fsy}) appear particularly simple if one replaces
$\,\th,\ha\,$ by new affine coordinates $\,a,b\,$ with $\,\th=a+1,\hs\ha=b+1$.
The equation of $\,\mathcal{H}\,$ then becomes $\,ab=1\,$ (cf.\ (\ref{ast}.b)), while (\ref{fsy}) reads $\,(a,b)\mapsto(1/b\hs,1/a)$. In other words, let us shift the
origin from $\,(\th,\ha)=(0\hs,\nh0)\,$ to $\,(\th,\ha)=(1\hs,\nh1)$, so as to
treat the original $\,\th\ha\hs$-plane as the $\,ab\hs$-plane $\,\rto$ with
linear coordinates $\,a,b$. We may now endow the latter with the indefinite
inner product corresponding to the quadratic function $\,(a,b)\mapsto ab$.
This makes $\,\mathcal{H}\,$ a unit pseu\-do\-cir\-cle centered at the origin,
while (\ref{fsy}) acts through division of any non-null vector $\,(a,b)\,$ by
its inner square $\,ab$. Thus, (\ref{fsy}) forms the pseu\-do-Euclid\-e\-an
analogue of the Euclidean conformal inversion
$\,\hs\mathbf{x}\hs\mapsto\hs\mathbf{x}/|\mathbf{x}|^2\nh$.
\begin{lem}\label{boldf}For\/ $\,F,E,\pv\hs$ as in {\rm(\ref{fet})} and\/ {\rm(\ref{puv})}
with a fixed integer\/ $\,m\ge2$, let\/ $\,\,\ef(\ps)=(\ps-1)\dot F(\ps)$,
$\,\,\ee(\ps)=(\ps-1)\dot E(\ps)$, where\/ $\,(\,\,)\dot{\,}=\,d/d\ps$. We
then have\/ $\,\ef(\hatt)=\ef(\ps)$, $\,\,\ee(\hatt)=\ef(\ps)-\ee(\ps)$, and\/
$\,\pv(\hu\nh,\hv)=\pv(\th,\ha)\,$ for any\/
$\,\ps,\th,\ha\in\bbR\smallsetminus\{1\}$, with\/ $\,\hatt\nh,\hu\nh,\hv$
defined as in {\rm(\ref{ast})}.
\end{lem}
In fact, $\,d\hatt\hskip-1.3pt/d\ps=-\hs1/(\ps-1)^2$ by (\ref{ast}.b). Applying
$\,d/d\ps\,$ and the chain rule to (\ref{fta}), we get
$\,\dot F(\hatt)=(\ps-1)^2\dot F(\ps)\,$ and
$\,\dot E(\hatt)=(\ps-1)^2[\hs\dot F(\ps)-\dot E(\ps)\hs]$, which yields the
first two relations, and (with (\ref{puv}), (\ref{fta}) and (\ref{ast}.b)), also the
third.
\begin{prop}\label{tildt}Given an integer\/ $\,m\ge2$, let\/
$\,T(\th,\ha)\,$ be as in {\rm(\ref{umv})}, and let\/
$\,\th,\ha\in\bbR\smallsetminus\{1\}$. Then, with\/ $\,\hatt,\hu\nh,\hv\,$
defined as in {\rm(\ref{ast})}, and\/ $\,\hs\sgn\,$ as in {\rm\S\ref{rati}},
\begin{enumerate}
\item[a)] $\wt(\hu\nh,\hv)=\wt(\th,\ha)$, where\/
$\,\wt(\th,\ha)=(\th-1)^{2-m}(\ha-1)^{2-m}\hs T(\th,\ha)$.
\item[b)] The involution {\rm(\ref{fsy})} leaves invariant the function
$\,\hs\sgn\,T$.
\item[c)] The set $\,\{(\th,\ha)\in\rto\hskip3pt\big|\hskip3ptT(\th,\ha)=0
\hskip5pt\mathrm{and}\hskip5pt\th\ne1\ne\ha\}\,$ is invariant under
{\rm(\ref{fsy})}.
\end{enumerate}
\end{prop}
This is clear as $\,(\th-1)(\ha-1)\hs\rh(\hu\nh,\hv)=-\hs\rh(\th,\ha)\,$ both
for $\,\rh(\th,\ha)=\ha-\th\,$ and for $\,\rh(\th,\ha)=\th\ha-\th-\ha$. Using
(\ref{umv}), Lemma~\ref{boldf} and (\ref{ast}.b), we now obtain
$\,(\th-1)^{2m-4}(\ha-1)^{2m-4}\,T(\hu\nh,\hv)=T(\th,\ha)$, and (\ref{ast})
gives (a), while (a) with $\,T(\th,\ha)=T(\ha,\th)\,$ yields (b), and (b)
implies (c).
\begin{rem}\label{cntrs}
The involution (\ref{fsy}) admits an interesting
algebraic-geometric interpretation in terms of the projective plane
in which the $\,ab\hs$-plane corresponding (as described above) to the
original $\,\th\ha\hs$-plane is canonically embedded. Namely, the
homogeneous-coordinate form of (\ref{fsy}) is
$\,[\hs a:b:c\hs]\mapsto[\hs ab:ac:bc\hs]$, so that (\ref{fsy}) is a {\it
quadratic transform with the centers} $\,\,[\hs0:0:1\hs]\hs$,
$\,[\hs0:1:0\hs]\,$ and $\,[\hs1:0:0\hs]$. (See \cite[pp.\
496--498]{griffiths-harris}.) In other words, (\ref{fsy}) consists of a
blow-up at the three centers, followed by a blow-down of the three projective
lines through each pair of centers. The three centers lie in the projective
closure $\,\overline{\mathcal{T}}\,$ of the
curve $\,\mathcal{T}\,$ with the equation $\,T\nh=0$, and one can show that
$\,m-2\,$ is their common algebraic multiplicity. Hence they are singularities
of $\,\overline{\mathcal{T}}\,$ for odd $\,m>3$, while $\,\overline{\mathcal{T}}\,$ also
has a fourth singularity, of multiplicity $2(m-2)$, at the point
$\,[\hs a:b:c\hs]\,$ with $\,a=b=-\hs1$, $\,c=1$, that is, at the origin of
the $\,\th\ha\hs$-plane.

Proposition~\ref{tildt} thus states that $\,\overline{\mathcal{T}}\,$ is invariant
under a quadratic transform with centers which, for odd $\,m>3$, are its
low-order singularities.
\end{rem}
\begin{rem}\label{invol}
A well-known involution (cf.\ \cite{apostolov-calderbank-gauduchon}) assigns to a
quadruple $\,\mgmt$ with (\ref{zon}) or (\ref{zto}), such that
$\,g\,$ is not locally reducible as a K\"ahler metric, its {\it dual\/}
$\,\hatmgmt$. Here $\,\hat M,M\,$ coincide as real manifolds, but have
different complex structures, so that $\,\hat m=m\hs$, while
$\,\hat g=\y^{\hs2}g/(\vp-\y)^2$ and $\,\hat\vp=\y\hs\vp/(\vp-\y)$, for
$\,\y\,$ determined by $\,\mgmt\,$ as in Remark~\ref{detrm}. (Cf.\ Theorem~\ref{cnver}
and \S\ref{redu}.) This duality involution induces the mapping (\ref{fsy}) on
$\,\mathcal{C}\smallsetminus\{(0\hs,\nh0)\}\,$ when $\,\mgmt\,$ is replaced by the
point $\,(\th,\ha)\in\mathcal{C}\,$ associated with it as in \S\ref{intr} or Remark~\ref{detrm}.
\end{rem}

\section{Equation $\,T(\th,\ha)=0\,$ with $\,\th<0<\ha<1$}\label{tzer}
\setcounter{equation}{0}
\begin{prop}\label{crvsg}Given an integer\/ $\,m\ge2$, let\/ $\,\vg\hs$
be the set of all\/ $\,(\th,\ha)\in\rto$ with\/ $\,\th<0<\ha<1\,$ and\/
$\,T(\th,\ha)=0$, for\/ $\,T$ as in {\rm Lemma~\ref{polyt}}.
\begin{enumerate}
\item If\/ $\,m\,$ is even, $\,\vg\,$ is empty.
\item If\/ $\,m\,$ is odd, $\,\vg\,$ is a real-analytic submanifold of\/
$\,\rto\nh$, diffeomorphic to\/ $\,\bbR\hs$, containing the point\/
$\,(\xj,\xj^*)\,$ described in {\rm Lemma~\ref{xxast}(b)}, and invariant
under the involution {\rm(\ref{fsy})}. Furthermore, {\rm(\ref{fsy})} keeps\/
$\,(\xj,\xj^*)\,$ fixed and interchanges the two connected components of\/
$\,\vg\smallsetminus\{(\xj,\xj^*)\}$. One of these components is unbounded,
while the closure of the other is the compact curve segment\/ $\,\dz\hs$
appearing in {\rm Theorem~\ref{cvseg}(b)}. For $\,\qx\,$ and\/ $\,\qx_*$ as
in\/ {\rm(\ref{qvu})}, the function $\,\qx\,$ sends either component
homeomorphically onto $\,(-\hs1,\qx_*)$, while its value at\/
$\,(\xj,\xj^*)\,$ is $\,\qx_*$.
\end{enumerate}
\end{prop}
\begin{proof}Let $\,\vk=\dz\cap(\bbR\times\{1\})$, with $\,\dz\,$ as in Theorem~\ref{cvseg}. Thus, $\,\vg\,$ is the union of the set $\,\dz\smallsetminus\vk\,$
and its image under the involution (\ref{fsy}). (In fact, Proposition~\ref{tildt}(c)
and Lemma~\ref{uuast} show that (\ref{fsy}) sends the set
$\,(-\infty,0)\times(0,1)\,$ onto itself by interchanging its subsets lying
``above'' and ``below'' the hyperbola $\,\ha=\hu$ and keeping each point of
the hyperbola fixed.) Now (i) follows from Theorem~\ref{cvseg}(a).

If $\,m\,$ is odd, Theorem~\ref{cvseg} gives $\,\vk=\{(\yj,1)\}$, while both
$\,\dz\smallsetminus\{(\yj,1)\}\,$ and its image under (\ref{fsy}) contain
$\,(\xj,\xj^*)\,$ (Theorem~\ref{cvseg}(b)), so that their union $\,\vg\,$ is
connected. Thus, according to Remark~\ref{submf}, $\,\vg\,$ is contained in
$\,\rto$ as a (connected) one\diml\/ real-analytic submanifold without
boundary. Finally, by (\ref{ast}.a), $\,\hv\to-\infty\,$ and $\,\hu\to\yj^*$ as
$\,(\th,\ha)\in\dz\smallsetminus\{(\yj,1)\}\,$ approaches $\,(\yj,1)$, which
shows that $\,\vg\,$ is unbounded. This completes the proof.
\end{proof}
\begin{lem}\label{mumil}Given an integer\/ $\,m\ge3$, let\/
$\,T,\lx,\ax,\bx\,$ be the rational functions of\/ $\,\th,\ha\,$ defined in
{\rm Lemma~\ref{polyt}}, {\rm(\ref{lba})} and\/ {\rm(\ref{abc})} with\/
$\,E,F,\vs\,$ as in {\rm(\ref{fet})} -- {\rm(\ref{esg})}, and let\/ $\,\mx\,$
be the rational function of\/ $\,\th,\ha\,$ with\/
$\,\mx\,=\,\th\ha/(\th+\ha-2)$. Then
\begin{equation}
\mx\,-\,\lx\,=\,{(\ha-\th)^3\th(\th-2)\ha(\ha-2)\hs T(\th,\ha)
\over\,m\hs(\th+\ha-2)(\th-1)^m(\ha-1)^m\ax\,}\label{mml}
\end{equation}
in the sense of equality between rational functions. Also, $\,0<\mx<1\,$ at
any $\,(\th,\ha)\in(-\infty,0)\times(0,1\hs]$. Finally, when $\,m\,$
is odd and\/ $\,(\th,\ha)\in(-\infty,0)\times(0,1)$,
\begin{enumerate}
\item $\ax<0\,$ at\/ $\,(\th,\ha)$, with\/ $\,\ax\,$ as
in {\rm(\ref{abc})}.
\item $T(\th,\ha)=0\,$ if and only if\/ $\,\lx=\mx\,$
at\/ $\,(\th,\ha)$.
\end{enumerate}
\end{lem}
In fact, $\,0<\mx<1\,$ on $\,(-\infty,0)\times(0,1\hs]\,$ as
$\,(\th-1)(\ha-1)\ge0>-\hs1\,$ there, and hence also $\,\th+\ha-2<\th\ha<0$.
Next, dividing (\ref{vuu}.i) by $\,m\av\fj$, with $\,\av,\fj\,$ as in (\ref{aph}),
and using (\ref{auv}), (\ref{vuu}.ii) and (\ref{lba}), we obtain (\ref{mml}).
Finally, if $\,m\,$ is odd,
$\,F(\th)<0<F(\ha)\,$ as $\,\th<0<\ha<1\,$ (see (\ref{mon})). Thus,
$\,e_\th-e_\ha>0\,$ for $\,e_\th=E(\th)/F(\th)$, $\,e_\ha=E(\ha)/F(\ha)\,$ (as
$\,e_\ha<0<1<e_\th\,$ by (\ref{mon})), and so
$\,\ax=(e_\th-e_\ha)\hs F(\th)F(\ha)<0$, which proves (i). Now (ii) follows:
the denominators involved (including those in (\ref{lba})), and the factor
$\,(\ha-\th)\th(\th-2)\ha(\ha-2)$, are all nonzero when $\,\th<0<\ha<1$.
(Cf.\ (i) and the obvious relation $\,\th+\ha-2<0$.)
\begin{prop}\label{tzrok}Given an odd integer\/ $\,m\ge3\,$ and real
numbers\/ $\,\th,\ha\,$ with\/ $\,\th<0<\ha<1$, let\/ $\,T\hs$ be the
polynomial defined by {\rm(\ref{umv})}, and let\/ $\,Q\in\jv\smallsetminus\{0\}\,$
be a function with {\rm(\ref{qdq}.b)} on the interval\/ $\,I=[\th,\ha]$, where\/
$\,\jv\hs$ is the space {\rm(\ref{spv})}. In view of\/ {\rm Lemma~\ref{onedi}}, such\/
$\,Q\,$ exists and is unique up to a nonzero constant factor.

If\/ $\,T(\th,\ha)=0$, then $\,Q\,$ and $\,I$ satisfy all five conditions
{\rm(\ref{qdq})}.
\end{prop}
In fact, $\,\lx=\mx>0\,$ at $\,(\th,\ha)\,$ in view of Lemma~\ref{mumil}(ii) and
relation $\,\mx>0\,$ in Lemma~\ref{mumil}. Thus, by Lemma~\ref{mumil}(i) and
(\ref{lba}), $\,\lx\ax<0\,$ and $\,\bx\hs\vs(0)>\ax\,$ at $\,(\th,\ha)$, with
$\,\ax,\bx,\vs\,$ given by (\ref{abc}), (\ref{esg}). However, (\ref{abc}) also
yields $\,\bx<0$, as $\,F\,$ is strictly increasing on $\,(-\infty,1)\,$ (see
\S\ref{mono}). Hence $\,\ax/\bx>\vs(0)=-\hs E(0)$, cf.\ (\ref{sto}.ii). Proposition~\ref{dtefu} now implies that $\,Q\,$ and $\,I\hs$ satisfy (\ref{qdq}.a), (\ref{qdq}.c) and (\ref{qdq}.d), while relation $\,T(\th,\ha)=0\,$ yields (\ref{qdq}.e) (see Theorem~\ref{facto}).

\section{Expansions of $\,T\,$ about $\,(0\hs,\nh0)\hs$ and
$\,(1\hs,\nh1)$}\label{expa}
\setcounter{equation}{0}
For $\,E,\cj\,$ as in (\ref{fet}) and (\ref{vuu}.ii) with an integer $\,m\ge2$,
we have $\,\cj={2m-1\choose m-1}\,$ (see (\ref{sto}.i,\hs iii)), and there
exists a polynomial $\,D\,$ such that, for all $\,\th\in\bbR\hs$,
\begin{equation}
\begin{array}{rl}
\mathrm{a)}&
(\th-1)^m\hs E(\th)\,=\,(\th-2)\th^{2m-1}\,+\,\hs D(\th)\hs,\\
\mathrm{b)}&
D(\th)\,-\,\cj\hs(\th-1)^m\,=\,(\th-2)\,\sum_{j=0}^{m-1}(-1)^{m+j}
{2m-1\choose j}{2m-j-2\choose m-1}\th^j\nh.\end{array}
\label{ume}
\end{equation}
Namely, $\,D_1(\th)=1\,$ and, by (\ref{fmt}.ii),
$\,D_m(\th)=\th^2\hs D_{m-1}(\th)-m^{-1}{2m-2\choose m-1}(\th-1)^m$ whenever
$\,m\ge2$, for $\,D=D_m$ {\it defined\/} by (\ref{ume}.a). Now
$\,D(\th)=m^{-1}\sum_{j=0}^{m-1}(-1)^{m+j+1}{2m\choose j}
{2m-j-2\choose m-1}\th^j$, which follows from easy induction on $\,m\hs$, and
amounts to (\ref{ume}.b).
\begin{lem}\label{sympo}Every symmetric polynomial in the variables\/
$\,\th,\ha\,$ can be uniquely written as a combination of the products\/
$\,(\th\ha)^j\,\varTheta_k(\th,\ha)$, where\/ $\,j\ge0\,$ and\/ $\,k\ge2\,$
are integers, and\/
$\,\varTheta_k(\th,\ha)=\sum_{j=1}^{k-1}j\hs(k-j)\hs \th^{j-1}\ha^{k-j-1}$.
\end{lem}
In fact, the space of degree $\,m\,$ homogeneous symmetric polynomials in
$\,\th,\ha\,$ has the obvious basis
$\,\varPhi_\ro=(\th\ha)^\ro(\th^{m-2\ro}+\ha^{m-2\ro})\,$ with
$\,\ro\in\bbZ\,$ and $\,0\le\ro\le m/2$. Then
$\,\tilde\varPhi_\ro=(\th\ha)^\ro\,\varTheta_{m-2\ro+2}(\th,\ha)\,$ form
another basis of that space. Namely, $\,\tilde\varPhi_\ro$ equals
$\,(m-\ro+1)\varPhi_\ro$ plus a combination of $\,\varPhi_\sj$ with
$\,\ro<\sj\le m/2$, as one sees pairing up, for each $\,j$, the $\,j$th and
$\,(k-j)$th terms in the formula for $\,\varTheta_k(\th,\ha)$. Thus, the
triangular matrix expressing the $\,\tilde\varPhi_\ro$ through the
$\,\varPhi_\ro$ is invertible.
\begin{rem}\label{cubth}
By Lemma~\ref{sympo}, any symmetric polynomial
$\,\varPhi\,$ in the variables $\,a,b\,$ has an expansion
$\,\sum_{j,k}\mathbf{c}_{j,k}\hs(ab)^j\hs\varTheta_k(a,b)\,$ with some
unique coefficients $\,\hs\mathbf{c}_{j,k}$, indexed by integers $\,j,k$,
and such that $\,\hs\mathbf{c}_{j,k}\ne0\,$ for at most finitely many pairs
$\,j,k$, all of which have $\,j\ge0\,$ and $\,k\ge2$. Expanding
$\,(a-b)^3\hs\varTheta_k(a,b)\,$ into powers of $\,a\,$ and $\,b$, we get
\begin{equation}
(a-b)^3\hs\varTheta_k(a,b)\,
=\,(k-1)a^{k+1}-\hs(k+1)a^kb\hs+(k+1)ab^k-\hs(k-1)b^{k+1}\nh,
\label{abt}
\end{equation}
and so, for any $\,\ro,\sj\in\bbZ$, the coefficient of $\,a^\ro b^\sj$ in the
monomial expansion of $\,(b-a)^3\hs\varPhi\,$ is
\begin{equation}
\begin{array}{lllll}\arraycolsep9pt
(\sj-\ro\hs-2)\hs\mathbf{c}_{\ro,\hs\sj-\ro\hs-1}
&+&(\ro\hs-\sj-2)\hs\mathbf{c}_{\ro\hs-1,\hs\sj-\ro\hs+1}&&\\
&+&(\sj-\ro\hs+2)\hs\mathbf{c}_{\sj,\hs\ro\hs-\sj-1}
&\hskip-3pt+&(\ro\hs-\sj+2)\hs\mathbf{c}_{\sj-1,\hs\ro\hs-\sj+1}\,.
\end{array}
\label{src}
\end{equation}
\end{rem}
\begin{lem}\label{eaplu}Let\/ $\,\lf(\th,\ha)=E(\th)\hs\fy-\cj\hs\fj\,$
with\/ $\,E,\cj\,$ as in {\rm(\ref{fet})}, {\rm(\ref{vuu}.ii)} and\/
$\,\fj,\fy\,$ depending on $\,\th,\ha\in\bbR\,$ as in {\rm(\ref{aph})} for a
fixed integer\/ $\,m\ge2$. Then, for\/ $\,a,b\in\bbR\hs$,
\begin{enumerate}
\item $E(a+1)\,=\,2\hs\sum_{j=1}^mj(m+j)^{-1}\gm(j)\hs a^j$, \ and
\item $\lf(a+1,b+1)=(a^2-1)\hs\sum_{j=1}^m(jb-j+1)\hs\gm(j)\hs a^{j-1}$,
\hskip6ptwhere\/ $\,\gm(j)={2m-1\choose m-j}$.
\end{enumerate}
\end{lem}
\begin{proof}We prove (i) by induction on $\,m\ge1$, writing, as in
(\ref{fmt}), $\,E_m,\gm_m$ for $\,E,\gm$. If $\,m=1$, (i) is trivial as
$\,E_1(a+1)=a$. The inductive step: (\ref{fmt}.ii) multiplied by $\,2m-1\,$
gives  $\,(2m-1)E_m(a+1)=(2m-1)\hs a^{-1}(a+1)^2\hs E_{m-1}(a+1)-\gm_m(1)$.
Assuming (i) for $\,m-1\,$ rather than the given $\,m\ge2$, let us first
substitute for $\,E_{m-1}(a+1)$, in the last equality, the sum as in (i)
involving the $\,\gm_{m-1}(j)\,$ instead of $\,\gm_m(j)$, then multiply by
$\,(a+1)^2$ and make the summation index $\,j\,$ coincide with the exponent in
$\,a^j$, next replace$\,(2m-1)(m+j-1)^{-1}\hs\gm_{m-1}(j)\,$ by
$\,(2m-2)^{-1}(m-j)\hs\gm_m(j)$, as well as $\,(m-j+1)\hs\gm_m(j-1)\,$ by
$\,(m+j-1)\hs\gm_m(j)\,$ and $\,\gm_m(j+1)\,$ by
$\,(m+j)^{-1}(m-j)\hs\gm_m(j)$. This yields (i) for $\,m\hs$, as required.
Assertion (ii) will in turn follow once we establish the identities
\begin{enumerate}
\item[iii)] $(m+a-ma)\hs E(a+1)\,-\,a\hs\gm(1)\,
=\,(a^2-1)\sum_{j=1}^m(1-j)\hs\gm(j)\hs a^{j-1}$,
\item[iv)] $(ma-m+1)\hs E(a+1)\,-\,\gm(1)\,
=\,(a^2-1)\sum_{j=1}^mj\hs\gm(j)\hs a^{j-1}$.
\end{enumerate}
In fact, to get (ii) one can add (iv) multiplied by $\,b\,$ to (iii), as (\ref{aph}.c) for $\,(\th,\ha)=(a+1,b+1)\,$ gives
$\,\fy=(ma-m+1)\hs b+(m+a-ma)\,$ and $\,\fj=a+b$.

To prove (iii) and (iv), we use (i) to rewrite all four expressions as
combinations of the powers $\,a^j$ and verify that the corresponding
coefficients agree. Such a coefficient, found by expressing $\,\gm(j\pm1)\,$
through $\,\gm(j)\,$ (as above, with the subscript $\,m$), turns out to be
$\,0,\hs\gm(2),\hs1-m\,$ for $\,j=0,\hs1,\hs m+1\,$ in (iii) and
$\,-\hs\gm(1),\hs m\,$ for $\,j=0,m+1\,$ in (iv), while for (iii) and
$\,j=2,\dots,m\,$ (or, respectively, (iv) and $\,j=1,\dots,m$) it equals
$\,2\hs\gm(j)/[(m+j)(m-j+1)]\,$ times $\,3mj-2mj^2+m^2+j^2-m-j\,$ (or,
respectively, times $\,2mj^2-mj-m^2-j^2+j$). This completes the proof.
\end{proof}
\begin{lem}\label{xpant}Let\/ $\,m\ge2\,$ be an integer. For the polynomials\/
$\,T,\,\varTheta_k$ defined as in {\rm Lemmas~\ref{polyt}} and\/
{\rm\ref{sympo}}, and\/ $\,\th,\ha,a,b\in\bbR\hs$, we then have\/
$\,\hs\mathrm{deg}\hskip3ptT=3(m-2)\,$ and
\begin{enumerate}
\item[a)] $\,T(\th,\ha)\,=\,{\displaystyle{1\over m-1}\,
\sum_{\hs k=m}^{\hs2m-2}(-1)^{m+k}{2m-1\choose k+1}{k-2\choose m-2}
(\th\ha)^{2m-k-2}\hs\varTheta_k(\th,\ha)}$.
\item[b)] $\,T(a+1,b+1)\,=\,\,{\displaystyle\sum_{\hs j=0}^{\hs m-2}
\sum_{k=\hs m-j}^{2m-j-2}}
\mathbf{c}_{jk}\hs(ab)^j\hs\varTheta_k(a,b)$, with $\,\hs\mathbf{c}_{jk}$
given by
\item[c)] $\,\mathbf{c}_{jk}\,
=\,{\displaystyle{(2m-1)(m-j-1)(k-m+j+1)\over(m-1)(k^2-1)}\,{2m-2\choose j}
{2m-2\choose k+j}}$.
\end{enumerate}
\end{lem}
\begin{proof}Formula (\ref{vuu}.i) with $\,(\th,\ha)=(a+\nh1,b+\nh1)\,$ gives
$\,(b-a)^3(a^2-1)(b^2-1)\hs T(a+1,b+1)
=(b^2-1)(b+1)^{2m-2}\hs a^m\hs\lf(a+1,b+1)
-(a^2-1)(a+1)^{2m-2}\hs b^m\hs\lf(b+1,a+1)$, for $\,\lf\,$ as in Lemma~\ref{eaplu}. By Lemma~\ref{eaplu}(ii), this polynomial equality yields
\begin{equation}
\begin{array}{rl}\arraycolsep9pt
\mathrm{i)}&
(b-a)^3\hs T(a+1,b+1)\,=\,\varPsi(a,b)\,-\,\varPsi(b,a)\hs,\hskip12pt
\mathrm{as\ well\ as}\\
\mathrm{ii)}&
\varPsi(a,b)=\sum_{\ro\hs=\hs m}^{2m-1}\sum_{\sj=\hs0}^{2m-1}
\hs\mathbf{d}_{\ro\sj}\hs a^\ro b^\sj\hs,\hskip7pt\mathrm{where}\hskip6pt
\varPsi\hskip6pt\mathrm{is\ the\ polynomial\ with}\\
\mathrm{iii)}&
(a^2-1)\hs\varPsi(a,b)\,=\,(b+1)^{2m-2}\hs a^m\hs\lf(a+1,b+1)\hs,\hskip12pt
\mathrm{and}\\
\mathrm{iv)}&
\mathbf{d}_{\ro\sj}\,=\,\hs{2m-1\choose\ro}
{2m-1\choose\sj}\hs[\hs m-(\ro+\sj)+2\ro\sj/(2m-1)]\hs.\end{array}
\label{bat}
\end{equation}
According to Lemma~\ref{sympo},
$\,T(a+1,b+1)=\,\sum_{j,k}\mathbf{c}_{j,k}\hs(ab)^j\hs\varTheta_k(a,b)\,$
with $\,\hs\mathbf{c}_{j,k}$ as in Remark~\ref{cubth}. In terms of the conditions
\begin{equation}
\mathrm{i)}\quad0\le\sj\le m-1<\ro<2m\hs,\hskip24pt
\mathrm{ii)}\quad0\le\ro\le m-1<\sj<2m\label{osm}
\end{equation}
for $\,\ro,\sj\in\bbZ$, these $\,\hs\mathbf{c}_{j,k}$ satisfy the following
system of linear equations:
\begin{equation}
\begin{array}{l}
\mathrm{expression\ (\ref{src})\ equals\hskip-3.2pt:}\\
\mathbf{d}_{\ro\sj}\mathrm{\ in\ case\
(\ref{osm}.i),\ }\,\,-\hs\mathbf{d}_{\ro\sj}\mathrm{\ in\ case\
(\ref{osm}.ii),\ and}\hskip5pt0\hskip5pt\mathrm{otherwise,}
\end{array}
\label{srt}
\end{equation}
with $\,\ro,\sj\in\bbZ\,$ and $\,\hs\mathbf{d}_{\ro\sj}$ as in (\ref{bat}.iv);
see (\ref{bat}.i) and Remark~\ref{cubth}. Moreover, the
$\,\hs\mathbf{c}_{j,k}$
form the {\it unique\/} solution to (\ref{srt}), under the finiteness
requirement of Remark~\ref{cubth}. (The coefficients of $\,a^\ro b^\sj$ in the
expansion of the right-hand side (\ref{bat}.i) are zero for those $\,\ro,\sj\,$
which both lie in the range $\,m,\dots,2m-1$, due to the subtraction in (\ref{bat}.i) and the obvious symmetry relation
$\,\hs\mathbf{d}_{\ro\sj}=\,\mathbf{d}_{\sj\hskip-1pt\ro}$.)

We now show that by setting $\,\hs\mathbf{c}_{j,k}=\,\mathbf{c}_{jk}$ for
$\,\hs\mathbf{c}_{jk}$ as in (c) when $\,j\,$ and $\,j+k-m\,$ lie in
$\,\{0,\dots,m-2\}$, and $\,\hs\mathbf{c}_{j,k}=0\,$ otherwise, one obtains
a solution to (\ref{srt}).

First, let us assume (\ref{osm}.i). If $\,\sj=0$, all but the third term on
the left-hand side then must vanish, as they involve $\,\hs\mathbf{c}_{j,k}$
with $\,k<0\,$ or $\,j<0$, and we get (\ref{srt}) by considering the separate
cases $\,m=\ro=2\,$ and $\,\ro>2$. If $\,\sj=m-1$, only the fourth term may be
nonzero (as the others involve $\,\hs\mathbf{c}_{j,k}$ with $\,k<0\,$ or
$\,k=0\,$ or $\,j=m-1$), and (\ref{srt}) is easily verified. If, however,
$\,0<\sj<m-1<\ro<2m\hs$, the first two terms vanish as they have $\,j>m-2$, and
two cases are possible: $\,\ro=m=\sj+2\,$ (so that the third term, with
$\,k=1$, is zero, and (\ref{srt}) easily follows), or $\,\ro-\sj>2$, and a
simple calculation again gives (\ref{srt}).

Now let (\ref{osm}.ii) be satisfied. Interchanging $\,\ro\,$ and $\,\sj\,$ we
reduce this case to (\ref{osm}.i), since both sides of (\ref{srt}) are
antisymmetric in $\,\ro,\sj\,$ (as the second line arises from the first by
switching $\,\ro,\sj\,$ and changing the sign, while
$\,\hs\mathbf{d}_{\ro\sj}=\,\mathbf{d}_{\sj\hskip-.3pt\ro}$).

However, if we have neither (\ref{osm}.i) nor (\ref{osm}.ii), the right-hand side
vanishes, and so do all four terms on the left-hand side, as our bounds on
$\,j\,$ and $\,j+k-m\,$ in the definition of $\,\hs\mathbf{c}_{j,k}$ show
that condition (\ref{osm}.i) (or, (\ref{osm}.ii)) is necessary for the second (or,
first) line to contain nonzero term. Uniqueness of the solution to (\ref{srt})
now proves (b).

Next, with $\,\cj,\fj,\fy\,$ and $\,\lf(\th,\ha)=E(\th)\hs\fy-\cj\hs\fj\,$ as
in Lemma~\ref{eaplu}, (\ref{bat}.iii) for $\,(a,b)=(\th-1,\ha-1)\,$ and
(\ref{ume}.a) show that $\,(\th-2)\th\hs\varPsi(\th-1,\ha-1)\,$ equals
$\,\ha^{2m-2}\hs\lj\hs(\th-2)\th^{2m-1}\hs+\hs D(\th)
-\cj\hs(\th-1)^m\hs\rj\hs\fy\,+\,\cj\hs\ha^{2m-2}(\th-1)^m(\fy-\fj)$.
Using (\ref{ume}.b) and the relation $\,\fy-\fj=m\hs(\th-2)(\ha-2)$, we can
now divide both sides by $\,\th-2$, obtaining a polynomial expression for
$\,\th\hs\varPsi(\th-1,\ha-1)\,$ that involves two expansions into powers of
$\,\th$, coming from (\ref{ume}.b) and the binomial formula for
$\,(\th-1)^m\nh$. As $\,m\hs\cj=(2m-1){2m-2\choose m-1}\,$ by (\ref{vuu}.ii)
and (\ref{sto}.i), the overall coefficient of the $\,0$th power of $\,\th\,$
is zero and so we can divide both sides by $\,\th$, getting
\begin{equation}
\begin{array}{rcl}\arraycolsep9pt
\varPsi(\th-1,\ha-1)&=&(\th\ha)^{2m-2}\hs[m\hs\th\ha-(2m-1)(\th+\ha-2)]\\[2pt]
&+&\,\hs\sum_{j=0}^{m-2}
(-1)^{m+j}{2m-1\choose j}{2m-j-3\choose m-1}\th^j\ha^{2m-1}\\[4pt]
&+&{2m-1\over m-1}\,\sum_{j=1}^{m-1}(-1)^{m+j}{2m-2\choose j-1}
{2m-j-3\choose m-2}\th^j\ha^{2m-2}\hs.
\end{array}
\label{psu}
\end{equation}
Assertion (a) claims that $\,T=\varPhi\,$ for a polynomial $\,\varPhi\,$ given
by $\,\varPhi(\th,\ha)=\,\sum_{j,k}\mathbf{c}_{j,k}\hs(\th\ha)^j\hs
\varTheta_k(\th,\ha)\,$ with a specific {\it new meaning} of the
coefficients $\,\hs\mathbf{c}_{j,k}$. Clearly, (a) will follow if we prove
that the equality
$\,(\ha-\th)^3\hs T(\th,\ha)=\varPsi(\th-1,\ha-1)-\varPsi(\ha-1,\th-1)$,
immediate from (\ref{bat}.i), still holds when $\,T\,$ is replaced by
$\,\varPhi$. This, according to Remark~\ref{cubth}, amounts to showing that
(\ref{src}), for any $\,\ro,\sj\in\bbZ$, equals the coefficient of
$\,\th^\ro\ha^\sj$ in the expansion of
$\,\varPsi(\th-1,\ha-1)-\varPsi(\ha-1,\th-1)$. The latter coefficient is
easily obtained from the right-hand side of (\ref{psu}), with no contribution
from the first line (which is symmetric in $\,\th,\ha$). Also,
$\,\hs\mathbf{c}_{j,k}$ are nonzero only for integers $\,j,k\,$ with
$\,j+k=2m-2\,$ and $\,0\le j\le m-2$, and so, for any given
$\,\ro,\sj\in\bbZ$, at most one of the four expressions
$\,\hs\mathbf{c}_{j,k}$ occurring in (\ref{src}) may be nonzero. Specifically,
the first (or second, third, fourth) $\,\hs\mathbf{c}_{j,k}$ in (\ref{src}) is
nonzero only if $\,0\le\ro\le m-2\,$ and $\,\sj=2m-1$, or
$\,1\le\ro\le m-1\,$ and $\,\sj=2m-2$, or
$\,0\le\sj\le m-2\,$ and $\,\ro=2m-1$, or, respectively,
$\,1\le\sj\le m-1\,$ and $\,\ro=2m-2$. The required equality is easily
verified in each of these four cases, while both sides are zero when none of
the four cases occurs. We have thus established (a). Finally,
$\,\hs\mathrm{deg}\hskip3ptT=3(m-2)\,$ as the $\,k$th term in (a) is
homogeneous of degree $\,4m-k-6$. This completes the proof.
\end{proof}
\begin{rem}\label{graph}
Let $\,m=3$. Lemma~\ref{xpant}(a) then yields
$\,T(\th,\ha)=5(\th\ha^2+\ha\th^2)-(3\th^2+3\ha^2+4\th\ha)
=2\hs[5(\xi-1)\xi^2-(5\xi+1)\eta^2]$, with $\,\xh,\yh\,$ given by (\ref{uxe}).
Fig.~1 thus conveys the correct idea of what the set given by
$\,T(\th,\ha)=0\,$ looks like when $\,m=3$. Namely, the lines $\,\ha=\th\,$
and $\,\ha=-\hs\th\,$ are the $\,\xh\,$ and $\,\yh\,$ coordinate axes, while
$\,T(\th,\ha)=0\,$ amounts to $\,\xi\ne-\hs1/5\,$ and
$\,\yh^2=(\xh-1)\xh^2/(\xh+1/5))$. This gives $\,\xh\ge1\,$ or
$\,\xh<-\hs1/5$, and $\,\yh=\pm\hs\xh\hs[\hs1-6/(5\xh+1)\hs]^{1/2}\nh$.
\end{rem}

\section{The sign of $\,T\,$ on specific regions}\label{sgtr}
\setcounter{equation}{0}
Assertion (i) below is needed only to derive (ii) and (iii); however, (ii),
(iii) also follow from our formula (\ref{abt}) combined with Lemma~32.1 of
\cite{potentials} for $\,\beta=\th/\ha$.
\begin{lem}\label{sgnth}For an integer\/ $\,k\ge2$, any\/
$\,\th,\ha\in\bbR\hs$, and\/ $\,\varTheta_k\,$ as in {\rm Lemma~\ref{sympo}},
\begin{enumerate}
\item $\varTheta_k(\th,\ha)\,
=\,\sum_{\hs1\,\le\,j\,\le\,k/2}\,j{k+1\choose2j+1}\xh^{k-2j}\yh^{2j-2}$ for
$\,\xh,\yh\,$ given by\/ {\rm(\ref{uxe})}.
\item If\/ $\,k\,$ is even, $\,\varTheta_k(\th,\ha)>0\,$ unless
$\,\th=\ha=0$, while\/ $\,\varTheta_k(0\hs,\nh0)=0$.
\item If\/ $\,k\,$ is odd,
$\,\hs\sgn\hs(\varTheta_k(\th,\ha))=\,\sgn\hs(\th+\ha)$, with\/
$\,\hs\sgn\hs\,$ as in\/ {\rm\S\ref{rati}}.
\end{enumerate}
\end{lem}
Here (i) is straightforward if one multiplies both sides by
$\,(\th-\ha)^3=8\yh^3$ and verifies that the resulting polynomials in
$\,\xh,\yh\,$ coincide by applying (\ref{uxe}) and the binomial formula to
rewrite the right-hand side of (\ref{abt}) as a function of $\,\xh,\yh$.
Now (ii) and (iii) follow since, when $\,k\,$ is even (or, odd), (i) expresses
$\,\varTheta_k(\th,\ha)\,$ as a sum of squares (or, respectively, as $\,\xh\,$
times a sum of squares).

Let $\,T\hs$ be the polynomial defined by (\ref{umv}) with any given integer
$\,m\ge2$. Then $\,T\nh=1\,$ for $\,m=2\,$ (cf.\ Lemma~\ref{xpant}(a)), while, if
$\,m\ge3$,
\begin{equation}
\begin{array}{rll}\arraycolsep9pt
\mathrm{a)}&
\sgn\,T\,=\,(-1)^m\hskip14pt&\mathrm{at\ any}\hskip5pt(\th,\ha)\ne(0\hs,\nh0)
\hskip5pt\mathrm{with}\hskip5pt(\th+\ha)\hs\th\ha\,\le\,0\hs,\\
\mathrm{b)}&
\sgn\,T\,=\,(-1)^m&\mathrm{on}\hskip5pt[\hs0,1)\times[\hs0,1)
\hskip5pt\mathrm{except\ at}\hskip5pt(0\hs,\nh0)\hskip5pt\mathrm{or}
\hskip5pt(1\hs,\nh1)\hs,\\
\mathrm{c)}&
\sgn\,T\,=\,(-1)^m&\mathrm{on}\hskip5pt(-\infty,0\hs]\times
(-\infty,0\hs]\hskip5pt\mathrm{except\ at}\hskip5pt(0\hs,\nh0)\hs,\\
\mathrm{d)}&
\sgn\,T\,=\,1&\mathrm{on}\hskip5pt[\hs1,\infty)\times
[\hs1,\infty)\hskip5pt\mathrm{except\ at}\hskip5pt(1\hs,\nh1)\hs.
\end{array}
\label{sgt}
\end{equation}
In fact, let $\,\th+\ha\ge0\,$ and $\,\th\ha\le0$, or $\,\th+\ha\le0\,$ and
$\,\th\ha\ge0$. The even $\,k\,$ (or, odd $\,k$) summands in
Lemma~\ref{xpant}(a) are nonnegative by Lemma~\ref{sgnth}(ii) (or,
respectively, Lemma~\ref{sgnth}(iii)). This yields (\ref{sgt}.a), and hence
(\ref{sgt}.c). Next, as $\,\hu\nh,\hv\in(-\infty,0\hs]\,$ whenever
$\,\th,\ha\in[\hs0,1)\,$ (cf.\ Lemma~\ref{uuast}), (\ref{sgt}.c) and
Proposition~\ref{tildt}(b) imply (\ref{sgt}.b). Finally, Lemmas~\ref{xpant}(b)
and~\ref{sgnth}(ii),\hs(iii) give (\ref{sgt}.d).

\section{Third subcase of \hs{\rm(\ref{qdq}.e)}: condition
\hs{\rm(\ref{trc}.c)}}\label{sbtr}
\setcounter{equation}{0}
Let $\,\jv\hs$ be the space (\ref{spv}) for a given integer $\,m\ge2$, and let
us fix $\,\th,\ha\in\bbR\,$ with $\,\th\ne\ha$. By Lemma~\ref{onedi}, there exists
$\,Q\in\jv\smallsetminus\{0\}\,$ satisfying (\ref{qdq}.b) on the interval
$\,I\hs$ with the endpoints $\,\th,\ha$, and such $\,Q\,$ is unique up to a
constant factor.

For these $\,Q,I\nh$, (\ref{qdq}.e) holds if and only if $\,(\th,\ha)\,$ satisfies
one of the three conditions in (\ref{trc}). (See Theorem~\ref{facto}.) The question
of finding direct descriptions of the three sets in the $\,\th\ha\hs$-plane
$\,\rto$ defined by (\ref{trc}.a), (\ref{trc}.b) and, respectively, (\ref{trc}.c),
has an obvious answer for (\ref{trc}.b), the set being the hyperbola $\,\ha=\hu$
(that is, $\,\th\ha=\th+\ha$). As for (\ref{trc}.a), Proposition~\ref{sgrts} yields an
answer: (\ref{trc}.a) defines for odd $\,m\,$ the two-point set
$\,\{(1,\sx),(\sx,1)\}$, and for even $\,m\,$ the empty set.

Our real interest lies, however, in those $\,(\th,\ha)\,$ for which $\,Q,I\hs$
chosen above satisfy {\it all five\/} conditions in (\ref{qdq}). This leads to
switching our focus from the three solution sets in $\,\rto$ (see the last
paragraph) to their respective subsets obtained by imposing on $\,Q,I$ also
the positivity condition of \S\ref{posi}. For (\ref{trc}.a), the resulting
subset is empty (Proposition~\ref{noone}), while in the case of (\ref{trc}.b)
an explicit description of that subset is provided by Proposition~\ref{pyuzw}.

That new focus also explains why, unlike the approach to (\ref{trc}.a) -- (\ref{trc}.b) outlined above, our discussion of (\ref{trc}.c) bypasses the step of
first describing the set given by (\ref{trc}.c) alone. Instead, we proceed
directly to discuss the subset of the $\,\th\ha\hs$-plane $\,\rto$ defined by
requiring that $\,Q,I\hs$ corresponding to the given $\,(\th,\ha)\,$ with
$\,\th\ne\ha\,$ satisfy {\it both\/} (\ref{trc}.c) {\it and\/} (\ref{qdq}). This
subset turns out to be empty for even $\,m\hs$, while for odd $\,m\,$ it is the
union of the curve $\,\vg\,\subset\,(-\infty,0)\times(0,1)\,$ described in
Proposition~\ref{crvsg}(ii) and the image of $\,\vg$ under the symmetry
$\,(\th,\ha)\mapsto(\ha,\th)$.

In fact, a point $\,(\th,\ha)\,$ in this subset can never satisfy any of the
following conditions: $\,\th=1\,$ or $\,\ha=1\,$ (by (\ref{trc}.c)),
$\,\th\ha=0$, or $\,\th<1<\ha$, or $\,\ha<1<\th\,$ (see (ii) in \S\ref{posi}). Also,
$\,\th,\ha\,$ cannot both lie in $\,(-\infty,0)$, $\,(0,1)$, or
$\,(1,\infty)\,$ (by (\ref{trc}.c) and (\ref{sgt}.b\hs-\hs d)). This leaves
$\,\th<0<\ha<1\,$ or $\,\ha<0<\th<1\,$ as the only possibilities, so that our
claim follows from Proposition~\ref{crvsg} and symmetry of $\,T$.

\section{A synopsis of conditions \hs{\rm(\ref{qdq})}}\label{syno}
\setcounter{equation}{0}
The {\it moduli curve\/} defined in \S\ref{modu} is the subset $\,\mathcal{C}\hs$ of
$\,\rto\nh$, depending on an integer $\,m\ge2$, and formed by all pairs
$\,(\th,\ha)\,$ such that either $\,(\th,\ha)=(0\hs,\nh0)$, or $\,\th<\ha\,$
and all five conditions in (\ref{qdq}) are satisfied by the interval
$\,I=[\th,\ha]\,$ and some function $\,Q\in\jv\smallsetminus\{0\}$, with
$\,\jv\,$ as in (\ref{spv}). In view of Lemma~\ref{onedi}, $\,Q\,$ then is unique up
to a constant factor, while Lemma~\ref{abceq} and Remark~\ref{ptone} provide a choice
of such $\,Q\,$ for which $\,\ax,\bx,\cx\,$ in (\ref{qef}) are specific rational
functions of $\,\th,\ha$. Namely, they are given by (\ref{abc}) (if
$\,1\notin I$), or by $\,(\ax,\bx,\cx)=(-\hs E(\ps),\hs1,\hs0)\,$ (if
$\,\{\th,\ha\}=\{1,\ps\}$).

The next result clearly implies (\ref{xii}) for the sets $\,\,\mathsf{I}\hs$,
\ \di\ \ and $\,\,\mathsf{X}\,\,$ defined in \S\ref{stat}. In other words, the
definition of $\,\mathcal{C}\hs$ (see above) agrees with the explicit
description $\,\mathcal{C}\hs$ given in \S\ref{stat}.
\begin{thm}\label{modcu}Let\/ $\,\hu=\th/(\th-1)\,$ for\/ $\,\th\ne1$. A pair
$\,(\th,\ha)\in\rto$ lies in the moduli curve\/ $\,\mathcal{C}$, defined as
above for a fixed integer\/ $\,m\ge2$, if and only if one of the following
three cases occurs\/{\rm:}
\begin{enumerate}
\item[a)] $m\,$ is even, $\,\th\in(-\infty,0\hs]\cup(1,2)$, and\/
$\,\ha=\hu\nh$.
\item[b)] $m\,$ is odd, $\,\th\in(-\infty,\tz)\cup(\tw,0\hs]\cup(1,2)$, and\/
$\,\ha=\hu$, for\/ $\,\tz,\tw\,$ as in {\rm\S\ref{anth}}.
\item[c)] $m\,$ is odd and\/ $\,(\th,\ha)\in\vg$, where\/
$\,\vg\,\subset\,(-\infty,0)\times(0,1)\,$ is the connected real-analytic
curve diffeomorphic to\/ $\,\bbR\hs$, described in {\rm Proposition~\ref{crvsg}(ii)}.
\end{enumerate}
\end{thm}
\begin{proof}Let $\,\th<\ha$. Each of (a), (b), (c) separately implies (\ref{qdq}) with $\,Q,I\hs$ as above. For (a), (b) this is clear:
$\,1\notin I$ (cf.\ (\ref{ast})), which yields (\ref{qdq}.a); $\,\ha=\hu\nh$, i.e., (\ref{trc}.b), implies (\ref{qdq}.e) (Theorem~\ref{facto}); Proposition~\ref{pyuzw} and Lemma~\ref{uuast} give (\ref{qdq}.c) -- (\ref{qdq}.d). Also, (c) yields (\ref{qdq}) by Proposition~\ref{tzrok}.

Conversely, let $\,Q,I\hs$ as above satisfy (\ref{qdq}). Theorem~\ref{facto} then
gives (\ref{trc}.b) or (\ref{trc}.c), as Proposition~\ref{noone} excludes (\ref{trc}.a).
In case (\ref{trc}.b), we must have (a) or (b) (Proposition~\ref{pyuzw} and
Lemma~\ref{uuast}), while (\ref{trc}.c) yields (c), cf.\ the fourth paragraph
of \S\ref{sbtr}. This completes the proof.
\end{proof}
\begin{cor}\label{compo}If the integer\/ $\,m\ge2\,$ is even, or odd, then the
two sets $\,\,\di\hs$, $\,\mathsf{I}\,\,$ or, respectively, three sets
$\,\,\mathsf{X}\hs$, $\,\di\hs$, $\,\mathsf{I}\hs$, defined in\/
{\rm\S\ref{stat}}, are the connected components of the moduli curve
$\,\mathcal{C}$.
\end{cor}
This is clear from (\ref{xii}): the two/three sets are connected, relatively
open in $\,\mathcal{C}\hs$ and, for odd $\,m\hs$, the $\,\mathcal{T}\nh$-beam
$\,\xat\,$ coincides with the set $\,\vg\hs$ in Proposition~\ref{crvsg}, and
so, by Lemma~\ref{xxast}(b), $\,\xat\,$ intersects the $\,\mathcal{H}$-beam
$\,\xah\,$ only at $\,(\xj,\xj^*)$, and does not intersect
$\,\,\mathsf{I}\,\,$ or $\,\,\di\hs$.
\begin{rem}\label{trnsv}
For odd $\,m\hs$, the intersection of the two beams of
$\,\,\mathsf{X}\,\,$ at their unique common point $\,(\xj,\xj^*)\,$ is
transverse, by Lemma~\ref{xxast}(d).
\end{rem}

\begin{rem}\label{apole}
Let (\ref{qdq}) be satisfied by $\,I=[\th,\ha]\,$ and a
rational function $\,Q\,$ of the form (\ref{qef}) with some
$\,\ax,\bx,\cx\in\bbR\hs$. Then $\,\cx\nh\ne0$, so that, by (\ref{fet}), $\,Q\,$
is analytic on $\,\bbR\smallsetminus\{1\}\,$ and has a pole at $\,1$.

In fact, suppose on the contrary that $\,\cx\nh=0$. Thus, $\,m\,$ is odd and
$\,(\th,\ha)\in\vg$, for $\,\vg\,$ as in Proposition~\ref{crvsg}(ii) (or else, as
our assumption gives $\,(\th,\ha)\in\mathcal{C}\smallsetminus\{(0\hs,\nh0)\}$,
Theorem~\ref{modcu} would yield $\,\ha=\hu$ with $\,0\ne\th\ne2$, and hence
$\,\cx\nh\ne0\,$ by (\ref{fug}) and (\ref{fet})). Also, (\ref{qdq}.b), (\ref{qdq}.c) and (\ref{qef}) with $\,\cx\nh=0\,$ imply (\ref{qta}) with $\,\bx\ne0$. Therefore,
since $\,\th\ne1\ne\ha$, Lemma~\ref{czero}(i) shows that
$\,E(\th)=E(\ha)=-\hs\ax/\bx\,$ as well as
$\,\dot Q(\ps)\,=\,(\ps-1)\hs\bx\dot E(\ps)$, with
$\,(\,\,)\dot{\,}=\,d/d\ps$, for both $\,\ps=\th\,$ and $\,\ps=\ha$. Now (\ref{qdq}.e) yields $\,(\th-1)\hs\dot E(\th)=(1-\ha)\hs\dot E(\ha)$. Since
$\,\th<0<\ha<1\,$ (due to the definition of $\,\vg$), (\ref{emo}) along with
the three lines preceding it (from now on referred to simply
as \S\ref{mopr}) give $\,E(0)<E(\ha)<0<\dot E(\ha)$, while $\,\th-1<0<1-\ha$. Thus,
by the last equality, $\,\dot E(\th)<0$, and $\,E(0)<E(\th)<0\,$ as
$\,E(\th)=E(\ha)$. From $\,\th<0\,$ and $\,\dot E(\th)<0\,$ we in turn get
$\,\th>\wu>\zu\,$ (see \S\ref{mopr}), and so $\,\th>\wu>\xj$, as $\,\zu>\tz\,$
(Remark~\ref{nalgy}) and $\,\tz>\xj\,$ (Lemma~\ref{xxast}(b)). Similarly, using (\ref{ast}.a) for $\,\ps=\wu\,$ and noting that $\,\xj^*>\hatt$ as
$\,\ps=\wu>\xj$, cf.\ Lemma~\ref{uuast}, while $\,F<0\,$ on $\,(-\infty,0)\,$ (see
\S\ref{mono}), we get $\,E(\xj^*)>E(\hatt)=E(\ps)-F(\ps)>E(\ps)>E(\th)=E(\ha)\,$ for
$\,\ps=\wu\,$ from (\ref{fta}.ii) and the monotonicity properties of $\,E\,$
listed in \S\ref{mopr}, and so, again from \S\ref{mopr} and Lemma~\ref{uuast},
$\,0<\ha<\xj^*<1$.

The function $\,\qx\,$ given by (\ref{qvu}) is clearly increasing (or,
decreasing) as a function of $\,\th\,$ (or, $\,\ha$) alone in the region where
$\,\th<0<\ha<1$. As $\,\xj<\th<0\,$ and $\,0<\ha<\xj^*$, this implies that the
value of $\,\qx\,$ at $\,(\th,\ha)\,$ is {\it greater} than $\,\qx_*$, its
value at $\,(\xj,\xj^*)$. Since $\,(\th,\ha)\in\vg\,$ and $\,\qx\le\qx_*$ on
$\,\vg\,$ (see the last sentence in Proposition~\ref{crvsg}(ii)), we now obtain
a contradiction. Therefore, $\,\cx\nh\ne0$.
\end{rem}

\section{The rational function $\,\px$}\label{trfp}
\setcounter{equation}{0}
Given an integer $\,m\ge2$, we let $\,\px\,$ stand for the rational function
of the variables $\,\th,\ha$, defined by the formula in (\ref{pum})
with $\,\lx\,$ as in (\ref{lba}) for $\,\ax,\bx,\vs\,$ given by (\ref{abc}) and (\ref{fet}) -- (\ref{esg}). For later convenience, we modify this definition by
declaring the value of $\,\px\,$ at $\,(0\hs,\nh0)\,$ to be $\,0$.

Thus, $\,\px\,$ is real-analytic everywhere in $\,\rto$ with a possible
exception of those points $\,(\th,\ha)\in\rto\,$ for which $\,\th\in\{0,2\}\,$
(cf.\ (\ref{pum}), or $\,\th=1$, or $\,\ha=1\,$ (since (\ref{lba}) involves
$\,\ax,\bx\,$ with (\ref{abc}), and $\,F$, given by (\ref{fet}), has a pole at
$\,1$), or, finally, $\,\th\ne1\ne\ha\,$ and $\,\ax=0\,$ at $\,(\th,\ha)\,$
(as $\,\ax\,$ appears in the denominator of (\ref{lba})).

We are interested in the restriction of $\,\px\,$ to the moduli curve
$\,\mathcal{C}\hs$ (see \S\ref{syno}). Of the singularities just listed, only
$\,(0\hs,\nh0)\,$ lies on $\,\mathcal{C}\hs$ if $\,m\,$ is odd, and just two,
$\,(0\hs,\nh0)\,$ and $\,(\tz,\zj)$, lie on $\,\mathcal{C}\hs$ when $\,m\,$ is
even, with $\,\tz<0\,$ defined as in \S\ref{anth} and $\,\zj=\tz/(\tz-1)$. The
singularity at $\,(\tz,\zj)\,$ arises since $\,\ax=0\,$ there.

In fact, for $\,(\th,\ha)\in\mathcal{C}\smallsetminus\{(0\hs,\nh0)\}\,$ we have
$\,\th,\ha\notin\{0,1,2\}\,$ by Theorem~\ref{modcu} and Lemma~\ref{uuast}. Hence,
from (\ref{abc}) and (\ref{fet}), $\,\ax=0\,$ at such $\,(\th,\ha)\,$ if and only
if $\,E/F\,$ has equal values at $\,\th\,$ and $\,\ha$. This in turn excludes
the possibility that either $\,1<\th<2<\ha$, or $\,m\,$ is odd and
$\,\th<0<\ha<1\,$ (cf.\ the last two lines in (\ref{mon})), so that, by Lemma~\ref{uuast}, the only case still allowed in Theorem~\ref{modcu} is (a) with
$\,\th\in(-\infty,0\hs]\,$ (and so $\,m\,$ is even, while $\,\ha=\hu$). Our
claim about $\,(\tz,\zj)\,$ and $\,\ax\,$ now follows from (\ref{fug}) -- (\ref{god}) (and (\ref{fet})).
\begin{rem}\label{kemac}
Given an integer $\,m\ge2$, we have $\,\kx=\ve\hs m\ax/\y\,$ and (\ref{pde})
whenever
$\,(\th,\ha)\in\mathcal{C}\smallsetminus\{(0\hs,\nh0)\}\,$ and $\,Q\,$ is a
function in the space (\ref{spv}) satisfying (\ref{qdq}) on $\,I=[\th,\ha]$, while
$\,\ax\in\bbR\,$ is determined by $\,Q\,$ via (\ref{qef}), $\,\ve=\pm\hs1\,$ and
$\,\y,\ah\,$ are nonzero constants with $\,dQ/d\ps=-\hs2\ah\y\,$ at
$\,\ps=\th$, and, finally, either $\,\ah=-\hs\ve\hs\px\kx/2\,$ (where
$\,\px\,$ stands for the value at $\,(\th,\ha)\,$ of the rational function
$\,\px$), or $\,m\,$ is even, $\,(\th,\ha)=(\tz,\zj)\,$ and $\,\kx=0$. In
fact, if $\,m\,$ is even and $\,(\th,\ha)=(\tz,\zj)$, this follows since, as
we just saw, $\,\ax\,$ then equals $\,0\,$ at $\,(\tz,\zj)$. Otherwise,
$\,\ax\ne0\,$ at $\,(\th,\ha)\,$ and $\,\th,\ha\notin\{0,1,2\}\,$ (see above),
so that, by (\ref{pde}.i) and (\ref{qdq}.d), we have $\,\px\ne0\,$ and
$\,m\px=\dot Q(\th)/\ax=-\hs2\ah\y/\ax=\ve\hs\px\kx\y/\ax$, as required.
\end{rem}

By (\ref{xii}), a substantial part of the moduli curve $\,\mathcal{C}\hs$ is
contained in the hyperbola $\,\mathcal{H}\,$ given by $\,\ha=\hu$, where
$\,\hu=\th/(\th-1)$. It is therefore useful to introduce a rational function
$\,\pw\,$ of the real variable $\,\th\,$ which is the restriction of $\,\px\,$
to $\,\mathcal{H}$, that is, the result of substituting $\,\hu$ for $\,\ha\,$
in the rational expression for $\,\px\,$ in terms of $\,\th\,$ and $\,\ha$.
Then, with $\,\vs,G\,$ as in (\ref{esg}) and (\ref{gef}), we have
\begin{equation}
\th(\th-2)\,\pw(\th)\,=\,2(2\,-\,1/m)\hs(\th-1)\,\lj1\,
+\,\vs(0)/G(\th)\rj\,\,-\,\,\th^2\nh,\label{uup}
\end{equation}
in the sense of equality between rational functions of $\,\th$. In fact, for
those $\,\th\,$ for which $\,\pw(\th)\,$ defined by (\ref{uup}) makes sense,
it coincides with the number $\,\px\,$ in (\ref{pum}) for
$\,\lx=(2-1/m)\hs\lj1\,+\,\vs(0)/G(\th)\rj$, that is, for $\,\lx\,$ given by
(\ref{lba}) with $\,\ha=\hu$, as one sees evaluating $\,\bx/\ax\,$ from
(\ref{fug}).

A trivial argument (see (c) in \S\ref{sprp}) shows that the rational function
$\,\pw\,$ defined by (\ref{uup}) is analytic at $\,0\,$ and $\,\pw(0)=0$. Thus,
our convention about $\,\px\,$ at $\,(0\hs,\nh0)\,$ requires no further
modification of $\,\pw$.
\begin{rem}\label{patxx}
Let $\,\qx\,$ and $\,\qx_*$ be as in (\ref{qvu}) for a fixed odd integer
$\,m\ge3$. The $\,\mathcal{T}\nh$-beam of the moduli curve  (cf.\
\S\ref{stat}) is the set $\,\vg\,$ in Proposition~\ref{crvsg}, while
$\,\px=\qx\,$ at any point $\,(\th,\ha)\in\vg\,$ since, by
Lemma~\ref{mumil}(ii), $\,\lx\,$ in (\ref{pum}) then may be replaced with
$\,\mx$. Thus (cf.\ the final clause in Proposition~\ref{crvsg}(ii)), the
value of $\,\px\,$ at $\,(\xj,\xj^*)\,$ is $\,\qx_*\in(-\hs1,0)$.
\end{rem}

\section{Some properties of $\,\pw$}\label{sprp}
\setcounter{equation}{0}
Let $\,\dot\pw=\hs d\pw/d\th\,$ for the rational function $\,\pw\,$ given by (\ref{uup}) with a fixed integer $\,m\ge2$. Then $\,\pw\,$ satisfies the
differential equation
\begin{equation}
\th(\th-1)\,\dot\pw\,=\,m\th\hs\pw^2\hs\,+\,\,(m-1)(\th-2)\hs\pw\,\,
-\,\,\th\,,\label{udp}
\end{equation}
where $\,\pw,\,\dot\pw\,$ stand for $\,\pw(\th),\,\dot\pw(\th)$. Also, for
$\,\tz,\tw\,$ as in \rm\S\ref{anth} and $\,\hu=\th/(\th-1)$,
\begin{enumerate}
\item[a)] If $\,m\,$ is odd, $\,\pw\,$ is analytic everywhere in $\,\bbR\,$
and $\,\pw(\tw)=\pw(\wj)=0$.
\item[b)] If $\,m\,$ is even, $\,\pw\,$ has just two real poles, at $\,\tz\,$
and $\,\zj$.
\item[c)] $\pw(0)=\pw(2)=0$, $\,\,\pw(1)=1$, $\,\,\dot\pw(0)=1/(3-2m)$,
$\,\,\dot\pw(1)=2(1-m)/m\hs$.
\item[d)] $\pw(\hu)=-\hs\pw(\th)\,$ for all
$\,\th\in\bbR\smallsetminus\{1\}\,$ at which $\,\pw\,$ is analytic.
\end{enumerate}
In fact, (d) is immediate if one replaces $\,\th\,$ in (\ref{uup}) by $\,\hu$
and then uses (\ref{ast}) and (\ref{gta}.a). If $\,m\,$ is odd,
$\,\dot G(\tw)=0\,$  by the definition of $\,\tw$, and so (\ref{gta}.c) with
$\,\ps=\tw\,$ gives $\,\vs(0)/G(\tw)=\varLambda(\tw)/[2(2m-1)(\tw-1)]
=-\hs1+m\tw^2/[2(2m-1)(\tw-1)]$, where we also used (\ref{dtf}.b). Replacing
$\,\vs(0)/G(\th)\,$ in (\ref{uup}) by this last expression for $\,\th=\tw$, we
obtain $\,\pw(\tw)=0\,$ (as $\,\tw<0$) and, from (d), $\,\pw(\wj)=0$. Next,
dividing (\ref{uup}) by $\,\th(\th-2)\,$ we find, using l'Hospital's rule and
(\ref{gta}.d) -- (\ref{gta}.f), that $\,\pw(\th)\to0\,$ as $\,\th\to0\,$ or
$\,\th\to2$. Thus, $\,\pw(0)=\pw(2)=0$. The comment on the zeros of $\,G\,$
following (\ref{god}) now gives $\,\pw(1)=1\,$ (by (\ref{uup})), and hence
(a), (b).

With $\,G\,$ expressed in terms of the function
$\,\varphi(\th)=2(2m-1)(\th-1)\hs\vs(0)/G(\th)$, equation (\ref{gta}.c) becomes
$\,\th(\th-1)(\th-2)\dot\varphi+\lj\varLambda(\th)-\varphi-\th(\th-2)\rj
\hs\varphi=0$, where $\,\varphi,\dot\varphi\,$ stand for
$\,\varphi(\th),\hs\dot\varphi(\th)$. Replacing $\,\varphi\,$ by
$\,\varLambda(\th)\hs+\hs m\hs\th(\th-2)\hs\pw\,$ (which equals $\,\varphi\,$
in view of (\ref{uup}) and (\ref{dtf}.b)), and then
substituting for $\,\varLambda(\th)\,$ the expression in (\ref{dtf}.b), we can
further rewrite this as an equation imposed on $\,\pw$. That equation is
easily verified to be (\ref{udp}) with both sides multiplied by
$\,m\hs\th(\th-2)^2$.

Finally, the values of $\,\dot\pw\,$ required in (c) are easily obtained by
differentiating (\ref{udp}) at $\,\th=0\,$ or $\,\th=1\,$ and using the
relations $\,\pw(0)=0$, $\,\pw(1)=1$.

\section{Monotonicity intervals for $\,\pw$}\label{mivp}
\setcounter{equation}{0}
The rational function $\,\pw\,$ defined by (\ref{uup}) with a fixed integer
$\,m\ge2\,$ has a nonzero derivative at every $\,\th\in\bbR\,$ except
$\,\tz,\zj$ (for even $\,m$), or $\,\th_\pm,\th_\pm^*$ (for odd $\,m$), with
$\,\hu=\th/(\th-1)\,$ if $\,\th\ne1\,$ (cf.\ (\ref{ast})), $\,\tz\,$ as in (\ref{gev}), and $\,\th_\pm$ described below. The values/limits of $\,\pw\,$ at
selected points, along with its strict monotonicity types on the intervening
intervals, marked by slanted arrows, are listed below, with $\,\tz,\tw,\xj\,$
as in (\ref{gev}) -- (\ref{god}) and Lemma~\ref{xxast}(b). First, for even $\,m\hs$,
\begin{equation}
\begin{array}{lccccccc}\arraycolsep9pt
\mathrm{value\ or\ limit\ at\hskip-3.2pt:}\quad&0&&\zj&&1&&2\\
[1pt]
\hline
&&&&&&&\\[-9.7pt]
\mathrm{for}\hskip4.5pt\pw\hskip4pt\mathrm{(}m\hskip4pt\mathrm{even)
\hskip-3.2pt:}&0&\seam&{-\infty}\vum{+\infty}&\seam&1&\seam&0
\end{array}
\label{vlp}
\end{equation}
If $\,m\,$ is odd, $\,\pw\,$ restricted to $\,[\hs0,2\hs]\,$ reaches its
extrema at unique points $\,\th_\pm$ with
\begin{equation}
0<\th_-\hskip-1.3pt<\wj\hskip-1.3pt<\zj\hskip-1.3pt
<\xj^*\hskip-1.3pt<\th_+\hskip-1.3pt<1\hs\,\mathrm{\ and\ }
\,-1/m<\pw(\th_-)<0<1<\pw(\th_+)\hs.\label{uwz}
\end{equation}
(See also (vi) in \S\ref{thrd}.) Here is the corresponding diagram:
\begin{equation}
\begin{array}{lccccccccccccccc}\arraycolsep0pt
\mathrm{value\ at\hskip-3.2pt:}\quad&0&&\th_-&&\wj&&\zj&&\xj^*&&\th_+&&1&&2\\
[1pt]
\hline
&&&&&&&&&&&&&&&\\[-9.7pt]
\mathrm{of}\hskip4.5pt\pw\hskip4pt\mathrm{(}m\hskip4pt\mathrm{odd)\hskip-3.2pt:
}&0&\hskip-2pt\seam\hskip-2pt&\pw(\th_-)&\hskip-2pt\neam\hskip-2pt&0&\hskip-2pt
\neam\hskip-2pt&\tz/(\tz-2)&\hskip-2pt\neam\hskip-2pt&\qx_*&\hskip-2pt\neam
\hskip-2pt&\pw(\th_+)
&\hskip-2pt\seam\hskip-2pt&1&\seam&0
\end{array}
\label{vpo}
\end{equation}
\begin{rem}\label{behvp}
In view of (d),\hs(a),\hs(b) of \S\ref{sprp} and Lemma~\ref{uuast}, the monotonicity intervals of $\,\pw\,$ on the {\it whole\/} real
line can be easily determined using (\ref{vlp}) -- (\ref{vpo}). Specifically,
$\,\pw\,$ always
decreases from $\,0\,$ to $\,-\hs1\,$ on $\,[\hs2,\infty)$, that is, forms a
decreasing diffeomorphism $\,[\hs2,\infty)\to(-\hs1,0\hs]$. Similarly, when
$\,m\,$ is even, $\,\pw\,$ decreases on $\,(-\infty,\tz)\,$ (or, on
$\,(\tz,0\hs]$) from $\,-\hs1\,$ to $\,-\infty\,$ (or, respectively, from
$\,\infty\,$ to $\,0$). Finally, if $\,m\,$ is odd, $\,\pw\,$ decreases on
$\,(-\infty,\th_+^*]\,$ from $\,-\hs1\,$ to $\,\pw(\th_+^*)=-\pw(\th_+)$,
increases on $\,[\th_+^*,\tw]\,$ from $\,-\pw(\th_+)\,$ to $\,0$, and then
continues increasing on $\,[\tw,\th_-^*]$, from $\,0\,$ to
$\,\pw(\th_-^*)=-\pw(\th_-)$, while on $\,[\th_-^*,0\hs]\,$ it decreases from
$\,-\pw(\th_-)\,$ to $\,0$.
\end{rem}

\section{Proofs of the above claims}\label{pfac}
\setcounter{equation}{0}
For any $\,\th\ne0\,$ the right-hand side of (\ref{udp}) is a quadratic
polynomial in $\,\pw\,$ having real roots $\,\pw_\pm(\th)\,$ with
$\,\pw_-(\th)<\pw_+(\th)$. Clearly, $\,\th\mapsto\pw_\pm(\th)\,$ are
real-analytic functions on $\,(-\infty,0)\cup(0,\infty)$. As shown below, for
$\,\dot\pw_\pm=\hs d\pw_\pm/d\th\,$ and (in (iii), (iv)) for any
$\,\th\in\bbR\smallsetminus\{0,1\}\,$ at which $\,\pw(\th)\,$ is defined, cf.\
(a), (b) in \S\ref{sprp},
\begin{enumerate}
\item[i)] $\,\pw_-<0<\pw_+$,\hskip13pt(ii)\hskip6pt$\dot\pw_\pm<0$,
\hskip13pt(iii)\hskip6pt$\,(\th-1)\hs\dot\pw(\th)<0\,$\hskip8ptif and only
if\hskip8pt$\,\pw_-(\th)<\pw(\th)<\pw_+(\th)$.
\item[iv)] $\,\dot\pw(\th)=0\,$\hskip8ptif and only
if\hskip8pt$\,\pw(\th)=\pw_-(\th)\,$\hskip8ptor\hskip8pt$\,\pw(\th)
=\pw_+(\th)$.
\item[v)] $\,\pw_\mp(0^\pm)=0$\hskip8ptand\hskip8pt$\dot\pw_\mp(0^\pm)
=[2(1-m)]^{-1}$\hskip14pt(one-sided limits at $\,0$),
\item[vi)] $\,\pw_+(1)=1\hs$,\hskip8ptwhile\hskip8pt$\pw_-(1)
=-\hs1/m$\hskip8ptand\hskip8pt$\dot\pw_+(1)=-\hs2(m-1)/(m+1)$.
\end{enumerate}
In addition to (v), $\,\pw_+$ on $\,(-\infty,0)\,$ and $\,\pw_-$ on
$\,(0,\infty)\,$ are restrictions of a single analytic function on $\,\bbR\,$
with the value $\,0\,$ at $\,0$. The strict-monotonicity intervals (marked by
slanted arrows) and some limits of $\,\pw_\pm$ appear in the diagram
\begin{equation}
\begin{array}{lcccccccc}\arraycolsep9pt
\mathrm{value\ or\ limit\ at\hskip-3.2pt:}&{-\infty}&&0&&1&&+\infty\\
[1pt]
\hline
&&&&&&&&\\[-9.7pt]
\mathrm{for}\hskip4.5pt\pw_+\hskip4pt\mathrm{(any}\hskip4ptm
\mathrm{)\hskip-3.2pt:}&1/m&\seam&0\vum{+\infty}&\seam&1&\seam&1/m\\
\mathrm{for}\hskip4.5pt\pw_-\hskip4pt\mathrm{(any}\hskip4ptm
\mathrm{)\hskip-3.2pt:}&{-1}&\seam&{-\infty}\vum\hs0&\seam&{-1/m}&\seam&{-1}
\end{array}
\label{vpp}
\end{equation}
In fact, (i) is obvious as $\,\pw_-(\th)\pw_+(\th)=-\hs1/m<0$, while (\ref{udp})
gives (iii), (iv). Next, equating the right-hand side of (\ref{udp}) to zero
we obtain$\,(\pw+1)(m\pw-1)\hs\th=2(m-1)\hs\pw$. This defines a set in the
$\,\th\pw\hs$-plane, namely, the union of $\,\{(0\hs,\nh0)\}\,$ and the graphs
of $\,\pw_\pm$, which, at the same time, forms the graph of the rational
function $\,\th\,$ of the variable $\,\pw\,$ with
$\,\th=2(m-1)\hs\pw/[(\pw+1)(m\pw-1)]$. The latter function has a negative
derivative except at the two poles $\,\pw=-\hs1\,$ and $\,\pw=1/m\hs$, and tends
to $\,0\,$ as $\,\pw\to\pm\infty$, so that (ii), (v), (vi) and (\ref{vpp})
follow easily.

Next, $\,\pw(0),\pw(1),\pw(2)\,$ are given by (c) in \S\ref{sprp}, and, if
$\,m\,$ is even, $\,\pw(\th)\to\pm\infty\,$ as $\,\th\to\tz^\pm$ by
(\ref{gev}) and (\ref{uup}), so that (d) in \S\ref{sprp} gives
$\,\pw(\th)\to\pm\infty\,$ as $\,\th-\zj\to0^\pm\hskip-1.5pt$. If $\,m\,$ is
odd, $\,1+\vs(0)/G(\tz)=0\,$ (see (\ref{god}), (\ref{sto}.ii)), and so
$\,\pw(\zj)=\tz/(2-\tz)\,$ by (\ref{uup}). To prove (\ref{vlp}) -- (\ref{vpo}), we
now only need to show that $\,\dot\pw<0\,$ on
$\,(0,\zj)\cup(\zj,1)\cup(1,2\hs]\,$ for even $\,m\hs$, while, for odd $\,m\hs$,
there exist $\,\th_\pm\in\bbR\,$ with (\ref{uwz}) such that
$\,\dot\pw<0\,$ on $\,(0,\th_-)\cup(\th_+,1)\cup(1,2\hs]\,$ and
$\,\dot\pw>0\,$ on $\,(\th_-,\th_+)$. (Note that, according to (c) in \S\ref{sprp},
$\,\dot\pw<0\,$ at $\,0\,$ and $\,1$.)

By using (v) -- (vi) above, (\ref{vpp}) and (c) in \S\ref{sprp} to find the value of
$\,\varPhi=\pw-\pw_\pm$ at $\,1\,$ or its right-sided limit at $\,0\,$ (and
the same for $\,\dot\varPhi$, as needed), we see that $\,\pw-\pw_+$ changes
sign at $\,1$, from $\,+\,$ to $\,-\hs$, while $\,\pw-\pw_->0\,$ at $\,1\,$
and $\,\pw-\pw_\pm<0\,$ at every $\,\th>0\,$ close to $\,0$. Also,
$\,\pw-\pw_+<0<\pw-\pw_-$ at $\,2\,$ by (i), since (c) in \S\ref{sprp} states that
$\,\pw(2)=0$.

In view of (ii), (iv) above and (a), (b) in \S\ref{sprp}, the assumptions of Remark~\ref{zerpo} are satisfied by $\,\varPhi=\pw-\pw_\pm$ on $\,\mathcal{I}=(0,1)\,$ (for
odd $\,m$), as well as on $\,\mathcal{I}=(0,\zj)\,$ or $\,\mathcal{I}=(\zj,1)\,$ (for
even $\,m$), and on $\,\mathcal{I}=(1,2\hs]\,$ (for all $\,m$). In each case, the inequalities of the last paragraph lead, as shown below, to a
unique choice between the two alternatives allowed in the conclusion of Remark~\ref{zerpo}.

First, as $\,\pw-\pw_+$ (or, $\,\pw-\pw_-$) restricted to $\,(1,2\hs]\,$ is
negative (or, respectively, positive) near both endpoints, Remark~\ref{zerpo}
implies that $\,\hs\sgn\hs(\pw-\pw_\pm)\,$ is constant on $\,(1,2\hs]$, and
hence $\,\pw_-<\pw<\pw_+$ on $\,(1,2\hs]$. By (iii), this yields
$\,\dot\pw<0\,$ on $\,(1,2\hs]$.

Secondly, as $\,\pw-\pw_\pm\,$ restricted to $\,(0,1)\,$ is negative near the
endpoint $\,0\,$ and positive near the endpoint $\,1$, Remark~\ref{zerpo} gives
rise to two different cases, depending on $\,m\hs$. If $\,m\,$ is even, the
infinite limits of $\,\pw\,$ at $\,\zj\hskip-1.5pt$, already verified to be
those required in (\ref{vlp}), show that $\,\pw-\pw_\pm\,$ on
$\,\mathcal{I}=(0,\zj)\,$ (or, $\,\mathcal{I}=(\zj,1)$) is negative (or, respectively,
positive) near both endpoints, and hence, by Remark~\ref{zerpo}, it is so
everywhere in $\,\mathcal{I}$. Since this applies to both signs $\,\pm\hs$, (iii)
and (iv) give $\,\dot\pw<0$, for even $\,m\hs$, both on $\,(0,\zj)\,$ and
$\,(\zj,1)$. If $\,m\,$ is odd, however, $\,\pw-\pw_\pm\,$ is of class $\,C^1$
everywhere in $\,(0,1)$, and so Remark~\ref{zerpo} implies the existence of unique
points $\,\th_\pm\in(0,1)\,$ such that $\,\hs\sgn\hs(\pw-\pw_\pm)\,$ at any
$\,\th\in(0,1)\,$ equals $\,\hs\sgn\hs(\th-\th_\pm)$. Hence, by (i),
$\,\pw=\pw_+>\pw_-$ at $\,\th_+$, and so $\,\hs\sgn\hs(\th_+-\th_-)=1$, that
is, $\,\th_-<\th_+$. Now, by (iii), (iv) and the last paragraph,
$\,\dot\pw<0\,$ on $\,(0,\th_-)\cup(\th_+,2\hs]\,$ and $\,\dot\pw>0\,$ on
$\,(\th_-,\th_+)$. Thus, as $\,\th\,$ increases from $\,0\,$ to $\,\th_-$,
then to $\,\th_+$, then to $\,1$, and finally to $\,2$, the value of $\,\pw\,$
decreases from $\,0\,$ to $\,\pw(\th_-)$, then increases to $\,\pw(\th_+)$,
then decreases to $\,1\,$ and after that continues decreasing to $\,0$. Any
point $\,\th\in(0,1)\,$ with $\,0\le\pw(\th)\le1\,$ must therefore lie in
$\,(\th_-,\th_+)$. This includes $\,\th=\wj\,$ and $\,\th=\xj^*\nh$, as
$\,\pw(\wj)=0\,$ (see (a) in \S\ref{sprp}) and
$\,\pw(\xj^*)=-\hs\qx_*\in(0,1)\,$ (by (d) in \S\ref{sprp} and
Remark~\ref{patxx} as $\,\pw(\xj)\,$ is the value of $\,\px\,$ at
$\,(\xj,\xj^*)$, cf.\ \S\ref{trfp}). We have thus proved
(\ref{vlp}), (\ref{vpo}) and the first part of (\ref{uwz}), since
$\,\wj<\zj<\xj^*$ by Lemma~\ref{uuast} with $\,\xj<\tz<\tw<0\,$ (cf.\ the lines preceding (\ref{god})). The
description just given of the monotonicity intervals of $\,\pw\,$ on
$\,[\hs0,2\hs]\,$ also shows that $\,\pw\,$ assumes its extrema in
$\,[\hs0,2\hs]\,$ at $\,\th_-$ and $\,\th_+$, while
$\,\pw(\th_-)<0<1<\pw(\th_+)$. Also, $\,\pw(\th_-)>-\hs1/m\,$ for odd
$\,m\hs$, since the minimum $\,\pw(\th_-)\,$ equals, by (iv), the value at
$\,\th_-$ of $\,\pw_-$ (not of $\,\pw_+$, as $\,\pw_-<\pw_+$ by (i)), and so,
by (\ref{vpp}), $\,\pw(\th_-)=\pw_-(\th_-)>\pw_-(1)=-\hs1/m\hs$.

\section{The values assumed by $\,\px\,$ on the moduli curve}\label{valu}
Let $\,m\ge2\,$ be a fixed integer. The restrictions to the connected
components of the moduli curve $\,\mathcal{C}\hs$ (see Corollary~\ref{compo})
of the functions $\,\dv\,$ and $\,\px\,$ defined in \S\ref{stat} and
\S\ref{trfp} have the following properties, which, as explained below, are
easy consequences of the results of the preceding sections. In the case of
$\,\dv$, (a) -- (e) simply repeat its definition from \S\ref{stat}, to provide
a convenient reference.
\begin{enumerate}
\item[a)] On the $\,\mathcal{T}\nh$-beam $\,\xat\,$ of $\,\,\mathsf{X}\hs$, when
$\,m\,$ is odd: $\,\dv=1\,$ and the range of $\,\px\,$ is $\,(-1,\qx_*]$, for
$\,\qx_*\in(-\hs1,0)\,$ as in (\ref{qvu}). Every value in $\,(-1,\qx_*)\,$ is
assumed by $\,\px\,$ exactly twice in $\,\xat$, while $\,\qx_*$ is assumed
just once, at $\,(\xj,\xj^*)$. Two different points of the $\,\xat\,$ have the
same value of $\,\px\,$ if and only if they are each other's images under the
involution (\ref{fsy}).
\item[b)] On the $\,\mathcal{H}$-beam $\,\xah$, when $\,m\,$ is odd: $\,\dv=1$, the
range of $\,\px\,$ is $\,[\pw(\th_+^*),\hs\tz/(2-\tz))$, while
$\,\pw(\th_+^*)<-\hs1<\tz/(2-\tz)<0\,$ by (\ref{uwz}) with (d), (a) in \S\ref{sprp},
and as $\,\tz<0\,$ (see \S\ref{anth}). Cf.\ (iii) in \S\ref{thrd}. The values in
$\,(\pw(\th_+^*),-\hs1)\,$ are assumed by $\,\px\,$ in $\,\xah\,$ twice, those
in the union $\,\{\pw(\th_+^*)\}\cup[-\hs1,\hs\tz/(2-\tz))\,$ just once.
\item[c)] On \ \di\hs, when $\,m\,$ is odd: $\,\dv=-\hs1\,$ and the range of
$\,\px\,$ is $\,[\hs0,\pw(\th_-^*)]$, with $\,0<\pw(\th_-^*)<1/m\hs$. Every value
in $\,[\hs0,\pw(\th_-^*)]\,$ is assumed by $\,\px\,$ twice
in\hskip5pt\di\hskip-.5pt, except
for $\,0\,$ and $\,\pw(\th_-^*)$, assumed just once, at $\,(0\hs,\nh0)\,$ and
$\,(\th_-^*,\th_-)$.
\item[d)] On \ \di\hs, when $\,m\,$ is even, $\,(\tz,\zj)\,$ is the only point
at which $\,\dv=0$, while $\,\dv=1\,$ (or, $\,\dv=-\hs1$) on the subset of \
\di\ \ formed by all $\,(\th,\hu)\in\,\hs\di\,\,$ with $\,\th<\tz\,$ (or,
respectively, $\,\tz<\th\le0$); that subset is mapped by $\,\px\,$ bijectively
onto $\,(-\infty,-1)\,$ (or, respectively, $\,[\hs0,\infty)$).
\item[e)] On $\,\,\mathsf{I}\hs$, for every $\,m\hs$, we have $\,\dv=1\,$ and
$\,\px:\mathsf{I}\to(0,1)\,$ is bijective.
\end{enumerate}
In fact, (a) is obvious from the last sentence in Proposition~\ref{crvsg}(ii),
along with the easily-verified invariance of $\,\qx\,$ under (\ref{fsy}) and
Remark~\ref{patxx}.

Moreover, $\,\ha=\hu$ for every point $\,(\th,\ha)\,$ of the moduli curve that
does not lie in the $\,\mathcal{T}$-beam $\,\xat\,$ (see \S\ref{syno}). Now (b) -- (e) are
immediate, since $\,\px\,$ at $\,(\th,\ha)\,$ then equals $\,\pw(\th)$, for
the function $\,\pw\,$ defined by (\ref{uup}), which has the limits/values and
monotonicity intervals are described in (\ref{vlp}) and (\ref{vpo}). The
inequalities in (b), (c) easily follow from (\ref{uwz}) and (d) in \S\ref{sprp}.
\begin{proof}[Proof of Theorem~\ref{ctbly}]Assertions (i), (ii) are obvious
from Definition~\ref{pratl} and (c), (d) above; that the set of
$\,\px\hs$-rational points in $\,\,\di\,\,$ is countably infinite for odd
$\,m\,$ as well follows from (c), as $\,[\hs0,\pw(\th_-^*)]\cap\bbQ\,$ is
infinite.
\end{proof}

\section{More on $\,\px\hs$-rationality}\label{more}
\setcounter{equation}{0}
Let $\,\mathcal{C}\hs$ again denote the moduli curve for a fixed integer
$\,m\ge2\,$ (see \S\ref{syno}). In \S\ref{stat}\ and \S\ref{trfp} we introduced two
functions on $\,\mathcal{C}$, namely, $\,\dv:\mathcal{C}\to\{-1,0,1\}\,$ and the
rational function $\,\px$. We also observed that $\,\px$, declared to be
$\,0\,$ at $\,(0\hs,\nh0)\in\mathcal{C}$, is defined everywhere in $\,\mathcal{C}$,
except at $\,(\tz,\zj)\,$ when $\,m\,$ is even. Both functions are involved in
Definition~\ref{pratl}, which describes a subset of the $\,\th\ha\hs$-plane
$\,\rto\nh$, contained in $\,\mathcal{C}$, and consisting of what we call the
$\,\px\hs$-rational points.
\begin{rem}\label{prati}
For a fixed integer $\,m\ge2$, a point $\,(\th,\ha)\in\rto$ with $\,\th<\ha\,$
is $\,\px\hs$-rational if and only if it can be used to construct a quadruple
$\,\mgmt\,$ with (\ref{zon}) or (\ref{zto}) as described in \S\ref{stat}. This
is in turn equivalent to the existence of objects required in (\ref{cea}) --
(\ref{nho}) for the given $\,\th,\ha\,$ and the function $\,Q\in\jv$, unique
up to a factor (cf.\ Lemma~\ref{onedi}), which satisfies (\ref{qdq}) on
$\,I=[\th,\ha]\,$ and is positive on the interior of $\,I\nh$.

Since the objects with (\ref{cea}) always exist, while the existence of those
in (\ref{nho}) is equivalent to (\ref{edo}), the above assertion will follow
once we show that the constants $\,\dv,\px\,$ in (\ref{edo}), defined by
(\ref{pde}), coincide with the values at $\,(\th,\ha)\,$ of the functions
$\,\dv\,$ and $\,\px\,$ defined in \S\ref{trfp} and \S\ref{stat}. This is
obvious for $\,\px\,$ when $\,\ax\,$ in (\ref{qef}) is nonzero (see the line
preceding (\ref{pum})), and for both $\,\px,\dv\,$ if $\,\ax=0\,$ (as
$\,(\th,\ha)\in\mathcal{C}\smallsetminus\{(0\hs,\nh0)\}\,$ due to the
definition of $\,\mathcal{C}$, and hence, according to \S\ref{trfp}, $\,m\,$
then is even and $\,(\th,\ha)=(\tz,\zj)$, so that $\,\dv=0\,$ while $\,\px\,$
is undefined, for either meaning of $\,\dv\,$ and $\,\px$). Finally, when
$\,\ax\ne0$, the inequality in (\ref{cea}) gives, for $\,\dv\,$ as in
(\ref{pde}.ii), $\,\dv=\,\sgn\,\varphi(\ps)\,$ whenever $\,\th<\ps<\ha$, with
$\,\varphi(\ps)=(\ps-1)\hs\ax Q(\ps)$. As $\,\th\ne1\,$ (cf.\
Remark~\ref{uvone}) and $\,\varphi(\th)=0\ne\dot\varphi(\th)\,$ unless
$\,\ax=0\,$ (see (\ref{qdq})), this yields $\,\dv=\,\sgn\,\dot\varphi(\th)$,
that is, $\,\dv=\,\sgn\hs[(\th-1)\hs\ax\dot Q(\th)]$, where
$\,(\,\,)\dot{\,}=\,d/d\ps$. Hence, by (\ref{pde}.i),
$\,\dv=\,\sgn\hs[(\th-1)\hs\px\hs]$, with either definition of $\,\px\,$ at
$\,(\th,\ha)$. (We already showed that both definitions agree.) This last
formula for $\,\dv\,$ is clearly consistent with (a) -- (e) of \S\ref{valu}:
$\,\th<0\,$ and $\,\px<0\,$ on both beams of $\,\,\mathsf{X}\hs$, while on \
\di\ \ we have $\,\th<0\,$ and $\,\px\,$ is
represented by $\,\pw\,$ with (\ref{vlp}) -- (\ref{vpo}), and, finally,
$\,\th>1\,$ and $\,\px>0\,$ on $\,\,\mathsf{I}\hs$.
\end{rem}
We will now describe some results of \cite{potentials} and use them to prove Theorem~\ref{cnver}. First, according to the discussion following
\cite[Proposition~33.1 in~\S33]{potentials},
\begin{enumerate}
\item Any quadruple $\,\mgmt\,$ satisfying (\ref{zon}) or (\ref{zto}) must belong to one of four disjoint {\it types\/} (a), (b), (c1),
(c2) defined in \S33 of \cite{potentials}.
\item Type (b) cannot occur, as it contradicts the compactness assumption
made in (\ref{zon}) and (\ref{zto}). (In \S33 of \cite{potentials} compactness of
$\,M\,$ is not assumed.)
\item Quadruples of type (a) (or, (c1)) all arise from the construction
in \S\ref{redu}\ (or, \S\ref{main}) of this paper.
\end{enumerate}
Namely, (ii) -- (iii) follow from \cite[Theorems~33.2,~33.3,~34.3 and
Remark~2.4]{potentials}.
\begin{lem}\label{asttt}One of the following two assertions holds for any
given quadruple $\,\mgmt\,$ with {\rm(\ref{zon})} or {\rm(\ref{zto})}.
\begin{enumerate}
\item[($*$)] Up to a $\,\vp$-preserving biholomorphic isometry, $\,\mgmt\,$
arises from the construction in\/ {\rm\S\ref{redu}}, or from that in\/
{\rm\S\ref{main}}.
\item[($*\hskip-.3pt*$)] Conditions {\rm(\ref{qdq})} in {\rm\S\ref{modu}} are
satisfied by some nontrivial closed interval\/ $\,I$ with $\,1\in I$ and some
$\,Q\in\jv\nh$, with\/ $\,\jv\,$ as in {\rm(\ref{spv})} for this given
$\,m\hs$.
\end{enumerate}
\end{lem}
In fact, by (i) -- (iii) above, type (b) is excluded, types (a) and (c1) lead
to ($*$), while, for type (c2), Corollary~35.1 in \cite{potentials} yields
($*\hskip-.3pt*$).

\begin{proof}[Proof of Theorem~\ref{cnver}]Case ($*\hskip-.3pt*$) in Lemma~\ref{asttt} is made impossible by Proposition~\ref{noone}, while the construction in
\S\ref{main} amounts to that described before Proposition~\ref{quadr}: for
$\,I=[\th,\ha]\,$ used in \S\ref{main}, $\,\px\hs$-rationality of $\,(\th,\ha)\,$ is
obvious from Remark~\ref{prati} .
\end{proof}

\section{Bounds on $\,\tz\,$ and $\,\xj$}\label{bnds}
\setcounter{equation}{0}
Let $\,R_m(\ps)=-\hs G_m(\ps)/\vs_m(0)$, with $\,F_m,E_m,\vs_m,G_m$ standing
for $\,F,E,\vs,G$, as in (\ref{fmt}). Also, let $\,\vs_0(0)=-\hs1/2\,$ and
$\,R_0(\ps)=1$. Then, for any integer $\,m\ge1$,
\begin{equation}
\begin{array}{rl}\arraycolsep9pt
\mathrm{i)}&
\vs_m(0)\,=\,(4-6/m)\hs\vs_{m-1}(0)\,,\hskip18pt\vs_1(0)\,=\,1\,,\\
\mathrm{ii)}&
R_m(\ps)=1-m\sa\hs R_{m-1}(\ps)/(4m-6)\hs,\hskip12ptR_1(\ps)
=1+\sa/2\hs,\hskip5pt\mathrm{with}\hskip5pt\sa=\ps^2/(1-\ps)\hs,\\
\mathrm{iii)}&
(R_1(-\hs2),\hs R_2(-\hs2),\hs R_3(-\hs2),\hs R_4(-\hs2))
=(5/3\hs,\hs-11/9\hs,\hs49/27\hs,\hs13/405),\\
\mathrm{iv)}&
\vs_m(0)R_m(\ps)\,=\,\vs_m(0)\,-\,\vs_{m-1}(0)\hs\sa\,+\,\ldots\,
+\,(-1)^m\vs_0(0)\hs\sa^m\hskip-1pt,\end{array}
\label{smo}
\end{equation}
whenever $\,\ps\in(-\infty,0)$, where 
$\,\sa=\ps^2/(1-\ps)\,$ in (iv) as well. Namely, (\ref{sto}.i) gives (i), and (ii)
follows since (\ref{fmt}.ii) clearly remains valid, for $\,\ps\ne1$, even if one
replaces $\,E_m$ by $\,G_m=E_m-F_m/2$. Finally, (i) and (ii) easily imply (iii)
and (iv).

Any given $\,\sa>0\,$ corresponds as in (\ref{smo}) to a unique $\,\ps<0$. In
fact, since $\,\sa=-\hs\ps-1+1/(1-\ps)$, we have $\,d\sa/d\ps<0$, and taking
the limits of $\,\sa\,$ we see that
\begin{equation}
(-\infty,0)\ni\ps\,\mapsto\,\sa=\ps^2/(1-\ps)\in(0,\infty)\hskip14pt
\mathrm{is\ a\ decreasing\ diffeomorphism.}\label{tts}
\end{equation}
With $\,\vs_m(0)\,$ again denoting the sequence given by (\ref{smo}.i) (or (\ref{sto}.i)),
\begin{equation}
\begin{array}{rl}\arraycolsep9pt
\mathrm{a)}&
a_m\hs=\,m^{3/2}4^{1-m}\hs\vs_m(0)\hskip8pt\mathrm{is\ a\ positive\
decreasing\ function\ of}\hskip6ptm\ge1\hs,\\
\mathrm{b)}&
a_m\,\to\,1/\sqrt{\pi\,}\hskip7pt\mathrm{as}\hskip7ptm\to\infty\hs,\\
\mathrm{c)}&
\vs_m(0)\,=\,(4-6/2)(4-6/3)\hs\ldots\hs(4-6/m)
\hskip9pt\mathrm{for\ any\ integer}\hskip7ptm\ge2\hs.\end{array}
\label{ame}
\end{equation}
Namely, (\ref{smo}.i) yields both (c) and
$\,4(a_m/a_{m-1})^2=m(2m-3)^2/(m-1)^3<4\,$ for $\,m\ge2$, which gives (a),
and Wallis's formula
$\,\lim_{\hskip1.2ptm\to\infty}\hskip.4pt2\hs m^{1/2}
\prod_{j=1}^{m-1}[2j/(2j+1)]=\sqrt{\pi\,}\,$ implies (b) since
$\,1/\prod_{j=1}^{m-1}[2j/(2j+1)]
=2^{1-m}[(m-1)\hs!\hs]^{-1}\prod_{j=1}^{m-1}(2j+1)
=4^{1-m}[(m-1)\hs!\hs]^{-2}(2m-1)\hs!$, which, by (\ref{sto}.i), equals
$\,4^{1-m}(2m-1)m\hs\vs_m(0)=(2m-1)m^{-1/2}a_m$.

For a real variable $\,\rz\,$ and for $\,\vs_m(0)\,$ as in (\ref{smo}.i) with
$\,\vs_0(0)=-\hs1/2$, let
\begin{equation}
\varPsi_m(\rz)=\vs_0(0)-\nh\vs_1(0)\hs\rz+
\ldots+(-1)^m\vs_m(0)\hs\rz^m\hskip-1.5pt,\hskip8pt
\varPsi_\infty(\rz)=-\,{\sqrt{1+4\rz\,}\over2}\,,\label{psm}
\end{equation}
with $\,\rz\ge-1/4\,$ in $\,\varPsi_\infty(\rz)$. By (\ref{smo}.i),
$\,\varPsi_m$ satisfies the initial value problem
\begin{equation}
(2\rz+1/2)\hs d\hs\varPsi_m/d\rz\,=\,\varPsi_m(\rz)\,
+\,(2m-1)(-1)^m\hs\vs_m(0)\hs\rz^m\hskip-1.5pt,
\hskip7pt\varPsi(0)=-\hs1/2\hs,\label{trp}
\end{equation}
while, by (\ref{smo}.ii), the power series whose partial sum appears in (\ref{psm}) has the convergence radius $\,1/4$. The sum of this series is
$\,\varPsi_\infty(\rz)$, as one sees either noting that the sum satisfies on
$\,(-1/4,1/4)\,$ the initial value problem
$\,(1+4\rz)\hs d\hs\varPsi_\infty/d\rz=2\varPsi_\infty(\rz)\,$ with
$\,\varPsi(0)=-\hs1/2\,$ (which one may derive from (\ref{trp}) and (\ref{smo}.i), as well as directly from (\ref{smo}.i)), or using (\ref{smo}.iv) to
verify that $\,(-1)^m\vs_m(0)\,$ equals, for every $\,m\ge0$, the $\,m$th
Taylor coefficient of $\,\varPsi_\infty(\rz)\,$ at $\,\rz=0$. Also, by (\ref{smo}.i),
\begin{equation}
\varPsi_m\,\to\,\varPsi_\infty\hskip8pt\mathrm{as}\hskip8ptm\to\infty\hs,
\hskip5pt\mathrm{uniformly\ on}\hskip7pt[-\hs1/4,1/4\hs]\hs,\label{mti}
\end{equation}
since $\,4\hs|\vs_m(0)\hs\rz^m|\le m^{-3/2}a_m$ whenever $\,|\rz|\le1/4$.

If $\,m\,$ is odd and $\,m\ge3$, since
$\,\vs_{2k-1}(0)=k\hs\vs_{2k}(0)/(4k-3)\,$ by (\ref{smo}.i), we get
\begin{equation}
\varPsi_{m-1}(\rz)=-\hs{1\over2}\hs
+\hskip-4.5pt\sum_{k=1}^{(m-1)/2}\hskip-4.5pt\varXi_k(\rz)\hs,
\hskip4pt\mathrm{where}\hskip3pt\varXi_k(\rz)=\vs_{2k}(0)\hskip1pt\rz^{2k-1}
\hskip-1.2pt\left[\rz-\hs{k\over4k-3}\right]\hskip-1.2pt,\label{pmo}
\end{equation}
by grouping terms in (\ref{psm}). Obviously,
$\,\varXi_k\le0\,$ on the interval $\,[\hs0,k/(4k-3)]\,$ and, since these
intervals form a descending sequence, $\,\varXi_k\le0\,$ on
$\,[\hs0,(m-1)/(4m-10)]\,$ for all $\,k=1,\dots,(m-1)/2$. Thus, by (\ref{pmo}),
\begin{equation}
\varPsi_{m-1}(\rz)\le-\hs1/2\hskip5pt\mathrm{whenever}\hskip4ptm\ge3
\hskip4.5pt\mathrm{is\ odd\ and}\hskip4pt0\le\rz\le(m-1)/(4m-10)\hs.
\label{psl}
\end{equation}
\begin{rem}\label{kmsig}
Let $\,R_m(\ps)\,$ and $\,\hi_m$ be defined by (\ref{smo}) and
$\,\hi_m=\xj^2/(1-\xj)$, with $\,\xj<0\,$ as in Lemma~\ref{xxast}(b) for any
$\,\ps\in(-\infty,0)\,$ and an odd integer $\,m\ge3$. Setting
$\,\sa=\ps^2/(1-\ps)$, we have
\begin{enumerate}
\item[a)] if $\,R_m(\ps)\le1$, then $\,\hi_m>\sa$,
\item[b)] $\sgn\hs(\hi_m-\sa)=\hs\sgn\hs[R_{m-1}(\ps)+K_m(\sa)]$, where
$\,K_m(\sa)=2(2m-3)/[(m-1)\sa+2(2m-1)]$.
\end{enumerate}
In fact,
$\,\xo(\ps)/[\vs(0)\hs(1-\ps)]=[(m-1)\sa+2(2m-1)]\hs[1-R_m(\ps)]\,+\,m\sa\,$
as $\,\ps^2=(1-\ps)\sa$, with $\,\xo\,$ as in (\ref{tuu}.b). Since
$\,m,\sa,1-\ps\,$ and $\,\vs(0)\,$ are positive (cf.\ (\ref{esg})), this,
combined with Lemma~\ref{xxast}(c) and (\ref{tts}), gives (a) and, due to (\ref{smo}.ii), also (b).
\end{rem}
\begin{lem}\label{bonds}Given an odd integer\/ $\,m\ge3$, let\/
$\,\zh_m=\tz^2/(1-\tz)\,$ and\/ $\,\hi_m=\xj^2/(1-\xj)\,$ for
$\,\tz,\xj\in(-\infty,0)\,$ defined in {\rm(\ref{god})} and\/ {\rm
Lemma~\ref{xxast}(b)}. Then
\begin{enumerate}
\item $0<4-34/(m^{1/4}+\hs8)<\zh_m\nh<\hi_m\nh\le4-6/(m-1)<4$, with all
inequalities strict if $\,m\ge5$,
\item $\zh_m$ and\/ $\,\hi_m$ are strictly increasing functions of\/
$\,m$,
\item $\hi_m\,\to\,4\,$ and\/ $\,\zh_m\,\to\,4\,$ as $\,m\to\infty$.
\end{enumerate}
\end{lem}
\begin{proof}Let $\,m\ge3\,$ be odd. For $\,\rz\in(0,\infty)\,$ and
$\,\varPsi_{m-1},\hs\sgn\hs\,$ as in (\ref{psm}) and \S\ref{rati},
\begin{equation}
\begin{array}{rl}\arraycolsep9pt
\mathrm{a)}&\sgn\hs(\rz-1/\hi_m)\,=\,\hs\sgn\,\dl_m(\rz)\hs,
\hskip13pt\mathrm{b)}\hskip7pt\sgn\hs(\rz-1/\zh_m)\,
=\,\hs\sgn\,\varPsi_{m-1}(\rz)\hs,\\
&\mathrm{where}\hskip14pt\dl_m(\rz)\,
=\,\varPsi_{m-1}(\rz)\,+\,m\hs\vs_m(0)\hs\rz^m\nh/\hs[2(2m-1)\rz+m-1].
\end{array}
\label{sgr}
\end{equation}
In fact, $\,\dl_m(\rz)\,
=\,\sa^{1-m}\hs\vs_{m-1}(0)\hs[R_{m-1}(\ps)+K_m(\sa)]\,$ for $\,\sa=1/\rz$, by
(\ref{smo}.iv) (with $\,m-1\,$ rather than $\,m$) and (\ref{smo}.i), and so
a) is clear from Remark~\ref{kmsig}(b), Lemma~\ref{xxast}(c) and (\ref{tts}).
Similarly, (b) follows as $\,\hs\sgn\hs(\ps-\tz)\,$ for any
$\,\ps<0\,$ equals $\,\hs\sgn\hs(\zh_m-\sa)=\,\sgn\hs(\rz-1/\zh_m)\,$, where
$\,\rz=1/\sa\,$ with $\,\sa\,$ as in (\ref{tts}); on the other
hand, by (\ref{god}), $\,\hs\sgn\hs(\ps-\tz)\,$ also coincides with
$\,\sgn\hs[G(\ps)-E(0)]$, that is, $\,\sgn\hs[G(\ps)+\vs_m(0)]\,$ (cf.\ (\ref{sto}.ii)), and hence equals $\,\sgn\hs[\vs_m(0)-\vs_m(0)R_m(\ps)]\,$ (see
the definition of $\,R_m$, preceding (\ref{smo})), while
$\,\vs_m(0)-\vs_m(0)R_m(\ps)=\sa^m\varPsi_{m-1}(1/\sa)\,$ by (\ref{smo}.iv) and (\ref{psm}).

Next, $\,\zh_m>0\,$ as $\,\tz<0$, while $\,\zh_m<\hi_m$ since $\,\xj<\tz\,$
(Lemma~\ref{xxast}(b)) or, equivalently, by (\ref{sgr}) for $\,\rz=1/\zh_m$ (as
$\,\dl_m(\rz)>\varPsi_{m-1}(\rz)$). The lower bound on $\,\zh_m$ in (i) is
obtained by finding (as described below) a constant $\,\alpha>0\,$ and a
positive increasing function $\,\k_{\nh m}$ of the odd integer $\,m\ge3\,$
such that $\,\alpha\hs\varPsi_{m-1}(1/4)\,$ and
$\,\alpha\hs\varPsi\hs'_{m-1}(1/4)\,$ are both greater than $\,-\hs1$, while
$\,\alpha\hs\varPsi\hs''_{m-1}(\rz)>2\k_{\nh m}$ for every odd $\,m\ge3\,$
and every $\,\rz\in[\hs1/4,\infty)$, with $\,(\,\,)'=\,d/d\rz$. Then, clearly,
$\,\alpha\hs\varPsi_{m-1}(\rz)>f_m(\rz)\,$ at every $\,\rz\ge1/4$, where
$\,f_m(r)=\k_{\nh m}\hs(\rz-1/4)^2-(\rz-1/4)-1\,$ is the quadratic function of
$\,\rz\,$ whose value, derivative and second derivative at $\,\rz=1/4\,$ are
$\,-\hs1,-\hs1\,$ and $\,2\k_{\nh m}$. By (\ref{sgr}.b), this gives
$\,1/\rz<\zh_m$ for any $\,\rz\ge1/4\,$ at which $\,f_m(\rz)>0$. Next,
$\,f_m(\rz)>0\,$ for $\,\rz=(\beta+1)\k_{\nh m}^{\hs-1/2}+1/4$, where
$\,\beta=\k_{\nh3}^{\hs-1/2}\hskip-1.5pt$. Namely,
$\,f_m(\rz)=\beta^{\hs2}+2\beta-(\beta+1)\k_{\nh m}^{\hs-1/2}\ge\beta\,$ as
$\,\k_{\nh m}\ge \k_{\nh3}$, that is, $\,\k_{\nh m}^{\hs-1/2}\le\beta$. Thus,
$\,\zh_m>\,1/\rz\,=\,4\,-\,16(\beta+1)/[\k_{\nh m}^{\hs1/2}+4(\beta+1)]$. If,
in addition, $\,\k_{\nh m}=16\hs\sqrt{m\,}\hs/\gamma^2$ with a constant
$\,\gamma>0\,$ independent of $\,m$, this gives
$\,\zh_m>\,4\,-\,4(\beta+1)\gamma/[m^{1/4}+\hs(\beta+1)\gamma]$.

A choice of $\,\alpha\,$ and $\,\k_{\nh m}$ with the required properties is
$\,\alpha=\sqrt2\,$ and $\,\k_{\nh m}=\sqrt{\hs2m/\pi\,}$, that is,
$\,\k_{\nh m}=16\hs\sqrt{m\,}\hs/\gamma^2$ with $\,\gamma=(128\hs\pi)^{1/4}$.
In fact, the sequences $\,\varPsi_m(1/4)\,$ and $\,\varPsi\hs'_m(1/4)$, with
$\,m=1,2,3,\dots\,$, converge to $\,-\hs1/\sqrt2$, as one sees setting
$\,\rz=1/4\,$ in (\ref{mti}) and, respectively, (\ref{trp}) (where the term
involving $\,\vs_m(0)\,$ tends to
$\,0\,$ by (\ref{ame}.a)). On the other hand, $\,\varPsi_{m-1}(1/4)\,$ and
$\,\varPsi\hs'_{m-1}(1/4)\,$ are decreasing functions of the {\it odd\/}
integer $\,m\ge3$, since $\,\varXi_k(\rz)\,$ and $\,\varXi\hs'_k(\rz)\,$ are,
by (\ref{pmo}), negative for $\,\rz=1/4\,$ and any $\,k\ge1$. Also,
$\,\varXi\hs''_k(\rz)\,$ is a nondecreasing function of
$\,\rz\in[\hs1/4,\infty)$, as $\,\varXi\hs''_k(\rz)
=2k(2k-1)\vs_{2k}(0)\hskip1pt\rz^{2k-3}\hs[\rz-(k-1)/(4k-3)]$. (Note that
$\,\rz-(k-1)/(4k-3)\,$ is a {\it positive} increasing function of
$\,\rz\ge1/4$.) By (\ref{pmo}), the same is true of $\,\varPsi\hs''_{m-1}(\rz)$.
Thus, we just need to establish the inequality
$\,\sqrt{2\,}\hs\varPsi\hs''_{m-1}(1/4)>2\k_{\nh m}$ for
$\,\k_{\nh m}=\sqrt{\hs2m/\pi\,}$. Now (\ref{ame}) and the above formula
for $\,\varXi\hs''_k(\rz)\,$ give
$\,\varXi\hs''_k(1/4)=2a_{2k-1}/\sqrt{2k-1\,}>2/\sqrt{(2k-1)\pi\,}\,$ for
all $\,k\ge1$, and so, by (\ref{pmo}),
$\,\varPsi\hs''_{m-1}(1/4)>\sqrt{\hs4m/\pi\,}$, as required. (In fact, easy
induction gives
$\,1^{-1/2}+3^{-1/2}+5^{-1/2}+\ldots+(m-2)^{-1/2}\ge\sqrt{m\,}\,$ for any odd
integer $\,m\ge3$, as $\,(m^{1/2}+m^{-1/2})^2>m+2$.) Now
$\,(\beta+1)\gamma=(8\hs\pi/\sqrt3\hs)^{1/2}+\hs(128\hs\pi)^{1/4}\,
\approx\,8.29$, for $\,\beta=\k_{\nh3}^{\hs-1/2}\hskip-1.5pt$, so that
$\,(\beta+1)\gamma>8\,$ and $\,4(\beta+1)\gamma<34$. The inequality concluding
the last paragraph thus gives the lower bound for $\,\zh_m$ appearing in (i).

To prove the remainder of (i) we may
assume that $\,m\ge5$, since Example~\ref{xymtr} gives $\,\xj=-\hs(\sqrt5+1)/2\,$
for $\,m=3\,$ and hence $\,\hi_3=1=4-6/2$. The inequality
$\,\hi_m<4-6/(m-1)\,$ will be obvious from (\ref{sgr}.a) once we show that
\begin{equation}
m\hs\vs_m(0)\hs\rz^m\nh/\hs(m-1)\,<\,1\hskip8pt\mathrm{for}\hskip6pt
m\ge5\hskip6pt\mathrm{and}\hskip6pt\rz=(m-1)/(4m-10)\hs,\label{smo0}
\end{equation}
since the definition of $\,\dl_m$ in (\ref{sgr}) combined with (\ref{psl}) then
will give $\,\dl_m(\rz)<0$. Note that $\,1/\hs[2(2m-1)\rz+m-1]<1/[2(m-1)]$.

We derive (\ref{smo0}) from the fact that
$\,(4\rz)^{m-5/2}=(1+1/\nu)^{3\nu/2}<e^{3/2}<5$, for $\,m,\rz\,$ as in (\ref{smo0}), with $\,\nu=(2m-5)/3$. (That $\,\nu\hs\log\hs(1+1/\nu)$, and hence
also $\,(1+1/\nu)^\nu$, is an increasing function of $\,\nu>1$, is clear as
$\,\mu^{-1}\log\hs(1+\mu)\,$ is a decreasing function of
$\,\mu=1/\nu\in(0,1)$, which follows since $\,(1+\mu)\mu^2$ times its
derivative decreases on $\,[\hs0,1\hs]\,$ from $\,0\,$ to $\,\log\hs(e/4)$.)
Thus, $\,(4\rz)^m<5(4\rz)^{5/2}<20$, as $\,\rz\le2/5$. Also, (\ref{ame}) for
$\,m\ge5\,$ gives $\,a_m<a_2=\sqrt{2\,}/8<1/4\,$ and
$\,\vs_m(0)<m^{-3/2}4^{m-1}$, so that $\,m\hs\vs_m(0)\hs\rz^m\nh/\hs(m-1)<
m^{-1/2}(4\rz)^m\nh/\hs[4(m-1)]$, which is less than
$\,5m^{-1/2}\nh/\hs(m-1)<1$. This yields (i) and, consequently, (iii).

To prove (ii) for $\,\zh_m$, let us fix $\,m\ge5\,$ and set
$\,\rz=1/\zh_m$. By (\ref{pmo}) and (\ref{sgr}.b),
$\,-\hs\varPsi_{m-3}(\rz)=\varPsi_{m-1}(\rz)-\varPsi_{m-3}(\rz)
=\vs_{m-1}(0)\hs\rz^{m-2}[\hs\rz-(m-1)/(4m-10)]$, and so
$\,\varPsi_{m-3}(\rz)<0\,$ since (i) gives $\,\zh_m<\,4-6/(m-1)$, that is,
$\,\rz>(m-1)/(4m-10)$. Now (\ref{sgr}.b), for $\,m-2\,$ rather than $\,m$, gives
$\,\zh_m=1/\rz>\zh_{m-2}$.

Finally, to obtain (ii) for $\,\hi_m$, note that, for any $\,\rz>0\,$ and
$\,m\ge5$,
\begin{equation}
{\dl_m(\rz)\,-\,\dl_{m-2}(\rz)\over(2m-5)(m-1)\hs\vs_{m-1}(0)}\,
=\,{\rz^{m-2}(1+4\rz)^2\hs[\rz\,-\,(m-1)/(4m-10)]
\over[(4m-2)\rz+m-1]\hs[(4m-10)\rz+m-3]}\,,
\end{equation}
since (\ref{pmo}) leads to a massive cancellation of terms in
$\,\dl_m(\rz)-\hs\dl_{m-2}(\rz)\,$ and, by (\ref{smo}.i),
$\,\vs_{m-2}(0)=(m-1)\hs\vs_{m-1}(0)/(4m-10)$,
$\,\hs\vs_m(0)=(4-6/m)\hs\vs_{m-1}(0)$. Applied to
$\,\rz=1/\hi_m$, this gives, by (\ref{sgr}.a),
$\,-\hs\dl_{m-2}(\rz)=\dl_m(\rz)-\hs\dl_{m-2}(\rz)>0$, as the inequality
$\,\hi_m<4-6/(m-1)\,$ in (i) amounts to $\,\rz>(m-1)/(4m-10)$. Thus,
$\,\dl_{m-2}(\rz)<0\,$ and (\ref{sgr}.a) with $\,m\,$ replaced by $\,m-2\,$
yields $\,\rz<1/\hi_{m-2}$, that is, $\,\hi_m>\hi_{m-2}$, completing the proof.
\end{proof}
\begin{rem}\label{xinft}
Let $\,\xj_\infty=-\hs2(1+\sqrt2)\hs\approx-\hs4.83$.
Combined with (\ref{tts}), Lemma~\ref{bonds} immediately implies that both $\,\tz\,$
and $\,\xj\,$ are strictly decreasing functions of the odd integer $\,m\ge3$,
such that $\,\tz>\xj>\xj_\infty$, while $\xj\to\xj_\infty\,$ and
$\,\tz\to\xj_\infty$ as $\,m\to\infty$.
\end{rem}

\section{Further inequalities}\label{ftiq}
\setcounter{equation}{0}
Let $\,\qx_*\in(-\hs1,0)\,$ depend on an odd
integer $\,m\ge3\,$ as in (\ref{qvu}), and let $\,\lk_\infty=9+6\sqrt2$.
Then, with $\,\lkm\in(1,\infty)\,$ such that $\,\qx_*=-\hs1+1/\lkm$,
\begin{enumerate}
\item[a)] $\lkm\,<\hs18\,$ for every odd $\,m\ge3$. More precisely,
$\,\lkm<\lk_\infty\hs\approx\hs17.485$.
\item[b)] $\lkm\,<\,m\,$ whenever $\,m\,$ is odd and $\,m\ge9$.
\item[c)] $\lkm\,>\,m\,$ if $\,m\in\{3,5,7\}$, while $\,8<\lk_9<9<\lk_{11}<10$.
\item[d)] $\lk_3=(9\,+\,3\sqrt{5})/4\hs\approx\hs3.927$,
\item[e)] $\lk_m$ is a strictly increasing function of $\,m\,$ and
$\,\lk_m\to\lk_\infty$ as $\,m\to\infty$.
\end{enumerate}
In fact, by (\ref{qvu}), $\,\lkm=\xi\hs(\xi+\sqrt{\xi^2-1}\hs)\,$ for
$\,\xi=\xi_m$ given by $\,\xi=1+\hi/2\,$ with $\,\hi=\hi_m$ as in Lemma~\ref{bonds}. Thus, if $\,m=3$, Example~\ref{xymtr} gives $\,\xj=-\hs(\sqrt 5+1)/2$, so
that $\,\hi_3=1\,$ and $\,\xi_3=3/2$, which yields (d), and hence (c) for
$\,m=3$.

Obviously, $\,\lkm$ is an increasing function of $\,\xi_m>1$, and hence of
$\,\hi_m>0$. This has several consequences. First: (e), (a) are obvious from
Lemma~\ref{bonds}, and (a) gives (b) for $\,m\ge19$. Next, $\,R_4(-\hs2)>0\,$
in (\ref{smo}.iii), and so $\,R_5(-\hs2)<1\,$ by (\ref{smo}.ii); hence
$\,\hi_5>4/3$, from Remark~\ref{kmsig}(a) for $\,\ps=-\hs2$,
so that $\,\xi_5>5/3$, and (c) for $\,m=5\,$ follows. Finally, (b) and (c) for
$\,m=7,9,\dots,17\,$ are, similarly, numerical consequences of the bounds on
$\,\hi_m$ provided by (\ref{chi}) below.
\begin{rem}\label{bdslm}
Lemma~\ref{bonds}(i) clearly implies some explicit, though complicated bounds
for $\,\lkm$. Replacing them with weaker but simpler estimates, we get
$\,1-17/(m^{1/4}+\hs8)<\lkm/\lk_\infty<1-1/m$.
(The lower bound is of interest only when it is positive, that is, for very
large $\,m$.) Namely, as $\,0<\hi<4-6/m\,$ for $\,\hi=\hi_m$, setting
$\,\xi=1+\hi/2\,$ we get $\,\xi<3-3/m$, and so
$\,\sqrt{\xi^2-1}<2\sqrt{2\,}\xi/3\,$ as $\,1<\xi<3$. Thus,
$\,\lkm<\lk_\infty\xi^2/9<\lk_\infty\xi/3<(1-1/m)\lk_\infty$, since
$\,\lkm=\xi\hs(\xi+\sqrt{\xi^2-1}\hs)\,$ (see above). Similarly, for
$\,L(\eta)=\eta\hs(\eta+\sqrt{\eta^2-1}\hs)-\lk_\infty(\eta-2)$, with
$\,\eta\in(1,\infty)$, we have $\,L(3)=0\,$ and $\,dL/d\eta<0\,$ at
$\,\eta=3$, while $\,\dsq L/d\eta^2>0\,$ whenever $\,\eta>2/\sqrt{3\,}$, as
$\,\dsq L/d\eta^2=(2-\phi)(1+\phi)^2$ with $\,\phi=\eta/\sqrt{\eta^2-1\hs}$,
so that $\,d\phi/d\eta<0\,$ and condition $\,\eta>2/\sqrt{3\,}\,$ is
equivalent to $\,\phi<2$. Therefore, $\,dL/d\eta<0$, and hence $\,L(\eta)>0\,$
for all $\,\eta\in[\hs2/\sqrt{3\,},\hs3)$, and, with $\,\eta=3-17/(m^{1/4}+\hs8)$,
Lemma~\ref{bonds}(i) now leads to our lower bound.
\end{rem}

\section{Some simple facts from number theory}\label{nmth}
\setcounter{equation}{0}
The following lemma is a variation on the $\,s=2\,$ case of the well-known
fact that, for any integer $\,s\ge2$, the probability that $\,s\,$ randomly
chosen positive integers have a common divisor other than $\,1\,$ equals
$\,1/\zeta(s)$, where $\,\zeta\,$ is the Riemann zeta function.

We allow $\,\lkm$ and $\,\lk_\infty$ to be much more general here than in
\S\ref{ftiq} or \S\ref{frth}.
\begin{lem}\label{mrtns}Given a real constant\/ $\,\lk\ge1\,$ and an integer
$\,m\ge1$, let\/ $\,\p(\lk,m)\,$ be the set of all pairs $\,(k,\el\hs)\,$ of
relatively prime integers with $\,1\le k\lk<\el\le m\hs$. Then
\begin{enumerate}
\item[a)] $\pi^2\hs|\p(\lk,m)|\hs/m^2\,\to\,3\hs/\lk\,\,$ as $\,\,m\to\infty$,
with $\,|\hskip3pt|\,$ denoting cardinality,
\item[b)] $\pi^2\hs|\p(\lkm,m)|\hs/m^2\,\to\,3\hs/\lk_\infty$ as
$\,m\to\infty$, whenever $\,\lk_\infty\in\bbR\,$ is the limit, as
$\,m\to\infty$, of a function $\,m\mapsto\lkm\in[\hs1,\infty)\,$ defined on an
infinite set of positive integers $\,m\hs$.
\end{enumerate}
\end{lem}
In fact, one obtains (a) by modifying a standard proof (see \cite{chandrasekharan}) of
Mer\-tens's theorem. For the reader's convenience, we provide details in an
appendix (\S\ref{mert}).

Next, $\,\p(\lk\hs'\hskip-1.2pt,m)\subset\p(\lk,m)\,$ if
$\,\lk\hs'\ge\lk\ge1$, and so, if $\,\lk\hs'\ge\lkm\ge\lk\ge1\,$ for $\,\lkm$
as in (b), then $\,|\p(\lk\hs',m)|/m^2\le|\p(\lkm,m)|/m^2\le|\p(\lk,m)|/m^2$.
This obviously is the case for any fixed $\,\lk,\lk\hs'$ with either
$\,\lk\hs'>\lk_\infty>\lk\ge1\,$ (when $\,\lk_\infty>1$), or
$\,\lk\hs'>\lk_\infty=\lk=1$, and sufficiently large $\,m\,$ for which
$\,\lkm$ is defined. Taking the upper/lower limits, we now get (b) from (a)
for $\,\lk,\lk\hs'$ arbitrarily close to $\,\lk_\infty$.
\begin{rem}\label{farey}
For $\,\p(\lk,m)\,$ as above, $\,(k,\el\hs)\mapsto k/\el\,$ is a bijective
correspondence between
$\,\p(\lk,m)\,$ and the set of all those rational numbers in
$\,(0,1/\lk\hs)\,$ which can be written as fractions with a denominator in
$\,\{1,\dots,m\}$. When $\,\lk=1$, the elements of the latter set along with
the numbers $\,0\,$ and $\,1$, listed in increasing order, form what is called
the {\it Farey sequence of order} $\,m\,$ (cf.\ \cite{chandrasekharan}).

It is clear that the Farey sequence of order $\,m\,$ has
$\,2+\varphi(2)+\ldots+\varphi(m)\,$ elements, where $\,\varphi\,$ is the {\it
Euler function}, assigning to a positive integer $\,k\,$ the number of
integers $\,j\,$ such that $\,0<j<k\,$ and $\,j,k\,$ are relatively prime.
\end{rem}

\section{Examples with locally reducible metrics}\label{redu}
\setcounter{equation}{0}
This and the next three sections describe constructions of the four families,
mentioned in \S\ref{intr}, of quadruples $\,\mgmt\,$ satisfying (\ref{zon}) or (\ref{zto}).

The {\it first family} is represented by just one $\,\px\hs$-rational point
$\,(0\hs,\nh0)\,$ on the moduli curve $\,\mathcal{C}\hs$ (see the end of
\S\ref{modu}), and
consists of those $\,\mgmt\,$ with (\ref{zon}) or (\ref{zto}) for which $\,g\,$ is
locally reducible as a K\"ahler metric. Such $\,\mgmt\,$ seem to be
well-known, and can be constructed as follows.

Given an integer $\,m\ge2$, real constants $\,\kx<0\,$ and $\,\y\ne0$, and a
compact K\"ahler-Einstein manifold $\,(N,h)\,$ of complex dimension $\,m-1\,$
with the Ricci tensor $\,\hs\rih=\kx\hs h$, let $\,\mathcal{L}\,$ be any
holomorphic line bundle over $\,N\,$ carrying a fixed flat
$\,\hs\mathrm{U}\hs(1)\,$ connection. Next, let $\,\mathcal{E}=N\times\bbR\,$ be
the product real-line bundle over $\,N\hs$ with the obvious flat connection
and Riemannian fibre metric, and let $\,M\,$ be the unit-sphere bundle of the
direct sum $\,\mathcal{L}\oplus\mathcal{E}$. Thus, $\,M\,$ is a $\,2$-sphere bundle over
$\,N\nh$. Since the direct-sum connection in $\,\mathcal{L}\,\oplus\,\mathcal{E}\,$ is
flat and compatible with the direct-sum metric, its horizontal distribution is
both integrable and tangent to the submanifold $\,M$, and so it gives rise to
an integrable distribution which may also be called {\it horizontal\hs}, and
whose leaves, along with the $\,\bbCP^1$ fibres form, locally in $\,M$, the
factor manifolds of a Cartesian-product decomposition.

Let $\,g\,$ now be a metric on $\,M\,$ such that the horizontal distribution
is $\,g$-normal to the fibres and $\,g\,$ restricted to it is the pullback of
$\,h\,$ under the bundle projection $\,M\to N\nh$, while $\,g\,$ on each fibre
equals $\,(3-2m)/\kx\,$ times the standard unit-sphere metric. Thus,
$\,(M,g)\,$ is a K\"ahler manifold since, locally, it is a Riemannian product
with the factors manifolds which are a $\,2$-sphere $\,S^2$ of constant
Gaussian curvature $\,\kx/(3-2m)\,$ and $\,(N,h)$, while the $\,2$-sphere
factors can be coherently oriented, which makes them K\"ahler manifolds of
complex dimension one.

Finally, let $\,\vp:M\to\bbR\,$ be the composite
$\,M\to\mathcal{L}\oplus\hs\mathcal{E}\to\mathcal{E}\to\bbR\to\bbR\,$ of the inclusion
mapping of $\,M$, followed by the direct-sum projection morphism, followed by
the Cartesian-product projection $\,\mathcal{E}=N\times\bbR\to\bbR\hs$, followed by
the multiplication by the nonzero constant $\,\y\,$ in $\,\bbR\hs$. In terms
of a local Riemannian-product decomposition just described, with the $\,S^2$
factor treated as the sphere of radius $\,\sqrt{(3-2m)/\kx\,}\,$ about $\,0\,$
in a Euclidean$\,3$-space $\,V$, our $\,\vp\,$ is a function on $\,M$,
constant in the direction of the $\,N\hs$ factor, and, as a function
$\,\vp:S^2\to\bbR$, it is the restriction to $\,S^2$ of a nonzero linear
homogeneous function $\,V\to\bbR\hs$.

The quadruple $\,\mgmt\,$ then satisfies (\ref{zon}) or (\ref{zto}); see
\cite[\S30]{potentials} or \cite[\S25]{local}.
\begin{rem}\label{irred}
Unlike the examples just described, any quadruple
$\,\mgmt\,$ constructed as in \S\ref{main} is {\it locally irreducible} in
the sense that no open submanifold of $\,(M,g)\,$ is biholomorphically
isometric to a Cartesian product of lower\diml\ K\"ahler manifolds.

In fact, the horizontal and vertical distributions on $\,M\,$ then consist of
eigenvectors of both the Ricci tensor $\,\hs\ri\hs\,$ of $\,g\,$ and the
second covariant derivative $\,\nabla d\vp\,$ of $\,\vp\,$ relative to $\,g$,
with some eigenvalue functions $\,\la,\my\,$ for $\,\hs\ri\hs\,$ and
$\,\si,\ta\,$ for $\,\nabla d\vp$, all of which are also functions of $\,\ps$,
that is, of $\,\vp$. (We borrow these notations from
\cite[formula~(7.4)]{local}, cf.\ \cite[\S8]{local}.) One then has
$\,Q=2(\vp-\y)\hs\si\,$
and $\,(\la-\my)\vp=2(m-1)(\vp-\y)\hskip2ptd\si/d\vp$. Thus, by (\ref{qdq}.a,\hs b,\hs c), $\,\si\,$ is not constant on any nonempty open subset
of $\,M$, and so $\,g$, restricted to any such subset, cannot be Einstein.
Local irreducibility of $\,(M,g)\,$ now follows from
\cite[Corollaries~13.2(iii) and~9.3]{local}.
\end{rem}

\section{B\'erard Bergery's and Page's examples}\label{bbpg}
\setcounter{equation}{0}
Given an integer $\,m\ge2$, let $\,\smi$ be the set of all
$\,\px\hs$-rational points in the $\,\,\mathsf{I}\,\,$ component of the
moduli curve (cf.\ Definition~\ref{pratl} and the lines following (\ref{xii})). By
applying the construction of \S\ref{stat}\ to all points $\,(\th,\ha)\in\smi$ and
appropriate additional data, we obtain the {\it second family} of quadruples
$\,\mgmt\,$ satisfying (\ref{zon}) or (\ref{zto}). Since $\,\th,\ha\,$ then are
both positive, so is the function $\,\ps=\vp/\y\,$ on $\,M\,$ (cf.\ Remark~\ref{detrm}); thus, all compact K\"ahler manifolds $\,(M,g)\,$ obtained here are
{\it globally\/} conformally Einstein.

The second family has been known for over two decades: the (essentially
unique) quadruple with (\ref{zto}) was found by Page \cite{page}, and those
with (\ref{zon}) by B\'erard Bergery \cite{berard-bergery}. More precisely, they both
described the corresponding conformally related Einstein manifolds. (See also
\cite[\S26]{local}.)
\begin{proof}[Proof of Theorem~\ref{asymp}, parts $\mathrm{(a),\hskip2pt(d)}$]Assertion
(d), which also implies finiteness of $\,\smi$, is obvious from (e) in \S\ref{valu} and Definition~\ref{pratl}, while (a) in Theorem~\ref{asymp} follows from (d), Lemma~\ref{mrtns}(a) for $\,\lk=1$, and Remark~\ref{farey}.
\end{proof}
Of particular interest is the case where $\,m=2$, with $\,(M,g)\,$
conformal to the Page manifold. By Theorem~\ref{asymp}(d),
$\,\mathcal{S}_2^{\hs\mathsf{I}}$ has a single element $\,(\th,\ha)\,$ with
$\,\ha=\th/(\th-1)$. Explicitly, $\,\th\,$ then is given by
$\,\th=[(3-\alpha+2\beta)^{1/2}+(3+\alpha)^{1/2}]\hs(\alpha+\beta)/6$, where
$\,\alpha=(13+\sqrt{142\,})^{1/3}+(13-\sqrt{142\,})^{1/3}-1\,$ and
$\,\beta=(6+\alpha^2)^{1/2}$, so that $\,\alpha\,$ is the unique real root of
the equation $\,\alpha^3+3\alpha^2-6\alpha=34$. (In fact, relation
$\,\pw(\th)=1/2$, with $\,\pw\,$ as in (\ref{uup}) for $\,m=3$, amounts to
$\,(\th-2)(3\th^4-8\th^3+6\th^2-2)=0$.) Approximately, $\,\th\approx1.560\,$
and $\,\ha\approx2.786$.

\section{A third family of examples}\label{thrd}
\setcounter{equation}{0}
The {\it third family} of quadruples $\,\mgmt\,$ satisfying (\ref{zon}) or (\ref{zto}) is obtained by applying the construction preceding Proposition~\ref{quadr} to $\,\px\hs$-rational points of the moduli curve $\,\mathcal{C}\hs$ that
lie, for even $\,m$, in $\,\,\di\hs\smallsetminus\{(0\hs,\nh0)\}$, or, for odd
$\,m$, in the union of $\,\,\di\hs\smallsetminus\{(0\hs,\nh0)\}\,$ and the
$\,\mathcal{H}$-beam $\,\xah\,$ of $\,\,\mathsf{X}\,\,$ (see the lines
following (\ref{xii})). The second and third families may be thought of as
related to each other by a form of analytic continuation, since they both use
$\,\px\hs$-rational points $\,(\th,\ha)\,$ lying on the hyperbola
$\,\mathcal{H}\,$ given by $\,\th\ha=\th+\ha\,$ (cf.\ \S\ref{vert}).

In the complex dimension $\,m=2$, quadruples of the third family were first
found by Hwang and Simanca \cite{hwang-simanca}) and T\o nnesen-Friedman \cite{tonnesen-friedman}.

Just as we did for the second family in \S\ref{bbpg}, we will now obtain a rough idea
about the ``size'' of the third family by estimating the number of
$\,\px\hs$-rational points involved. In \S\ref{bbpg} that amounted to proving (a)
and (d) in Theorem~\ref{asymp}. Here, the corresponding results consist of Theorem~\ref{ctbly}, already established in \S\ref{valu}, and assertions (b), (e) in Theorem~\ref{asymp}, proved at the end of this section. (The other parts of Theorem~\ref{asymp} pertaining to the third family then are immediate from (b) in \S\ref{valu}:
$\,\smh$ is finite since a bounded interval contains only finitely many
rational numbers that can be written as fractions with a denominator in
$\,\{1,\dots,m\}$, while $\,\px\,$ maps $\,\smh\,$ into $\,(-\hs2,0)\,$ since
$\,\pw(\th_+^*)>-\hs2$, as shown in (iii) below.)

Let $\,\pw,\th_+,\xj,\xj_\infty$ and $\,\tz\,$ depend on an odd integer
$\,m\ge3\,$ as in (\ref{uup}), (\ref{uwz}), Lemma~\ref{xxast}(b), Remark~\ref{xinft} and
\S\ref{anth}, and let $\,\hatt=\ps/(\ps-1)\,$ if
$\,\ps\in\bbR\smallsetminus\{1\}$. Then
\begin{enumerate}
\item[i)] $\xj_\infty<\th_+^*<\xj$,\hskip6.1pt
ii)\hskip6pt$\th_+^*\to\xj_\infty\,$ as $\,m\to\infty$,\hskip6.1pt
iii)\hskip6pt$\pw(\th_+^*)>-\hs2\hs$,
\hskip4.3ptiv)\hskip6pt$\pw(\th_+^*)\to-\hs\sqrt2\,$ as $\,m\to\infty$,
\item[v)] $\tz/(2-\tz)\,$ is a decreasing function of the odd integer
$\,m\ge3$, such that $\,\tz/(2-\tz)>-\hs1+1/m\,$ and
$\,\tz/(2-\tz)\to-\hs1/\sqrt2\,$ as $\,m\to\infty$,
\item[vi)] $\xj^*<\th_+<\xj_\infty^*<1\,$ and $\,\pw(\th_+)<2$, which
improves on (\ref{uwz}),
\item[vii)] $3/2<\pw_+(2/3)<2$, with $\,\pw_\pm$ as in \S\ref{pfac},\hskip29pt
viii)\hskip7pt$(\xj_\infty\hskip-1pt-2)/\xj_\infty=\,\sqrt2$.
\end{enumerate}
In fact, (v) follows from Remark~\ref{xinft}:
$\,\tz/(2-\tz)=-\hs1+2/(2-\tz)>-\hs1+1/m\,$ since $\,\tz>2(1-m)$, which for
$\,m\ge5\,$ is clear as $\,\tz>\xj>\xj_\infty$ and so
$\,\tz>\xj_\infty>2(1-m)$, and for $\,m=3\,$ is obvious as
$\,\tz\hs\approx\hs-\hs1.3\,$ (see the end of \S\ref{anth}).

Let us now set $\,\th=\xj_\infty$. We obviously have $\,(\th-2)/\th=\sqrt2$,
that is, (viii), and $\,\th-1=-\hs\th^2/4$. Dividing (\ref{uup}) by
$\,-\hs\th^2$ and using the definition of $\,R_m$ preceding (\ref{smo}), we now
get $\,-\hs\sqrt{2\,}\hs\pw(\th)=[1-1/(2m)]\hs[1-1/R_m(\th)]+1$, for any fixed
odd integer $\,m\ge3$. Relation $\,R_m(\th)\ge3m\,$ for any odd $\,m\ge1$,
easily obtained from (\ref{smo}.ii) with $\,\ps=\th=\xj_\infty$ and $\,\sa=4\,$ using induction on
$\,m=1,3,5,\dots\,$, now yields $\,\pw(\th)\le\pw_{\!*}<0$ for $\,\pw_{\!*}$
given by $\,-\hs6\sqrt2\hs m^2\pw_{\!*}=12m^2-5m+1$. On the other hand,
dividing (\ref{udp}) by $\,\th\,$ and using (viii) we obtain
$\,(\th-1)\,\dot\pw=m\pw^2+\sqrt2\hs(m-1)\pw-1$, the right-hand side of
which is easily verified to be positive when $\,\pw\,$ (that is, $\,\pw(\th)$)
is replaced by $\,\pw_{\!*}$ and $\,m\ge3$. As $\,-\hs1<0$, that
right-hand side is a quadratic function of $\,\pw\,$ with roots of opposite
signs; thus, it is strictly decreasing on the subset of the negative $\,\pw\,$
axis on which it is positive. Therefore it is positive at $\,\pw(\th)\,$ as
well, and so $\,(\th-1)\,\dot\pw>0$, at $\,\th=\xj_\infty$, that is,
$\,\dot\pw(\th)<0$. Since we also have $\,\pw(\th)<0$, Remark~\ref{behvp} now
gives $\,\th<\th_+^*$ for $\,\th=\xj_\infty$, and, as $\,\xj^*<\th_+$ (see
(\ref{uwz})), we obtain (i) and (ii) from Remark~\ref{xinft}.

Next, let $\,\pw_\pm$ be as in \S\ref{pfac}, so that $\,\pw_\pm$ depend here on an
odd integer $\,m\ge3$, and let a sequence $\,\th_m<0$, $\,m=3,5,7,\dots\,$,
converge to a limit $\,\th_\infty<0$ as $\,m\to\infty$. The expression for
$\,\pw_-(\th)\,$ with any $\,\th<0$, provided by the quadratic formula, shows
that $\,\pw_-(\th_m)\,$ then has the limit
$\,-\hs(\th_\infty\hskip-1pt-2)/\th_\infty$, equal to $\,-\hs\sqrt2\,$ if
$\,\th_\infty=\xj_\infty$ (see (viii)). Using the sequence $\,\th_m=\th_+^*$
and (ii), we now get (iv).

For an odd integer $\,m\ge3$, let $\,\mathcal{F}(\pw)=m\pw^2-2(m-1)\pw-1$, which is
$\,3/2\,$ times the right-hand side of (\ref{udp}) with $\,\th=2/3$. As
$\,\mathcal{F}(3/2)<0<\mathcal{F}(2)$, (vii) follows, $\,\pw=\pw_+(2/3)\,$ being the
positive root of $\,\mathcal{F}$. Next, for any odd $\,m\ge5\,$ (or, $\,m=3$),
$\,1+\vs(0)/G(\th)\,$ in (\ref{uup}) with $\,\th=-\hs2\,$ is negative (or,
respectively, equal to $\,22/49$). This is clear from the definition of
$\,R_m$, since $\,R_3(-\hs2)=49/27\,$ by (\ref{smo}.iii), from which, using (\ref{smo}.ii) and induction on $\,k$, we obtain
$\,0<R_k(-\hs2)<1\,$ for every integer $\,k\ge4$. Hence, if $\,m\ge5\,$ (or,
$\,m=3$), (\ref{uup}) gives $\,\pw(-\hs2)>-\hs1/2\,$ (or,
$\,\pw(-\hs2)=-\hs52/49$), so that, (d) in
\S\ref{sprp} with $\,\th=-\hs2\,$ yields $\,\pw(2/3)<1/2\,$ (or, respectively,
$\,\pw(2/3)=52/49$). When $\,m>3\,$ this clearly yields $\,\th_+>2/3$, since
$\,2/3\notin[\th_+,1]$, as $\,\pw\ge1\,$ on $\,[\th_+,1]\,$ by (\ref{vpo}); if,
however, $\,m=3$, relations $\,\pw_+(2/3)\in(3/2,2)\,$ (see (vii)) and
$\,0<\pw(2/3)=52/49<3/2\,$ give, by (i) in \S\ref{pfac},
$\,\pw_-(\th)<\pw(\th)<\pw_+(\th)\,$ with $\,\th=2/3$, so that
$\,\dot\pw(2/3)>0\,$ (see (iii) in \S\ref{pfac}), and hence $\,\th_+>2/3$, as
$\,\dot\pw\le0\,$ on $\,[\th_+,1]\,$ by (\ref{vpo}). Thus, $\,\pw(\th_+)<2$,
which we verify in three steps. First, $\,\pw(\th_+)=\pw_+(\th_+)\,$ by
(iv),\hs(i) in \S\ref{pfac}, since $\,\pw(\th_+)\,$ is the maximum of $\,\pw\,$ on
$\,(0,1)$. Secondly, $\,\pw_+(\th_+)<\pw_+(2/3)$, since we just showed that
$\,\th_+>2/3$, and $\,\pw_+$ is decreasing on $\,(0,1)$, cf.\ (\ref{vpp}).
Thirdly, $\,\pw_+(2/3)<2\,$ by (vii). Now (iii) follows from (d) in \S\ref{sprp},
while (d),\hs(a) in \S\ref{sprp} and (iii),\hs(i) give (vi).
\begin{proof}[Proof of Theorem~\ref{asymp}, parts
$\mathrm{(b),\hskip2pt(e)}$]Assertion
(e) is immediate from (b) in \S\ref{valu}, since (iii) and (v) above give
$\,-\hs2<\pw(\th_+^*)<-\hs1\,$ and $\,\tz/(2-\tz)>-\hs1+1/m$. To prove (b),
recall (\S\ref{stat}) that the $\,\mathcal{H}$-beam of $\,\,\mathsf{X}\,\,$ is the
graph of the function $\,\ha=\hu$ on the interval $\,(-\infty,\tz)\,$ of the
variable $\,\th$. Dividing $\,(-\infty,\tz)\,$ into the three subintervals
$\,(-\infty,\th_+^*]$, $\,[\th_+^*,\th']$,
$\,[\th',\tz]$, for the unique $\,\th'<0\,$ with $\,\pw(\th')=-\hs1\,$ (cf.\
Remark~\ref{behvp}), we also divide the $\,\mathcal{H}$-beam into three segments
(subbeams). Since $\,\pw\,$ is the restriction of the function $\,\px\,$ to
the hyperbola $\,\ha=\hu$ (see \S\ref{trfp}), Remark~\ref{behvp} shows that $\,\px\,$
maps the set of all $\,\px\hs$-rational points in the first (or, second, or,
third) subbeam bijectively onto the set of all rational numbers that have
positive denominators not exceeding $\,m\,$ and lie in an interval with the
endpoints $\,\pw(\th_+^*)\,$ and $\,-\hs1\,$ (or, again, $\,\pw(\th_+^*)\,$
and $\,-\hs1$, or, respectively, $\,-\hs1\,$ and $\,\tz/(2-\tz)$). Using,
instead of $\,\px$, the function $\,1+1/\px\,$ for the first two subbeams,
and $\,1+\px\,$ for the third, we obtain an analogous property for new
intervals, with the lower endpoint $\,0\,$ and the upper endpoint
$\,1+1/\pw(\th_+^*)\,$ for the first two, or $\,1+\tz/(2-\tz)\,$ for the
third subbeam. Now (b) in Theorem~\ref{asymp} is immediate from Lemma~\ref{mrtns}(b)
(cf.\ Remark~\ref{farey}), where, for each subbeam, $\,\lk_\infty=1-1/\sqrt2$.
(See (iv) and (v) above.) This completes the proof.
\end{proof}
\begin{rem}\label{share}
According to the preceding three lines, asymptotically,
the three subbeams contribute the same number of $\,\px\hs$-rational points:
the share of each subbeam, divided by $\,m^2$, has the limit
$\,3\sqrt2\hs(\sqrt2+1)\hs/\pi^2\,\approx\,1.038\,$ as $\,m\to\infty$.
\end{rem}

\begin{rem}\label{every}
If $\,m\ge2\,$ is even, {\it every\/} compact
K\"ahler-Einstein manifold $\,(N,h)\,$ of the odd complex dimension $\,m-1\,$
appears as an ingredient of the construction of some quadruple $\,\mgmt\,$ of
the third family, except for one restriction: if $\,h\,$ is Ricci-flat, its
K\"ahler cohomology class must be a real multiple of an integral class. This
is immediate is one combines Theorem~\ref{ctbly}(ii) with the definition of
$\,\dv\,$ in \S\ref{stat}\ (second paragraph after (\ref{xii})), relations
(\ref{kpo}.ii) and $\,\dv=\,\sgn\,\kx\,$ in \S\ref{stat}, and the third
paragraph of Remark~\ref{speca}.
\end{rem}

\section{The fourth family: examples of a new type}\label{frth}
\setcounter{equation}{0}
Given an odd integer $\,m\ge3$, let $\,\smt$ be the set of those
$\,\px\hs$-rational points in the $\,\mathcal{T}\nh$-beam $\,\xat\,$ of the
$\,\,\mathsf{X}\,\,$ component of the moduli curve (\S\ref{stat}) which do not
lie in the $\,\mathcal{H}\nh$-beam $\,\xah$. Our {\it fourth family} of quadruples
$\,\mgmt\,$ with (\ref{zon}) or (\ref{zto}) is obtained from the construction of
\S\ref{stat}\ applied to points $\,(\th,\ha)\in\smt$. (See also Remark~\ref{xprat}.)

The fourth family exists only in the odd complex dimensions $\,m\ge9$, since,
according to Theorem~\ref{asymp}(f) (proved below), $\,\smt$ is empty for
$\,m=3,5,7$. Also,  for $\,\lkm$ and $\,\p(\lkm,m)\,$ as in \S\ref{ftiq} and Lemma~\ref{mrtns}, and with $\,|\hskip3pt|\,$ denoting cardinality,
\begin{equation}
|\hs\smt|\,=\,2\hs|\mathcal{P}_{\!m}|\,,\hskip16pt\mathrm{where}\hskip10pt
\mathcal{P}_{\!m}\,=\,\hs\p(\lkm,m)\,.\label{stm}
\end{equation}
(As $\,\lkm>1$, we thus have $\,|\hs\smt|\le m(m-1)$.) In fact, we can define
a two-to-one surjective mapping
$\,\smt\ni(\th,\ha)\mapsto(k,\el\hs)\in\mathcal{P}_{\!m}$ by
assigning to $\,(\th,\ha)\,$ the pair $\,(k,\el\hs)\,$ of relatively prime
positive integers such that $\,(k-\el\hs)/\el\,$ is the value at
$\,(\th,\ha)\,$ of the function $\,\px\,$ described in \S\ref{trfp}. (This,
including the fact that $\,(k,\el\hs)\in\p(\lkm,m)$, is clear from (a) in
\S\ref{valu} and Definition~\ref{pratl}.)
\begin{proof}[Proof of Theorem~\ref{asymp}, parts
$\mathrm{(c),\hskip2pt(f)}$]First, (\ref{stm}), (e) in \S\ref{ftiq} and
Lemma~\ref{mrtns}(b) give (c).

Next, let us recall that $\,\mathcal{P}_{\!m}$ in (\ref{stm}) is the set of all
relatively prime integer pairs $\,(k,\el\hs)\,$ with
$\,1\le k\lkm<\el\le m\hs$. Since that gives $\,\lkm\le k\lkm<m\hs$, no such
pair exists if $\,\lkm>m$. Thus, (c) in \S\ref{ftiq} and (\ref{stm}) show that
$\,\smt\hskip-1.5pt=$\hskip5pt\emp\hskip7ptfor $\,m=3,5,7$. However, if $\,\lkm<m\hs$,
the set $\,\mathcal{P}_{\!m}$ is nonempty, as $\,(1,m)\in\mathcal{P}_{\!m}$.
Therefore, $\,|\hs\smt|=2\hs|\mathcal{P}_{\!m}|\ge2\,$ for any odd $\,m\ge9$, in
view of (\ref{stm}) and (b) in \S\ref{ftiq}. If $\,m\,$ is odd and $\,m\ge19$, we
have $\,\lkm<18\,$ (see (a) in \S\ref{ftiq}), so that $\,\mathcal{P}_{\!m}$ contains the
$\,(m-17)$-element subset of all $\,(k,\el\hs)\,$ with $\,k=1\,$ and
$\,18\le d\le m\hs$, which yields $\,|\hs\smt|\ge2(m-17)\,$ by (\ref{stm}).

Finally, $\,8<\lk_9<9\,$ by (c) in \S\ref{ftiq}. Thus, $\,(k,\el\hs)=(1,9)\,$ is the
only integer pair for which $\,1\le k\lk_9<\el\le9$, with $\,m=9\,$ (since
that gives $\,8k<k\lk_9<9$, and so $\,k=1$). From (\ref{stm}) we now get
$\,|\hs\mathcal{S}_9^{\mathcal{T}}|=2\hs|\mathcal{P}_{\nh9}|=2$, completing the
proof.
\end{proof}
Note that, as $\,9<\lk_{11}<10\,$ (cf.\ (c) in \S\ref{ftiq}), the same argument as in
the last three lines gives $\,|\hs\mathcal{S}_{11}^{\mathcal{T}}|=4$, since
$\,\mathcal{P}_{\nh11}$ has just two elements: $\,(1,10)\,$ amd $\,(1,11)$.
\begin{rem}\label{xprat}
By Lemma~\ref{xxast}(b), $\,(\xat)\cap(\xah)=\{(\xj,\xj^*)\}$, so that
$\,\smt$, defined at the
beginning of this section, is also the set of all $\,\px\hs$-rational points
in $\,\xat\,$ {\it other than} $\,(\xj,\xj^*)$.

For some odd integers $\,m\ge3\,$ the phrase `other than $\,(\xj,\xj^*)$'
used here is redundant, since $\,(\xj,\xj^*)\,$ is not $\,\px\hs$-rational.
Actually, we do not know if $\,(\xj,\xj^*)\,$ can be $\,\px\hs$-rational {\it
for any\/} odd $\,m\ge3$.

However, if $\,m\ge3\,$ is odd and $\,2\hskip.3pt\lkm>m\hs$, then
$\,\px\hs$-rationality of $\,(\xj,\xj^*)\,$ implies that $\,\lkm$ is an
{\it integer} and $\,\lkm\le m\hs$. In fact, by (a) in \S\ref{valu},
$\,\qx_*=-\hs1+1/\lkm$ is the value of $\,\px\,$ at $\,(\xj,\xj^*)$, and, as
$\,\qx_*\in(-\hs1,0)\,$ (Remark~\ref{xtoze}), the number $\,\qx_*$, now
assumed rational, must have the form
$\,(k-\el\hs)/\el\,$ for some $\,k,\el\in\bbZ\,$ with $\,1\le k<\el\le m\,$
(cf.\ Definition~\ref{pratl} and (a) in \S\ref{valu}); hence
$\,0<k\lkm=\el\le m$, so that $\,k=1\,$ (as $\,k\lkm\le m<2\hskip.3pt\lkm$)
and $\,\lkm=\el\in\bbZ$.

For instance, $\,(\xj,\xj^*)\,$ is {\it not\/} $\,\px\hs$-rational for any
$\,m\in\{3,5,7,9,11\}$. Namely, (c) in \S\ref{ftiq} then gives
$\,2\hskip.3pt\lkm>m\,$ (also for $\,m=13,15,17$, as
$\,2\hskip.3pt\lkm>2\hskip.3pt\lk_{11}>18>m\,$ by (e) in \S\ref{ftiq}). However,
again by (c) in \S\ref{ftiq}, one of the two conditions just named, necessary for
$\,\px\hs$-rationality of $\,(\xj,\xj^*)$, fails: $\,\lkm>m\,$ if
$\,m\in\{3,5,7\}$, and $\,\lk_9,\lk_{11}\notin\bbZ$. (One can extend this argument and conclusion to $\,m\in\{13,15,\dots,23\}$,
since a numerical approximation of $\,\lkm$ then gives $\,2\hskip.3pt\lkm>m\,$
and $\,\lkm\notin\bbZ$.)
\end{rem}

\section{Appendix: a version of Mertens's theorem}\label{mert}
\setcounter{equation}{0}
Assertion (a) in Lemma~\ref{mrtns} will now be derived using a slightly modified
version of a standard proof of Mertens's theorem (cf.\ \cite[p.~59]{chandrasekharan}). Here $\,m\ge1\,$ is treated as a real variable, even though in Lemma~\ref{mrtns}(a) it
stands for an odd integer with $\,m\ge3$.

Given a real number $\,\lk\ge1$, we define $\,W_\lk(m)\,$ to be the set
of all pairs $\,(k,\el\hs)\,$ of integers with $\,1\le k\lk<\el\le m$. If,
in addition, $\,n\,$ is an integer with $\,1\le n\le m$, let
$\,W_\lk(m,n)\subset W_\lk(m)\,$ consist of those $\,(k,\el\hs)\in W_\lk(m)\,$
in which both $\,k\,$ and $\,\el\,$ are divisible by $\,n$. Obviously,
for positive integers $\,n_1,\dots,n_s$ that are pairwise relatively prime,
and with $\,|\hskip3pt|\,$ denoting cardinality,
\begin{equation}
\mathrm{a)}\hskip7.5ptW_\lk(m,\prod_{j=1}^s\nh n_j)\,
=\,{\displaystyle\bigcap_{j=1}^sW_\lk(m,n_j)}\hs,
\hskip17pt\mathrm{b)}\hskip7.5pt|W_\lk(m,n)|\,=\,|W_\lk(m/n)|\hs,\label{wlm}
\end{equation}
where (b) is due to the bijection $\,W_\lk(m,n)\to W_\lk(m/n)\,$ given by
$\,(k,\el\hs)\mapsto(k/n,\el/n)$. Next, given $\,\lk,m\in[\hs1,\infty)$, we
let $\,\p(\lk,m)\,$ denote, as in Lemma~\ref{mrtns}, the subset of $\,W_\lk(m)\,$
formed by those integer pairs $\,(k,\el\hs)\,$ which, in addition to having
$\,1\le k\lk<\el\le m$, are also {\it relatively prime}.

For any finite family $\,\mathcal{A}\,$ of finite sets, induction on its
cardinality $\,|\mathcal{A}\hs|\,$ easily implies that
$\,|\nh\bigcup\mathcal{A}\hs|
=\sum_{\mathcal{B}}(-1)^{|\mathcal{B}|-1}\,|\hskip-1pt\bigcap\mathcal{B}\hs|$, with summation
over all nonempty subfamilies $\,\mathcal{B}\,$ of $\,\mathcal{A}$. If $\,\mathcal{A}\,$
consists of all sets $\,W_\lk(m,p)$, where $\,m\in[\hs1,\infty)\,$ is fixed
and $\,p\,$ runs through all primes with $\,p\le m\hs$, then
$\,\bigcup\mathcal{A}\,$ clearly coincides with
$\,W_\lk(m)\smallsetminus \p(\lk,m)$, so that (\ref{wlm}) and our formula for
$\,|\nh\bigcup\mathcal{A}\hs|\,$ give
\begin{equation}
|\p(\lk,m)|\hskip12pt=\,\sum_{1\hs\le\hs j\hs\le\hs m}\hs\mu_j\,|W_\lk(m/j)|\,.
\label{plm}
\end{equation}
Here $\,\mu\,$ is the {\it M\"obius function} \cite{chandrasekharan}), assigning to every
integer $\,j\ge1\,$ the value $\,\mu_j\,$ with $\,\mu_j=0\,$ if $\,j\,$ is
divisible by the square of a prime and $\,\mu_j=(-1)^k$ when $\,j\,$ is the
product of $\,k\,$ distinct primes, for $\,k\ge0$. (Thus, $\,\mu_1=1$.) The
summation index $\,j$, here and below, is an integer.

Setting $\,f(m)=2\lk\hs|W_\lk(m)|-m^2$, we have $\,|f(m)|\le\alpha\hs m\,$ for
every $\,m\in[\hs1,\infty)$, with the constant $\,\alpha=2\lk+2$. In fact,
given $\,m\in[\hs1,\infty)$, let $\,n\,$ (or, $\,r$) be the largest integer
not exceeding $\,m\,$ (or, $\,m/\lk$) and, for any $\,k\in\{1,\dots,r\}$, let
$\,s_k$ be the largest integer with $\,s_k\le k\lk$. Thus, $\,W_\lk(m)\,$
contains exactly $\,n-s_k$ points with the first coordinate
$\,k$, namely, $\,(k,\el\hs)\,$ with $\,s_k<\el\le n$. Consequently,
$\,|W_\lk(m)|=\sum_{k=1}^r(n-s_k)$, and so the inequalities
$\,k\lk-1<s_k\le k\lk$, that is, $\,n-k\lk\le n-s_k<n+1-k\lk$, yield
$\,[2n-(r+1)\lk]\hs r\le
2|W_\lk(m)|\le[2(n+1)-(r+1)\lk]\hs r$. Hence
$\,-\hs(2\lk+2)m\le f(m)\le2m$. Namely,
$\,-\hs m-\lk\le-\hs(r+1)\lk\le-\hs m\,$ and $\,m-1\le n\le m\,$ due to our
choice of $\,m\,$ and $\,r$, so that $\,2n-(r+1)\lk\ge m-\lk-2\,$ and
$\,2(n+1)-(r+1)\lk\le m+2$, while $\,-\hs1+m/\lk\le r\le m/\lk$. Thus,
$\,-\hs\alpha\hs m\le f(m)\le\alpha\hs m$.

Given a bounded sequence $\,\mu_1,\mu_2,\dots\,$ of real numbers and a
function $\,f\,$ of the real variable $\,m\ge1\,$ such that
$\,|f(m)|\le\alpha\hs m\,$ for all $\,m\,$ and some constant $\,\alpha\ge0$,
we necessarily have $\,m^{-2}\sum_{1\le j\le m}\mu_j\hs f(m/j)\,\to\,0\,$ as
$\,m\to\infty$, which is obvious since
$\,\sum_{1\le j\le m}j^{-1}\le1+\log m\,$ due to a standard
area-under-the-graph estimate. For $\,f\,$ and $\,\alpha\,$ as in the last
paragraph, with the M\"obius function $\,\mu$, this shows that, by (\ref{plm}),
if $\,\sum_{1\le j\le m}(\mu_j/j^2)\,$ has a limit as $\,m\to\infty$, then so
does $\,2\lk\hs|\p(\lk,m)|/m^2$, and the limits coincide. However,
$\,\sum_{1\le j\le m}(\mu_j/j^2)\,$ clearly does converge, as $\,m\to\infty$,
to the product $\,\prod_p(1-1/p^2)$, where $\,p\,$ runs over all primes; now
(a) in Lemma~\ref{mrtns} follows since the inverse of the product equals
$\,\sum_{\hs n=1}^{\hs\infty}n^{-2}\hskip-1.2pt=\pi^2/6$ due to a special case
of {\it Euler's identity}, obtained by expanding each factor
$\,1/(1-1/p^2)\,$ into a geometric series. (Cf.\ \cite[pp.~61 and~76]{chandrasekharan}.)

\section{Appendix: decimal approximations}\label{dcax}
\setcounter{equation}{0}
We show here that, for $\,\hi_m$ defined in Lemma~\ref{bonds},
\begin{equation}
\hi_7>1.9,\hskip12pt2.2<\hi_9<2.3<2.4<\hi_{11}<2.5,
\hskip12pt\hi_{13}<2.6,\hskip12pt\hi_{15}<2.7,\hskip12pt
\hi_{17}<2.8\hs.\label{chi}
\end{equation}
This is verified using the following algorithm, designed for calculations that
can even be done by hand, and produce upper/lower bounds on $\,\hi_m$ having
the form $\,\hi_m<\sa_m$ or $\,\hi_m>\sa_m$ for a fixed odd integer
$\,m\ge3\,$ and a rational number $\,\sa_m>0\,$ with a simple decimal
expansion.

The values of $\,\vs_j(0)\,$ for $\,j=1,\dots,m-1\,$ can easily be found
from (\ref{smo}.i); when $\,m=17$, they are 1, 1, 2, 5, 14, 42, 132, 429,
1430, 4862, 16796, 58786, 208012, 742900, 2674440, 9694845. In formula
(\ref{smo}.iv), with $\,m\,$ replaced by
$\,m-1$, we now estimate each term $\,(-1)^j\vs_{m-j-1}(0)\hs\sa^j$ for
$\,\sa=\sa_m$ from above/below by the nearest integer. When these estimates
are added up, (\ref{smo}.iv) yields an up\-per/lower bound on
$\,\vs_{m-1}(0)R_{m-1}(\ps_m)$, for $\,\ps_m$ related to $\,\sa_m$ as in
(\ref{tts}).

Whenever this last bound happens to ensure negativ\-ity/pos\-i\-tiv\-ity of
$\hs R_{m-1}(\ps_m)\nh+\nh K_m(\sa_m)$, Remark~\ref{kmsig}(b) will imply
negativ\-ity/pos\-i\-tiv\-ity of $\,\hi_m-\sa_m$, as required.

Specifically, let $\,m\in\{7,9,11\}$. If $\,\sa_7=1.9$, $\,\sa_9=2.2$,
$\,\sa_{11}=2.4$, these steps give
$\,42\hs R_6(\ps_7)>42-27+18-14+13-25-24=-\hs17$,
$\quad429\hs R_8(\ps_9)>429-291+203-150+117-104+113-250-275=-\hs208$,
$\quad4862\hs R_{10}(\ps_{11})
>4862-3432+2471-1825+1393-1115+955-918+1100-2642-3171=-\hs2322$, so that
$\,R_{m-1}(\ps_m)>-\hs1/2\,$ for $\,m=7,9,11$, while
$\,K_7(\sa_7)>11/19>1/2\,$
as $\,\sa_7<2$, $\,K_9(\sa_9)>5/9>1/2\,$ as $\,\sa_9<2.5$, and
$\,K_{11}(\sa_{11})=19/33>1/2$. Hence $\,R_{m-1}(\ps_m)+K_m(\sa_m)>0\,$ and
Remark~\ref{kmsig}(b) yields the lower bounds on $\,\hi_7,\hi_9,\hi_{11}$ in
(\ref{chi}).

Next, for $\,m=9,11,13,15,17$, let us set, this time, $\,\sa_9=2.3$,
$\,\sa_{11}=2.5$, $\,\sa_{13}=2.6$, $\,\sa_{15}=2.7$, $\,\sa_{17}=2.8$. Our
algorithm now yields
$\,429\hs R_8(\ps_9)<429-303+223-170+140-128+149-340-391=-\hs391$,
$\quad4862\hs R_{10}(\ps_{11})
<4862-3575+2682-2062+1641-1367+1221-1220+1526-3814-4768=-\hs4874$,
$\quad58786\hs R_{12}(\ps_{13})
<58786-43669+32868-25133+19605-15683+12975-11244+10442-10859+14117-36703-47714
=-42212$,
$\quad742900\hs R_{14}(\ps_{15})
<742900-561632+428550-330595+258387-205189+166204-138076+118621-106758+102946
-111181+150095-405255-547094=-\hs438077$,
$\quad9694845\hs R_{16}(\ps_{17})
<9694845-7488432+5824336-4566279+3613317-2890653+2342951-1929488+1620771
-1396356+1244027-1161091+1161092-1300422+1820592-5097655-7136717=-\hs5645162$,
so that $\,R_8(\ps_9)<-\hs0.9$,
$\,R_{10}(\ps_{11})<-\hs1\,$ and $\,R_{m-1}(\ps_m)<-\hs0.58\,$ for
$\,m=13,15,17$, while
$\,K_9(\sa_9)<3/5=0.6\,$ as $\,\sa_9>2$,
$\,K_{11}(\sa_{11})=38/67<1$, and
$\,K_m(\sa_m)<0.57\,$ for $\,m=13,15,17$. Thus,
$\,R_{m-1}(\ps_m)+K_m(\sa_m)<0\,$ for $\,m=9,11,13,15,17$, and
Remark~\ref{kmsig}(b) gives the upper bounds on $\,\hi_m$ required in
(\ref{chi}).

\section{Appendix: vertical compactness}\label{vert}
\setcounter{equation}{0}
The following discussion provides a differential-geometric background for the
moduli curve, and is not used to derive any results of this paper. To save
space, our presentation is brief.

Let a quadruple $\,\mgmt\,$ have all the properties listed in (\ref{zon}) or (\ref{zto}) except for compactness of $\,M$. The set $\,M'\subset M\,$ on which
$\,d\vp\ne0\,$ then is open and dense in $\,M\,$ (cf.\
\cite[Remark~5.4]{local}, while the complex vector subbundle of $\,TM'$ spanned by the
$\,g$-gradient $\,\navp\,$ is an integrable real $\,2$\diml\ distribution on
$\,M'$ with totally geodesic leaves \cite[the end of \S7]{local}. We will say
that the quadruple $\,\mgmt\,$ satisfies the {\it vertically compact version}
of (\ref{zon}) or (\ref{zto})
if every such leaf is contained in a compact submanifold of real
dimension $\,2\,$ in $\,M$.

The constructions of \S\ref{main} and \S\ref{redu}\ applied to data satisfying all the
assumptions stated there except for compactness of $\,N\hs$ always leads to
$\,\mgmt\,$ satisfying the vertically compact version of (\ref{zon}) or (\ref{zto}). (In fact, the compactness assumption is never used in either
construction.) Similarly, the vertically compact version of (\ref{zon}) or (\ref{zto}) holds for $\,\mgmt\,$ obtained as in \S\ref{stat}\ from $\,(\th,\ha)\,$ and
additional data having all the properties listed in \S\ref{stat}\ except that
$\,(\th,\ha)\,$ is assumed merely to lie in the moduli curve $\,\mathcal{C}\,$ (and
is not required to be $\,\px\hs$-rational), while $\,N\hs$ is not necessarily
compact.

Conversely, every quadruple $\,\mgmt\,$ satisfying the vertically compact
version of (\ref{zon}) or (\ref{zto}) is obtained in this way from some
$\,(\th,\ha)\in\mathcal{C}\hs$ and additional data mentioned above. In fact,
Theorems~33.2,~33.3 and~34.3 in \cite{potentials} (which involve the four
types mentioned in (i) - (iii) of \S\ref{more} of this paper), as well as
Lemma~\ref{asttt} in \S\ref{more}, all remain valid, with essentially the same
proofs, also in the vertically compact case, and so our claim follows as in
the last three lines of \S\ref{more}.

\begin{acknowledgements}
We wish to thank David Farber and Joshua Lansky for comments and information
concerning Lemma~\ref{mrtns} and \S\ref{mert}.
\end{acknowledgements}

\end{document}